\documentclass[11pt]{article}
\usepackage{amsmath, amssymb, amscd, xypic}
\usepackage[dvips]{graphics}

\newtheorem{thm}{Theorem}
\newtheorem{prop}[thm]{Proposition}
\newtheorem{lem}[thm]{Lemma}

\newtheorem{cor}[thm]{Corollary}
\newtheorem{rem}[thm]{Remark}
\newtheorem{df}[thm]{Definition}

\renewcommand{\epsilon}{\varepsilon}
\renewcommand{\phi}{\varphi}

\newcommand{\BB}{\mathbb}

\newcommand{\pf}{\noindent {\it Proof. }}
\newcommand{\qed}{\nopagebreak $\qquad$ $\square$ \vskip5pt}
\newcommand{\separate}{\vskip5pt}

\newcommand{\im}{\operatorname{Im}}
\newcommand{\re}{\operatorname{Re}}
\newcommand{\tr}{\operatorname{Tr}}

\newcommand{\B}{\overline}
\newcommand{\HC}{\BB H_{\BB C}}
\newcommand{\HR}{\BB H_{\BB R}}

\newcommand{\pol}{\operatorname{Pol}}
\newcommand{\C}{\operatorname{C}}

\usepackage[OT2,T1]{fontenc}
\newcommand\textcyr[1]{{\fontencoding{OT2}\fontfamily{wncyr}\selectfont #1}}
\newcommand{\Zh}{\textit{\textcyr{Zh}}}

\newcommand{\M}{\operatorname{Mx}}
\newcommand{\xm}{\operatorname{Xm}}
\newcommand{\grad}{\overrightarrow{\nabla}}
\newcommand{\ha}{\operatorname{H}^2}

\begin{document}

\title{\bf Quaternionic Analysis, Representation Theory and Physics}
\author{Igor Frenkel and Matvei Libine}
\maketitle

\begin{abstract}
We develop quaternionic analysis using as a guiding principle
representation theory of various real forms of the conformal group.
We first review the Cauchy-Fueter and Poisson formulas and explain
their representation theoretic meaning.
The requirement of unitarity of representations leads us to the extensions
of these formulas in the Minkowski space, which can be viewed as
another real form of quaternions.
Representation theory also suggests a quaternionic version of the
Cauchy formula for the second order pole.
Remarkably, the derivative appearing in the complex case is replaced
by the Maxwell equations in the quaternionic counterpart.
We also uncover the connection between quaternionic analysis and various
structures in quantum mechanics and quantum field theory, such as the
spectrum of the hydrogen atom, polarization of vacuum, one-loop Feynman
integrals. We also make some further conjectures.
The main goal of this and our subsequent paper is to revive quaternionic
analysis and to show profound relations between quaternionic analysis,
representation theory and four-dimensional physics.
\end{abstract}


\section{Introduction}

It is well known that after discovering the algebra of quaternions
$\BB H = \BB R 1 \oplus \BB R i \oplus \BB R j \oplus \BB R k$
and carving the defining relations on a stone of Dublin's Brougham Bridge
on the 16 October 1843, the Irish physicist and mathematician
William Rowan Hamilton (1805-1865) devoted the remaining years
of his life developing the new theory which he
believed would have profound applications in physics. But one had to wait
another 90 years before von Rudolf Fueter produced a key result of
quaternionic analysis, an exact quaternionic counterpart of the Cauchy
integral formula
\begin{equation}  \label{Cauchy-intro}
f(w) = \frac 1{2\pi i} \oint \frac {f(z)\,dz}{z-w}.
\end{equation}
Because of the noncommutativity of quaternions, this formula comes in two
versions, one for each analogue of the complex holomorphic functions --
left- and right-regular quaternionic functions:
\begin{align}
f(W) &= \frac 1 {2\pi^2} \int_{\partial U}
\frac{(Z-W)^{-1}}{\det(Z-W)} \cdot *dZ \cdot f(Z), \label{Fueter-intro1}\\
g(W) &= \frac 1 {2\pi^2} \int_{\partial U}
g(Z) \cdot *dZ \cdot \frac {(Z-W)^{-1}}{\det(Z-W)},
\qquad \forall W \in U,  \label{Fueter-intro2}
\end{align}
where $U \subset \BB H$ is a bounded open set,
the determinant is taken in the standard matrix realization of $\BB H$,
$f,g: \BB H \to \BB H$ are differentiable functions satisfying the analogues
of the Cauchy-Riemann equations
\begin{align}
\nabla^+f &= \frac{\partial f}{\partial z^0} +
i\frac{\partial f}{\partial z^1} + j\frac{\partial f}{\partial z^2} +
k\frac{\partial f}{\partial z^3} =0,  \label{Dirac1-intro}  \\
g\nabla^+ &= \frac{\partial g}{\partial z^0} +
\frac{\partial g}{\partial z^1}i + \frac{\partial g}{\partial z^2}j +
\frac{\partial g}{\partial z^3}k =0,  \label{Dirac2-intro}
\end{align}
and $*dZ$ is a quaternionic-valued differential three form that is
Hodge dual to $dZ = dz^0 + idz^1 + jdz^2 + kdz^3$.

The Cauchy formula (\ref{Cauchy-intro}) also provides extensions of functions
on the unit circle to holomorphic functions on the discs $D_2^{\pm}$ inside
and outside the unit circle in $\BB C \cup \{\infty\}$.
Conversely, one can recover the original function on the circle from
the boundary values of the two holomorphic functions on $D_2^{\pm}$
and get a polarization
$$
L^2(S^1) \simeq \ha(D_2^+) \oplus \ha(D_2^-),
$$
where $\ha(D_2^{\pm})$ denote the Hardy spaces of analytic functions
on $D_2^{\pm}$. The dense subspace of trigonometric polynomials is identified
under the above isomorphism with the space of Laurent polynomials
on $\BB C^{\times} = \BB C \setminus \{0\}$ with the natural polarization.
Similarly, the Cauchy-Fueter formulas (\ref{Fueter-intro1}),
(\ref{Fueter-intro2}) yield decompositions
\begin{equation}  \label{decomposition-intro}
L^2 (S^3, \BB H) \simeq V(D_4^+) \oplus V(D_4^-), \qquad
L^2 (S^3, \BB H) \simeq V'(D_4^+) \oplus V'(D_4^-),
\end{equation}
where $D_4^{\pm}$ denote the balls inside and outside the unit sphere $S^3$
in $\BB H \cup \{\infty\}$, and $V(D_4^{\pm})$, $V'(D_4^{\pm})$ are the spaces
of left-, right-regular functions on $D_4^{\pm}$ with square integrable norm
induced from $L^2 (S^3, \BB H)$. This norm can also be described intrinsically,
as in the complex case.
Again, the dense subset of $K$-finite (with $K = SU(2) \times SU(2)$)
$\BB H$-valued functions on $S^3$ is isomorphic to polynomials on
$\BB H^{\times}$ with the polarization induced by the decompositions
(\ref{decomposition-intro}).

There is also a quaternionic counterpart of the Poisson formula
for harmonic functions on $\BB H$:
\begin{equation}  \label{Poisson-intro}
\phi(W) = \frac 1{2\pi^2} \int_{Z \in S^3}
\frac {(\widetilde{\deg}\phi)(Z)}{\det(Z-W)} \,dS
= \frac 1{2\pi^2} \int_{Z \in S^3}
\frac {1 - \det(W)}{\det(Z-W)^2} \cdot \phi(Z) \,dS,
\end{equation}
where
$\widetilde{\deg}=1 + z^0 \frac{\partial}{\partial z^0} +
z^1 \frac{\partial}{\partial z^1} + z^2 \frac{\partial}{\partial z^2}
+ z^3 \frac{\partial}{\partial z^3}$ and
$dS$ is the Euclidean measure on the three-dimensional sphere $S^3$
of radius $1$.
The spaces of harmonic and left- (or right-) regular functions are invariant
under the conformal (fractional linear) action of the group $SL(2,\BB H)$,
as defined in Proposition \ref{action-harmonic}.
This implies that the dense spaces of polynomial functions in the summands of
(\ref{decomposition-intro}) are invariant under the action of the
Lie algebra $\mathfrak{sl}(2,\BB H)$.
This is in complete parallel with the global conformal invariance of
meromorphic functions under $SL(2,\BB C)$ in the complex case.
Some contemporary reviews of the quaternionic analysis are given in
\cite{Su} and \cite{CSSS}.
Unfortunately, this promising parallel between complex and quaternionic
analysis essentially ends here, and, in spite of vigorous attempts,
the quaternionic analogues of the ring structure of holomorphic functions,
local conformal invariance, Riemann mapping theorem and many other classical
results of complex analysis have never been found.
Such a failure has even led R.~Penrose to say,
``[Quaternions] do have some very significant roles to play,
and in a slightly indirect sense their influence has been enormous,
through various types of generalizations. But the original `pure quaternions'
still have not lived up to what must undoubtedly have initially seemed to be
an extraordinary promise... The reason appears to be that there is no
satisfactory quaternionic analogue of the notion of a holomorphic function.''

In this paper we propose to approach the quaternionic analysis from the point
of view of representation theory of the conformal group $SL(2,\BB H)$ or,
better, the corresponding Lie algebra $\mathfrak{sl}(2,\BB H)$ and related
real forms.
Let us denote by ${\cal V}^+$ (respectively ${\cal V}'^+$) the spaces
of left- (respectively right-) regular functions defined on $\BB H$ with values
in the two-dimensional complex space $\BB S$ (respectively $\BB S'$)
obtained from $\BB H$ by fixing a complex structure.
We also denote by ${\cal H}^+$ the space of complex-valued harmonic functions
on $\BB H$.
Then ${\cal H}^+$, ${\cal V}^+$ and ${\cal V}'^+$ turn out to be irreducible
representations of $\mathfrak{sl}(2,\BB H)$ and its complexification
$\mathfrak{sl}(4,\BB C)$.
If we want to introduce unitary structures on these spaces, we must replace
$\mathfrak{sl}(2,\BB H)$ with another real form of $\mathfrak{sl}(4,\BB C)$,
namely $\mathfrak{su}(2,2)$. Then we are able to identify
${\cal H}^+$, ${\cal V}^+$ and ${\cal V}'^+$ with the irreducible unitary
representations of $\mathfrak{su}(2,2)$ of spin $0$, $1/2$ and $1/2$
respectively of the most degenerate series of representations with highest
weights. In order to develop further the quaternionic analysis following the
parallel with complex analysis, we should look for results that can be restated
in terms of representation theory of the complex conformal group $SL(2,\BB C)$
and its Lie algebra $\mathfrak{sl}(2,\BB C)$.
In particular, the Cauchy integral formula for the second order pole
\begin{equation}  \label{2pole-intro}
f'(w) = \frac 1{2\pi i} \oint \frac {f(z)\,dz}{(z-w)^2}
\end{equation}
can be viewed as an intertwining operator between certain representations of
$SL(2,\BB C)$. For illustration purposes, let us consider three actions of
$SL(2,\BB C)$ on meromorphic functions on $\BB CP^1$ with parameter $k=0,1,2$:
\begin{equation}  \label{rho-intro}
\rho_k(h): \: f(z) \mapsto (\rho_k(h)f)(z) =
\frac 1{(cz+d)^k} \cdot f \biggl( \frac {az+b}{cz+d} \biggr), \qquad
h^{-1} = \begin{pmatrix} a & b \\ c & d \end{pmatrix} \in SL(2, \BB C).
\end{equation}
Differentiating, we obtain three actions of $\mathfrak{sl}(2, \BB C)$
which preserve the space of polynomials on $\BB C$.
Let $V_0$, $V_1$ and $V_2$ denote this space of polynomials with
$\mathfrak{sl}(2, \BB C)$-actions $\rho_0$, $\rho_1$ and
$\rho_2$ respectively.
Then we have an $\mathfrak{sl}(2, \BB C)$-equivariant map
\begin{equation}  \label{M-intro}
M: (\rho_1 \otimes \rho_1, V_1 \otimes V_1) \to
(\rho_2, V_2), \qquad
f_1(z_1) \otimes f_2(z_2) \mapsto f_1(z) \cdot f_2(z).
\end{equation}
On the other hand, the derivation map is an intertwining operator between
$V_0$ and $V_2$:
\begin{equation}  \label{derivative-intro}
\frac d{dz}: (\rho_0, V_0) \to (\rho_2, V_2),
\qquad f(z) \mapsto \frac d{dz} f(z).
\end{equation}
This map has a kernel consisting of the constant functions, and its image
is all of $V_2$ which is irreducible.

Action $\rho_1$ of $SL(2,\BB C)$ is similar to the action of
$SL(2,\BB H)$
on the spaces of functions ${\cal H}^+$, ${\cal V}^+$ and ${\cal V}'^+$.
The quaternionic analogues of the multiplication map (\ref{M-intro}) lead
to the $\mathfrak{sl}(2,\BB H)$-equivariant maps
$$
{\cal H}^+ \otimes {\cal H}^+ \twoheadrightarrow \Zh^+, \qquad
{\cal V}^+ \otimes {\cal V}'^+ \twoheadrightarrow {\cal W}^+.
$$
It turns out that the representations $\Zh^+$ and ${\cal W}^+$ are irreducible,
unitary with respect to $\mathfrak{su}(2,2)$ and belong to another degenerate
series of irreducible unitary representations of $\mathfrak{su}(2,2)$
-- the so-called middle series. The middle series is formed from
the lowest component in the tensor product of dual
representations from the most degenerate series -- in our case of spin $0$
and spin $1/2$ -- so that $\Zh^+ \subset {\cal H}^+ \otimes {\cal H}^+$
and ${\cal W}^+ \subset {\cal V}^+ \otimes {\cal V}'^+$.
The quaternionic counterparts of the intertwining operator
(\ref{derivative-intro}) and its integral presentation (\ref{2pole-intro}) are
\begin{equation}  \label{Mx-intro}
\M f = \nabla f \nabla - \square f^+
\end{equation}
($f^+$ denotes the quaternionic conjugate of $f$)
and (note the square of the Fueter kernel)
\begin{equation}  \label{2pole1-intro}
(\M f)(W) = \frac {3i}{\pi^3} \int_{C_4} \frac {(Z-W)^{-1}}{\det(Z-W)} \cdot
f(Z) \cdot \frac {(Z-W)^{-1}}{\det(Z-W)} \,dZ^4,
\end{equation}
where $dZ^4$ is the volume form, $C_4$ is a certain four cycle in the
space of complexified quaternions
$\HC = \BB C 1 \oplus \BB C i \oplus \BB C j \oplus \BB C k$,
and $f: \HC \to \HC$ is a holomorphic function of four complex variables.
The image of $\M$ is ${\cal W}^+$.
Comparing (\ref{derivative-intro}) and (\ref{Mx-intro}) we can see that
the quaternionic analogue of the constant functions in the complex case
is the kernel of the operator $\M$, which turns out to be a Euclidean
version of the Maxwell equations for the gauge potential.
There is also a similar formula that involves the square of the Poisson
kernel for the scalar-valued functions that comprise $\Zh^+$.

Next we study a natural problem of describing the projectors
\begin{align*}
{\cal H}^+ \otimes {\cal H}^+ \twoheadrightarrow & \Zh^+ \hookrightarrow
{\cal H}^+ \otimes {\cal H}^+,  \\
{\cal V}^+ \otimes {\cal V}'^+ \twoheadrightarrow & {\cal W}^+ \hookrightarrow
{\cal V}^+ \otimes {\cal V}'^+.
\end{align*}
Our explicit description of these projectors as integral operators can be
viewed as the double pole formula in quaternionic analysis, presented in a
separated form for harmonic or regular functions of two quaternionic variables.
For example, for harmonic functions the projector is given by integrating
against a kernel
\begin{equation}  \label{2pole2-intro}
k(Z_1,Z_2;W_1,W_2) = \int_{C_4}
\frac {dT^4}{\det(Z_1-T)\cdot\det(Z_2-T)\cdot\det(W_1-T)\cdot\det(W_2-T)},
\end{equation}
where the integration again is done over a certain four cycle in $\HC$.
This is a very interesting function, it can be computed using the Cayley
transform that maps the cycle of integration $C_4$ into the flat Minkowski
space. The result turns out to be related to the hyperbolic volume of
an ideal tetrahedron, and is given by the dilogarithm function.
The appearance of these structures indicates a potential richness of
quaternionic analysis.

The representation theoretic approach explains why there cannot be a
natural ring structure in the space of left- or right-regular functions.
In fact, the latter spaces provide representations of $\mathfrak{sl}(2,\BB H)$
from the most degenerate series, while their tensor products belong
to the middle series that consists of functions depending on four rather
than three variables. In the complex case the representation $V_0$ has a
natural algebra structure that can be viewed as the algebra of holomorphic
functions. The quaternionic counterpart of $V_0$ is the dual space
${\cal W}'^+$, which consists of all quaternionic-valued functions on $\HC$.
As a representation of $\mathfrak{sl}(2,\BB H)$, it still does not have an
algebra structure, but contains a subspace which acquires a natural algebra
structure if one adjoins a unity.

It might seem unexpected that a generalization of the Cauchy-Fueter formula
to the second order pole requires an extension of the quaternionic analysis
to the complexified quaternions $\HC$. But this is not the only reason for the
ubiquitous presence of the complexified quaternions in the real quaternionic
analysis. We have seen that in order to introduce unitary structures on the
spaces of harmonic, left- and right-regular functions, we must replace the
quaternionic conformal group $SL(2,\BB H)$ with another real form of
$SL(2,\HC) = SL(4,\BB C)$,
namely $SU(2,2)$, which in turn can be identified with the conformal group
of the Minkowski space $\BB M$ realized as a real form of $\HC$.
Also the Minkowski space naturally emerges when we apply quaternionic
analysis to the spectral decomposition of the Hamiltonian of the hydrogen
atom studied in \cite{BI} (also see references therein).
It turns out that the Poisson formula in $\BB H$ yields the discrete part
of the spectrum, while the analysis of the continuous part leads to
a counterpart of the Poisson formula in $\BB M$.
In fact, the spaces ${\cal V}$, ${\cal V}'$ and ${\cal H}$ consisting of
regular and harmonic functions on $\BB H^{\times} = \BB H \setminus \{0\}$
can also be realized respectively as the spaces of regular
functions and solutions of the wave equation on $\BB M$.
Now the polarization is given by the support of the
Fourier transform on the future and past light cone.
There is also a conformal map that transforms the unit sphere
$S^3 \subset \BB H$ into the unit two-sheeted hyperboloid $H^3 \subset \BB M$
which induces natural isomorphisms
\begin{equation}  \label{decomposition2-intro}
L^2(S^3, \BB S) \simeq L^2 (H^3,\BB S), \qquad
L^2(S^3, \BB S') \simeq L^2 (H^3,\BB S').
\end{equation}
Thus the integrals over three-dimensional spheres in $\BB H$ can be
replaced by the integrals over two-sheeted hyperboloids in $\BB M$.
It is interesting to note that integrals over just one sheet --
the Lobachevski space -- do not give the correct Poisson formula.
Another advantage of the Minkowski
space over the space of quaternions is that the four-dimensional integrals in
the double pole formulas (\ref{2pole1-intro}) and (\ref{2pole2-intro})
will not need extensions to the complex domain $\HC$ if one uses generalized
functions or, in other words, need only an infinitesimal extension.
However, the Minkowski space formulation also brings some technical
difficulties related to the fact that the singularities of the kernels in
(\ref{Fueter-intro1}), (\ref{Fueter-intro2}) and (\ref{Poisson-intro})
are now concentrated on the light cone instead of just a single point in
the initial quaternionic picture.
These difficulties are resolved by infinitesimal extensions to the
complex domain $\HC$ or, equivalently, using generalized functions.

The Minkowski space reformulation of quaternionic analysis brings us into
a thorough study of Minkowski space realization of unitary
representations of the conformal group $SU(2,2)$ by H.~P.~Jakobsen and
M.~Vergne \cite{JV1,JV2}
who were motivated by the program of I.~E.~Segal \cite{Se} on the
foundational role of representation theory of $SU(2,2)$ in physics.
(Also see references in \cite{JV1,JV2}.)
In our paper we use their results extensively.
Taking the point of view of quaternionic analysis, we extend, give alternative
proofs and make more explicit some of their results for the degenerate
series representations.
We also would like to note that quaternionic analysis in $\BB H$ clarifies
certain aspects such as the unitary structures and the $K$-types, not directly
visible in the Minkowski space formulation.

The Minkowski space formulation also suggests that there should be analogues
of the Fueter and Poisson formulas where integration takes place over
the one-sheeted hyperboloid, usually called the imaginary Lobachevski space.
This will lead us to another change of the real forms of the conformal group
$SL(4,\BB C)$ and the complex domain $\HC$ -- to $SL(4,\BB R)$ and $\HR$
-- the split real quaternions consisting of $2 \times 2$ real matrices.
This case will be considered in our upcoming paper \cite{FL}.

The appearance of various real forms $\BB H$, $\BB M$, $\HR$ of $\HC$ in the
different developments of quaternionic analysis give rise to a more general
point of view than we indicated in the beginning of the introduction.
Namely one should consider the holomorphic functions in the open domains
$U^{\BB C} \subset \HC$, which, in addition, satisfy the quaternionic
left-/right-regularity (\ref{Dirac1-intro})/(\ref{Dirac2-intro}).
The Fueter formula, as well as the Poisson and double pole formulas,
can then be stated in a greater generality involving integration
over more general cycles in $U^{\BB C} \subset \HC$.
On the representation theory side considering such functions allows
us to relate representations of different real forms of $SL(4,\BB C)$
and their subgroups, which might suggest similar results for more general
classes of groups.

The Fueter and Poisson formulas allow generalizations to higher dimensions
and play a central role in Clifford analysis \cite{DSS}.
However, we believe that many subsequent results -- such as the double pole
formula -- are unique to quaternionic analysis, and the theory can be developed
in its own special direction.
In particular, the quaternionic analysis is deeply related to the harmonic
analysis associated with the simplest simple complex Lie group $SL(2,\BB C)$
studied in the classical work \cite{GGV}.
In fact, thanks to the isomorphisms (\ref{decomposition-intro}) and
(\ref{decomposition2-intro}) one can realize the spaces of functions on
$SU(2)$ and $SL(2,\BB C)/SU(2)$ via harmonic and regular functions on
$\BB H$ and $\BB M$ respectively. In the next paper \cite{FL} we will have
a similar realization of the spaces of functions on $SL(2,\BB R)$ and
$SL(2,\BB C)/SL(2,\BB R)$.
The use of the wave equation for the study of harmonic analysis on
hyperboloid goes back to \cite{St}, for more recent work in this direction
see \cite{KO}.

Another unique feature of quaternionic analysis is its deep relation to
physics, in particular, to the four-dimensional classical and quantum field
theories for massless particles. Many applications of quaternionic analysis
to physics are collected in the book  \cite{GT}.
In our paper we would like to add to their list a few more important examples.
We have already mentioned the implication of the Poisson formula -- both for
$\BB H$ and $\BB M$ -- to the spectral decomposition of the Hamiltonian
of the hydrogen atom.
In general, the Minkowski formulation of various results of quaternionic
analysis provides a link to the four-dimensional field theories.
This is hardly surprising since the Minkowski space is the playground
for these physical theories, but it is still quite remarkable
that we encounter some of the most fundamental objects of these theories.
It is certainly clear that the equations for the left- and
right-regular functions (\ref{Dirac1-intro}) and (\ref{Dirac2-intro}) are
nothing but the massless Dirac equation. But it comes as a surprise that the
quaternionic analogue of the Cauchy formula for the second order pole
(\ref{2pole-intro}) is precisely the Maxwell equations for the gauge potential.
Moreover, the integral itself appears in the Feynman diagram for vacuum
polarization and is responsible for the electric charge renormalization.
Also, the quaternionic double pole formula in the separated form
has a kernel (\ref{2pole2-intro}) represented by the one-loop Feynman integral.
There is no doubt for us that these relations are only a tip of the iceberg,
and the other Feynman integrals also admit an interpretation via quaternionic
analysis and representation theory of the conformal group.
In fact, we make some explicit conjectures at the end of our paper.
Thus we come to the conclusion that the quaternionic analysis is very much
alive and well integrated with other areas of mathematics, since it might
contain a great portion -- if not the whole theory -- of Feynman integrals.
On the other hand, the latter theory -- a vast and central subject of physics
-- might not seem so disconcerting and unmotivated anymore for mathematicians,
and many of its beautiful results should be incorporated in an extended version
of quaternionic analysis.

For technical reasons the paper is organized slightly differently from the
order of this discussion.
In Section \ref{review} we review the classical quaternionic analysis:
Cauchy-Fueter formulas, conformal transformations, bases of harmonic and
regular functions, realizations of representations of the most degenerate
series of $SU(2,2)$. We also give a new proof of the Poisson formula based
on a representation theoretic argument. We conclude this section with
the derivation of the discrete spectrum of the hydrogen atom.
In Section \ref{M-section} we extend the results of the classical
quaternionic analysis to the Minkowski space.
Using the deformation of contours from $\BB M$ to $\BB H$, we first derive
the Fueter formula in $\BB M$ for the bounded cycles. Then we consider
a generalization of this formula to the case of unbounded cycles,
most notably the hyperboloids of two sheets, which has a natural
representation theoretic interpretation. Then we prove the Poisson formula
for the hyperboloids of two sheets and apply it to the derivation of the
continuous spectrum of the hydrogen atom.
In Section \ref{middle} we begin to study the so-called middle series of
representations of $SU(2,2)$, which is another (less) degenerate series of
representations. It can be realized in the spaces of all matrix- or
scalar-valued functions on $\HC$. Then we derive quaternionic analogues of the
Cauchy formula for the second order pole -- both in the Euclidean
and Minkowski formulation.
We identify the differential operator that appears in the Cauchy formula
for the second order pole with the Maxwell equations and explain the relation
between this formula and the polarization of vacuum in four-dimensional quantum
field theory. In Section \ref{decomposition-section} we continue to study
the middle series of representations of $SU(2,2)$ by decomposing the tensor
product of representations of the most degenerate series into irreducible
components, and we study the projections onto these components.
We find an explicit expression for the integral kernel of the projector to
the top component, which can be interpreted as quaternionic analogue of
the Cauchy formula for the separated double pole.
It turns out that this kernel is given by the one-loop Feynman integral.
We conclude this section with a conjecture on the relation between
integral kernels of higher projectors and Feynman integrals.

Finally, we would like to thank P.~Etingof, P.~Gressman, R.~Howe,
A.~Zeitlin and G.~Zuckerman for helpful discussions.
Also we would like to thank the referee for careful reading of the manuscript
and useful suggestions.
The first author was supported by the NSF grant DMS-0457444.

\section{Review of Quaternionic Analysis}  \label{review}

\subsection{Regular Functions}

The classical quaternions $\BB H$ form an algebra over $\BB R$ generated by
the units $e_0$, $e_1$, $e_2$, $e_3$ corresponding
to the more familiar $1,i,j,k$
(we reserve the symbol $i$ for $\sqrt{-1} \in \BB C$).
The multiplicative structure is determined by the rules
$$
e_0 e_i = e_i e_0 = e_i, \qquad i=0,1,2,3,
$$
$$
e_ie_j=-e_ie_j, \qquad 1 \le i< j \le 3,
$$
$$
(e_1)^2=(e_2)^2=(e_3)^2=e_1e_2e_3=-e_0,
$$
and the fact that $\BB H$ is a division ring.
We write an element $X \in \BB H$ as
$$
X = x^0e_0 + x^1e_1 + x^2e_2 + x^3e_3, \qquad x^0,x^1,x^2,x^3 \in \BB R,
$$
then its quaternionic conjugate is
$$
X^+ = x^0e_0 - x^1e_1 - x^2e_2 - x^3e_3,
$$
and it is easy to check that
$$
(XY)^+ = Y^+X^+, \qquad X,Y \in \BB H.
$$
We denote by $\re X$ the $e_0$-component of $X$:
$$
\re X = x^0 =\frac {X+X^+}2.
$$
We also describe the bilinear form and the corresponding quadratic norm
on quaternions as follows. For
$X = x^0e_0 + x^1e_1 + x^2e_2 + x^3e_3$ and
$Y = y^0e_0 + y^1e_1 + y^2e_2 + y^3e_3$ in $\BB H$,
\begin{equation}  \label{bilinear}
\langle X, Y \rangle = x^0y^0+x^1y^1+x^2y^2+x^3y^3
= \re (XY^+) = \re (X^+Y) \in \BB R,
\end{equation}
$$
N(X)= \langle X,X \rangle = XX^+ = X^+X \in \BB R.
$$
Hence the inverse of an element $X \in \BB H$, $X \ne 0$, is given by
$$
X^{-1} = \frac {X^+}{XX^+} = \frac {X^+}{N(X)}.
$$
We denote by $\BB H^{\times}$ the group of invertible elements
$\BB H \setminus \{0\}$. We define the Euclidean norm
$$
|X| = \sqrt{N(X)}, \qquad X \in \BB H,
$$
then
$$
|XY|=|X| \cdot |Y|, \qquad X,Y \in \BB H.
$$

Next we turn to the notion of quaternionic derivative.
Let $U \subset \BB H$ be an open subset and consider a function
of one quaternionic variable $f: U \to \BB H$ of class ${\cal C}^1$.
The ``correct'' notion of a regular function is obtained
by mimicking the $\overline{\partial} f =0$ equation in complex variables.
We define a formal analogue of the $\overline{\partial}$ operator:
$$
\nabla^+ = e_0 \frac{\partial}{\partial x^0}+
e_1 \frac{\partial}{\partial x^1}+e_2 \frac{\partial}{\partial x^2}+
e_3 \frac{\partial}{\partial x^3}.
$$
Because quaternions are not commutative, this differential operator may be
applied to a function on the left and on the right yielding different results.

\begin{df}
A ${\cal C}^1$-function $f: U \to \BB H$ is {\em left-regular} if it satisfies
$$
(\nabla^+f)(X) = e_0 \frac{\partial f}{\partial x^0}(X)+
e_1 \frac{\partial f}{\partial x^1}(X)+e_2 \frac{\partial f}{\partial x^2}(X)+
e_3 \frac{\partial f}{\partial x^3}(X) =0, \qquad \forall X \in U.
$$
Similarly, $f$ is {\em right-regular} if
$$
(f\nabla^+)(X) = \frac{\partial f}{\partial x^0}(X)e_0 +
\frac{\partial f}{\partial x^1}(X)e_1+\frac{\partial f}{\partial x^2}(X)e_2+
\frac{\partial f}{\partial x^3}(X)e_3 =0, \qquad \forall X \in U.
$$
\end{df}

These functions are quaternionic analogues of complex holomorphic functions.
Note that in general the product of two left-regular (or right-regular)
functions is not left-regular (respectively right-regular).

\begin{prop}  \label{lrswitch}
A ${\cal C}^1$-function $f: U \to \BB H$ is left-regular if and only if
$g(X)=f^+(X^+)$ is right-regular on $U^+=\{X \in \BB H ;\: X^+ \in U\}$.
A ${\cal C}^1$-function $g: U \to \BB H$ is right-regular if and only if
$f(X)=g^+(X^+)$ is left-regular on $U^+$.
\end{prop}

Note that the four-dimensional Laplacian
\begin{equation}  \label{Laplacian}
\square = \frac{\partial^2}{(\partial x^0)^2}+
\frac{\partial^2}{(\partial x^1)^2} + \frac{\partial^2}{(\partial x^2)^2}+
\frac{\partial^2}{(\partial x^3)^2}
\end{equation}
can be factored as
$$
\square = \nabla\nabla^+ = \nabla^+\nabla,
$$
where
$$
\nabla = (\nabla^+)^+ = e_0 \frac{\partial}{\partial x^0}-
e_1 \frac{\partial}{\partial x^1} - e_2 \frac{\partial}{\partial x^2}-
e_3 \frac{\partial}{\partial x^3}.
$$
The operator $\nabla$ may be applied to functions on the left and on the right
and can be thought of as a quaternionic analogue of the complex $\partial$
operator.

\begin{prop}  \label{k}
\begin{enumerate}
\item
The function $k_0(X) =_{def} \frac 1{N(X)}$ is harmonic on
$\BB H^{\times} = \BB H \setminus \{0\}$.
\item
The function $k(X) =_{def} \frac {X^{-1}}{N(X)}$ satisfies
$$
k(X) = -\frac 12 (\nabla k_0)(X) = -\frac 12 (k_0\nabla)(X).
$$
\item
The function $k(X) = \frac {X^{-1}}{N(X)}$ is both left- and right-regular on
$\BB H^{\times} = \BB H \setminus \{0\}$.
\end{enumerate}
\end{prop}

\subsection{Cauchy-Fueter Formulas}

First we regard $\BB H$ as $\BB R^4$ and form the exterior algebra denoted by
$\Lambda_{\BB R}$. Thus $\Lambda_{\BB R}$ is a graded algebra generated by
$dx^0, dx^1, dx^2, dx^3$. Then we form
$\Lambda_{\BB H} = \Lambda_{\BB R} \otimes_{\BB R} \BB H$.
We think of $\Lambda_{\BB H}$ as quaternionic-valued differential forms on
$\BB R^4 \simeq \BB H$. By definition, $\Lambda_{\BB H}$ is the graded
$\BB H$-algebra generated by elements $e_idx^j$, and we have
$e_idx^j = dx^je_i$, for all $i,j =0,1,2,3$.
We introduce an element
$$
dX = e_0dx^0 + e_1dx^1 + e_2dx^2 + e_3dx^3 \quad \in \Lambda_{\BB H}^1.
$$
The bilinear form (\ref{bilinear}) determines the Hodge star operator on
$\Lambda_{\BB H}$, in particular we get an isomorphism
$*: \Lambda_{\BB H}^1 \tilde\longrightarrow \Lambda_{\BB H}^3$, and we set
$$
Dx = *(dX) \quad \in \Lambda_{\BB H}^3.
$$
In other words, set
$dV=dx^0 \wedge dx^1 \wedge dx^2 \wedge dx^3 \in \Lambda_{\BB H}^4$,
then $Dx$ is determined by the property
$$
\re( h_1^+ \cdot Dz(h_2,h_3,h_4)) = dV(h_1,h_2,h_3,h_4),
\qquad \forall h_1,h_2,h_3,h_4 \in \BB H.
$$
Hence it is easy to see that
\begin{equation}  \label{Dx-explicit}
Dx= e_0 dx^1 \wedge dx^2 \wedge dx^3 - e_1 dx^0 \wedge dx^2 \wedge dx^3
+ e_2 dx^0 \wedge dx^1 \wedge dx^3 - e_3 dx^0 \wedge dx^1 \wedge dx^2.
\end{equation}

\begin{prop}  \label{closed}
Let $f, g: U \to \BB H$ be two ${\cal C}^1$-functions. Then
$$
d(Dx \cdot f) = - Dx \wedge df = (\nabla^+ f)dV
\qquad \text{and} \qquad
d(g \cdot Dx) = dg \wedge Dx = (g \nabla^+)dV.
$$
Combining the two formulas we get
$$
d(g \cdot Dx \cdot f) = \bigl( (g \nabla^+)f + g(\nabla^+f) \bigr) dV.
$$
\end{prop}

\begin{cor}  \label{alt-reg}
A ${\cal C}^1$-function $f: U \to \BB H$ is left-regular if and only if
$Dx \wedge df =0$;
a ${\cal C}^1$-function $g: U \to \BB H$ is right-regular if
and only if $dg \wedge Dx =0$.
\end{cor}

Let $U \subset \BB H$ be an open region with piecewise smooth boundary
$\partial U$. We give a canonical orientation to $\partial U$ as follows.
The positive orientation of $U$ is determined by $\{e_0, e_1, e_2, e_3 \}$.
Pick a smooth point $p \in \partial U$ and let $\overrightarrow{n_p}$
be a non-zero vector in $T_p \BB H$ perpendicular to $T_p\partial U$ and
pointing outside of $U$.
Then $\{\overrightarrow{\tau_1}, \overrightarrow{\tau_2},
\overrightarrow{\tau_3}\} \subset T_p \partial U$ is positively oriented
in $\partial U$ if and only if
$\{\overrightarrow{n_p}, \overrightarrow{\tau_1}, \overrightarrow{\tau_2},
\overrightarrow{\tau_3}\}$ is positively oriented in $\BB H$.

Fix a $Y \in \BB H$, an $R>0$, and denote by $S^3_R(Y)$
the three-dimensional sphere of radius $R$ centered at $Y$:
$$
S^3_R(Y) = \{ X \in \BB H ;\: |X-Y|=R \}.
$$
We orient the sphere as the boundary of the open ball of radius $R$ centered
at $Y$.

\begin{lem}  \label{restriction}
We have:
$$
Dx \Bigl|_{T_p\partial U} = \overrightarrow{n_p} dS \Bigl|_{T_p\partial U},
$$
where $dS$ is the usual Euclidean volume element on $\partial U$.

In particular, the restriction of $Dx$ to the sphere $S^3_R(Y)$
is given by the formula
$$
Dx \Bigl|_{S^3_R(Y)} = \frac{X-Y}{|X-Y|} \,dS = \frac{X-Y}R \,dS,
$$
where $dS$ is the usual Euclidean volume element on $S^3_R(Y)$.
\end{lem}

First we state an analogue of the Cauchy theorem.

\begin{prop}  \label{CFthm}
Let $U \subset \BB H$ be an open bounded subset with piecewise ${\cal C}^1$
boundary $\partial U$. Suppose that $f(X)$ and $g(X)$ are
$\BB H$-valued ${\cal C}^1$-functions defined in a neighborhood of the
closure $\overline{U}$. Then
$$
\int_{\partial U} g \cdot Dx \cdot f =
\int_U \bigl( (g \nabla^+)f + g(\nabla^+f) \bigr) \,dV.
$$
\end{prop}

\begin{cor}  \label{zero}
Let $U \subset \BB H$ be an open bounded subset with piecewise ${\cal C}^1$
boundary $\partial U$. Suppose that $f(X)$ is left-regular and
$g(X)$ is right-regular on a neighborhood of the closure $\overline{U}$. Then
$$
\int_{\partial U} g \cdot Dx \cdot f =0.
$$
\end{cor}

Now we are ready to state the first main theorem of quaternionic analysis,
it is analogous to the Cauchy formula (\ref{Cauchy-intro}). Recall the function
$k$ introduced in Proposition \ref{k}. If $X_0 \in \BB H$ is fixed,
\begin{equation}  \label{k(z-w)}
k(X-X_0) = \frac {(X-X_0)^{-1}}{N(X-X_0)}
\end{equation}
is left- and right-regular on $\BB H \setminus \{X_0\}$.

\begin{thm} [Cauchy-Fueter Formulas, \cite{F1,F2}]  \label{Fueter}
Let $U \subset \BB H$ be an open bounded subset with piecewise ${\cal C}^1$
boundary $\partial U$. Suppose that $f(X)$ is left-regular on a
neighborhood of the closure $\overline{U}$, then
$$
\frac 1 {2\pi^2} \int_{\partial U} k(X-X_0) \cdot Dx \cdot f(X) =
\begin{cases}
f(X_0) & \text{if $X_0 \in U$;} \\
0 & \text{if $X_0 \notin \overline{U}$.}
\end{cases}
$$
If $g(X)$ is right-regular on a neighborhood of the closure $\overline{U}$,
then
$$
\frac 1 {2\pi^2} \int_{\partial U} g(X) \cdot Dx \cdot k(X-X_0) =
\begin{cases}
g(X_0) & \text{if $X_0 \in U$;} \\
0 & \text{if $X_0 \notin \overline{U}$.}
\end{cases}
$$
\end{thm}

The complex variable corollaries of the Cauchy formula have their
quaternionic analogues with nearly identical proofs.
Thus left- and right-regular functions are smooth, harmonic and, in fact,
real-analytic. If $f: \BB H \to \BB H$ is a left- (or right-) regular
function which is bounded on $\BB H$, then it is constant.

\subsection{Matrix Realizations, Spinors and Spinor-Valued Functions}

We realize $\BB H$ as a subalgebra of the algebra of $2 \times 2$
complex matrices by identifying the units $e_0$, $e_1$, $e_2$, $e_3$ with
$$
e_0 \longleftrightarrow \begin{pmatrix} 1 & 0 \\ 0 & 1 \end{pmatrix}, \qquad
e_1 \longleftrightarrow \begin{pmatrix} 0 & -i \\ -i & 0 \end{pmatrix}, \qquad
e_2 \longleftrightarrow \begin{pmatrix} 0 & -1 \\ 1 & 0 \end{pmatrix}, \qquad
e_3 \longleftrightarrow \begin{pmatrix} -i & 0 \\ 0 & i \end{pmatrix}.
$$
Thus
$$
X = x^0e_0 + x^1e_1 + x^2e_2 + x^3e_3
\quad \longleftrightarrow \quad
\begin{pmatrix} x^0-ix^3 & -ix^1-x^2 \\ -ix^1+x^2 & x^0+ix^3 \end{pmatrix}
= \begin{pmatrix} x_{11} & x_{12} \\ x_{21} & x_{22} \end{pmatrix},
$$
and we get an identification of quaternions:
$$
\BB H \:\simeq\:
\Bigl\{
\begin{pmatrix} a & b \\ -\overline{b} & \overline{a} \end{pmatrix}
\in \mathfrak{gl}(2,\BB C);\: a,b \in \BB C \Bigr\}.
$$
In these new coordinates
\begin{equation}  \label{conj-matrix-form}
\begin{pmatrix} x_{11} & x_{12} \\ x_{21} & x_{22} \end{pmatrix}^+
= \begin{pmatrix} x_{22} & -x_{12} \\ -x_{21} & x_{11} \end{pmatrix}
= \begin{pmatrix} \overline{x_{11}} & \overline{x_{21}} \\
\overline{x_{12}} & \overline{x_{22}} \end{pmatrix}^T
\end{equation}
-- matrix transpose followed by complex conjugation -- and
$$
N\begin{pmatrix} x_{11} & x_{12} \\ x_{21} & x_{22} \end{pmatrix}
= \det \begin{pmatrix} x_{11} & x_{12} \\ x_{21} & x_{22} \end{pmatrix}
=x_{11}x_{22}-x_{12}x_{21}.
$$
Under this identification, the unit sphere in $\BB H$ gets identified with
$SU(2)$ in $GL(2,\BB C)$:
\begin{multline}  \label{SU(2)}
S^3 = \{ X \in \BB H ;\: |X|=1 \}  \\
\simeq
\Bigl\{ \begin{pmatrix} a & b \\ -\overline{b} & \overline{a} \end{pmatrix}
\in GL(2,\BB C);\: a,b \in \BB C ,\: \det
\begin{pmatrix} a & b \\ -\overline{b} & \overline{a} \end{pmatrix} =1 \Bigr\}
= SU(2).
\end{multline}

Writing $\partial_i$ for $\frac{\partial}{\partial x^i}$ and
$\partial_{ij}$ for $\frac{\partial}{\partial x_{ij}}$, a simple change
of variables computation shows
\begin{align*}
\nabla^+ &= 
\begin{pmatrix} \partial_0-i\partial_3 & -i\partial_1-\partial_2 \\
-i\partial_1+\partial_2 & \partial_0+i\partial_3 \end{pmatrix}
=
2 \begin{pmatrix} \partial_{22} & -\partial_{21} \\
-\partial_{12} & \partial_{11} \end{pmatrix},  \\
\nabla &= 
2 \begin{pmatrix} \partial_{22} & -\partial_{21} \\
-\partial_{12} & \partial_{11} \end{pmatrix}^+
= 2 \begin{pmatrix} \partial_{11} & \partial_{21} \\
\partial_{12} & \partial_{22} \end{pmatrix},  \\
\square &=\nabla \nabla^+ =
4 \begin{pmatrix} \partial_{11} & \partial_{21} \\
\partial_{12} & \partial_{22} \end{pmatrix}
\begin{pmatrix} \partial_{22} & -\partial_{21} \\
-\partial_{12} & \partial_{11} \end{pmatrix}
= 4 (\partial_{11}\partial_{22} - \partial_{12}\partial_{21}).
\end{align*}

Let $\BB S$ be the natural two-dimensional complex representation of
the algebra $\BB H$; it can be realized as a column of two complex
numbers. We have the standard action of $\mathfrak{gl}(2,\BB C)$ on $\BB S$
by matrix multiplication on the left and hence left action of $\BB H$.
Similarly, we denote by $\BB S'$ the dual space of $\BB S$, this time
realized as a row of two numbers. We have the standard action of
$\mathfrak{gl}(2,\BB C)$ on $\BB S'$ by multiplication on the right and
hence right action of $\BB H$.
We define a complex bilinear pairing
$$
\langle \,\cdot\, ,  \,\cdot\, \rangle : \BB S' \times \BB S \to \BB C,
\qquad
\Bigl\langle (s_1', s_2'), \begin{pmatrix} s_1 \\ s_2 \end{pmatrix}
\Bigr\rangle = s_1's_1 + s_2's_2.
$$
The pairing is $\BB H$-invariant in the sense that
$$
\langle s'X,s \rangle = \langle s',Xs \rangle,
\qquad s' \in \BB S',\: s \in \BB S ,\: X \in \BB H.
$$
Next we identify $\BB S \otimes_{\BB C} \BB S'$ with
$\BB H_{\BB C} = \BB H \otimes_{\BB R} \BB C$.
Indeed, a vector space tensored with its dual is canonically isomorphic
to the space of endomorphisms of that space.
Pure tensors $s \otimes s'$ act on $\BB S$ by
$$
(s \otimes s')(t) = s \langle s', t \rangle,
\qquad s,t \in \BB S ,\: s' \in \BB S'.
$$
Heuristically, 
``$\BB S$ and $\BB S'$ are the square roots of $\BB H_{\BB C}$.''

We have $\BB C$-antilinear maps $\BB S \to \BB S'$ and $\BB S' \to \BB S$
-- matrix transposition followed by complex conjugation --
which are similar to quaternionic conjugation and so, by abuse of notation,
we use the same symbol to denote these maps:
$$
\begin{pmatrix} s_1 \\ s_2 \end{pmatrix}^+ = (\overline{s_1}, \overline{s_2}),
\qquad
(s'_1,s'_2)^+ =\begin{pmatrix} \overline{s'_1} \\ \overline{s'_2}\end{pmatrix},
\qquad
\begin{pmatrix} s_1 \\ s_2 \end{pmatrix} \in \BB S ,\: (s'_1,s'_2) \in \BB S'.
$$
Note that $(Xs)^+ = s^+ X^+$ and $(s'X)^+ = X^+ s'^+$
for all $X \in \BB H$, $s \in \BB S$, $s' \in \BB S'$.

Let $U \subset \BB H$ and consider ${\cal C}^1$-functions
$f: U \to \BB S$ and $g: U \to \BB S'$. We say that $f(X)$ is a
{\em left-regular spinor-valued} function if
$$
(\nabla^+f)(X)= e_0 \frac{\partial f}{\partial x^0}(X)+
e_1 \frac{\partial f}{\partial x^1}(X)+e_2 \frac{\partial f}{\partial x^2}(X)+
e_3 \frac{\partial f}{\partial x^3}(X) = 0, \qquad \forall X \in U,
$$
where multiplication by $e_0,e_1,e_2,e_3$ is understood in the sense of
left action of $\BB H$ on $\BB S$.
Dually, $g(X)$ is a {\em right-regular spinor-valued} function if
$$
(g\nabla^+)(X)= \frac{\partial f}{\partial x^0}(X) e_0 +
\frac{\partial f}{\partial x^1}(X)e_1 +\frac{\partial f}{\partial x^2}(X)e_2+
\frac{\partial f}{\partial x^3}(X) e_3 = 0, \qquad \forall X \in U,
$$
and the multiplication by $e_0,e_1,e_2,e_3$ is understood in the sense of
right action of $\BB H$ on $\BB S'$.
We have analogues of Propositions \ref{lrswitch}, \ref{CFthm},
Corollary \ref{zero} and the Cauchy-Fueter formula (Theorem \ref{Fueter})
for spinor-valued functions.

\subsection{Conformal Transformations}

We denote by $\widehat{\BB H} = \BB H \cup \{\infty\}$ the one-point
compactification of $\BB H$ as a four-dimensional sphere.
It can be realized as the one-dimensional
quaternionic projective space. However, since quaternions are not commutative,
there are two different kinds of projective spaces:
$$
\BB P^1 \BB H_l = \Bigl\{ \begin{pmatrix} X_1 \\ X_2 \end{pmatrix} ;\:
X_1, X_2 \in \BB H, \: X_1 \ne 0 \text{ or } X_2 \ne 0 \Bigr\} / \sim_l
$$
and
$$
\BB P^1 \BB H_r = \Bigl\{ \begin{pmatrix} X_1 \\ X_2 \end{pmatrix} ;\:
X_1, X_2 \in \BB H, \: X_1 \ne 0 \text{ or } X_2 \ne 0 \Bigr\} / \sim_r,
$$
where
$$
\begin{pmatrix} X_1 \\ X_2 \end{pmatrix} \sim_l
\begin{pmatrix} X_1a \\ X_2a \end{pmatrix}
\qquad \text{and} \qquad \begin{pmatrix} X_1 \\ X_2 \end{pmatrix} \sim_r
\begin{pmatrix} aX_1 \\ aX_2 \end{pmatrix}
\qquad \text{for all $a \in \BB H^{\times}$}.
$$
Then, as $\BB R$-manifolds,
$\BB P^1 \BB H_l  \simeq \widehat{\BB H} \simeq \BB P^1 \BB H_r$.

Next we consider the group $GL(2,\BB H)$ consisting of invertible $2 \times 2$
matrices with quaternionic entries.
This group contains $GL(2,\BB R)$ in the obvious way and also the
multiplicative group $\BB H^{\times} \times \BB H^{\times}$ as the
diagonal matrices $\begin{pmatrix} a & 0 \\ 0 & d \end{pmatrix}$,
$a, d  \in \BB H \setminus \{0\}$.
The group $GL(2,\BB H)$ acts on $\BB P^1 \BB H_l$ and $\BB P^1 \BB H_r$
by multiplications on the left and right respectively.
Hence we get two different-looking actions on
$\widehat{\BB H} = \BB H \cup \{\infty\}$
by conformal (fractional linear) transformations:
\begin{align*}
\pi_l(h): \: X \mapsto & (aX+b)(cX+d)^{-1}, \qquad
h^{-1} = \begin{pmatrix} a & b \\ c & d \end{pmatrix} \in GL(2, \BB H); \\
\pi_r(h): \: X \mapsto & (a'-Xc')^{-1}(-b'+Xd'), \qquad
h = \begin{pmatrix} a' & b' \\ c' & d' \end{pmatrix} \in GL(2, \BB H).
\end{align*}
These two actions coincide.

\begin{lem}  \label{X-Y_difference}
For $h = \begin{pmatrix} a' & b' \\ c' & d' \end{pmatrix}$ in $GL(2, \BB H)$
with $h^{-1} = \begin{pmatrix} a & b \\ c & d \end{pmatrix}$,
let $\tilde X = (aX+b)(cX+d)^{-1}$ and $\tilde Y = (aY+b)(cY+d)^{-1}$.
Then
\begin{align*}
(\tilde X - \tilde Y) &= (a'-Yc')^{-1} \cdot(X-Y) \cdot (cX+d)^{-1}  \\
&= (a'-Xc')^{-1} \cdot(X-Y) \cdot (cY+d)^{-1}.
\end{align*}
\end{lem}

\begin{prop} \label{Dx-pullback}
Let $h = \begin{pmatrix} a' & b' \\ c' & d' \end{pmatrix} \in GL(2, \BB H)$
with $h^{-1} = \begin{pmatrix} a & b \\ c & d \end{pmatrix}$.
Then the pull-back of $Dx$
$$
\pi_l(h)^*(Dx) =
\frac {(a'-Xc')^{-1}}{N(a'-Xc')} \cdot Dx \cdot \frac {(cX+d)^{-1}}{N(cX+d)}.
$$
\end{prop}

\begin{prop}  \label{conformal}
\begin{enumerate}
\item
The space of left-regular $\BB S$-valued functions on $\widehat{\BB H}$
with isolated singularities is invariant under the group of conformal
transformations
\begin{multline}  \label{spinor_left_action}
\pi_l(h): \: f(X) \mapsto (\pi_l(h)f)(X) =
\frac {(cX+d)^{-1}}{N(cX+d)} \cdot f \bigl( (aX+b)(cX+d)^{-1} \bigr),  \\
h^{-1} = \begin{pmatrix} a & b \\ c & d \end{pmatrix} \in GL(2, \BB H).
\end{multline}
\item
The space of right-regular $\BB S'$-valued functions on $\widehat{\BB H}$
with isolated singularities is invariant under the group of conformal
transformations
\begin{multline}  \label{spinor_right_action}
\pi_r(h): \: g(X) \mapsto (\pi_r(h)g)(X) =
g \bigl( (a'-Xc')^{-1}(-b'+Xd') \bigr) \cdot \frac {(a'-Xc')^{-1}}{N(a'-Xc')},
\\
h = \begin{pmatrix} a' & b' \\ c' & d' \end{pmatrix} \in GL(2, \BB H).
\end{multline}
\end{enumerate}
\end{prop}

Fix an open bounded subset $U \subset \BB H$ with piecewise ${\cal C}^1$
boundary $\partial U$.
We denote by ${\cal V}(\B{U})$ the space of left-regular
spinor-valued functions with isolated singularities defined in a neighborhood
of $\B{U}$ with no singularities on $\partial U$.
And we denote by ${\cal V}'(\B{U})$ the space of right-regular
spinor-valued functions with isolated singularities defined in a neighborhood
of $\B{U}$ with no singularities on $\partial U$.
We have a natural pairing between ${\cal V}(\overline{U})$ and
${\cal V}'(\overline{U})$:
\begin{equation}  \label{fgpairing}
\langle g, f \rangle_U = \int_{\partial U} \langle g(X), Dx \cdot f(X) \rangle
= \int_{\partial U} \langle g(X) \cdot Dx, f(X) \rangle \quad \in \BB C,
\qquad f \in {\cal V}(\B{U}), \: g \in {\cal V}'(\B{U}).
\end{equation}
Proposition \ref{CFthm} and Corollary \ref{zero} apply to spinor-valued
functions, and so it is clear that the integral (\ref{fgpairing})
stays unchanged if the contour of integration $\partial U$ is deformed
without crossing the singularities of $f$ and $g$.

Differentiating the actions $\pi_l$ and $\pi_r$ of $GL(2,\BB H)$ on
spinor-valued regular functions we get representations of the Lie algebra
$\mathfrak{gl}(2,\BB H)$ which we still denote by $\pi_l$ and $\pi_r$.

\begin{prop}  \label{pseudo-unitary}
Suppose that $f \in {\cal V}(\overline{U})$ and
$g \in {\cal V}'(\overline{U})$, then
$$
\langle \pi_r(x) g, f \rangle_U + \langle g, \pi_l(x) f \rangle_U =0,
\qquad \forall x \in \mathfrak{gl}(2,\BB H).
$$
\end{prop}

This result follows immediately from

\begin{prop}
The actions $\pi_l$ and $\pi_r$ preserve the pairing
$\langle \,\cdot\,, \,\cdot\, \rangle_U$:
$$
\langle \pi_r(h) g, \pi_l(h)f \rangle_U = \langle g, f \rangle_U,
\qquad f \in {\cal V}(\overline{U}), \: g \in {\cal V}'(\overline{U}),
$$
for all $h$ in a sufficiently small neighborhood of
$\begin{pmatrix} 1 & 0 \\ 0 & 1 \end{pmatrix} \in GL(2, \BB H)$
depending on $f$ and $g$.
\end{prop}

We also have two actions of $GL(2,\BB H)$ on the space of harmonic functions:

\begin{prop}  \label{action-harmonic}
We have two different actions of the group $GL(2, \BB H)$ acts on the space
of $\BB C$-valued harmonic functions on $\widehat{\BB H}$ with isolated
singularities:
\begin{multline}  \label{left_action}
\pi^0_l(h): \: \phi(X) \mapsto \bigl( \pi^0_l(h)\phi \bigr)(X) =
\frac 1{N(cX+d)} \cdot \phi \bigl( (aX+b)(cX+d)^{-1} \bigr),  \\
h^{-1} = \begin{pmatrix} a & b \\ c & d \end{pmatrix} \in GL(2, \BB H),
\end{multline}
\begin{multline}  \label{right_action}
\pi^0_r(h): \: \phi(X) \mapsto \bigl( \pi^0_r(h)\phi \bigr)(X) =
\frac 1{N(a'-Xc')} \cdot \phi \bigl( (a'-Xc')^{-1}(-b'+Xd') \bigr),  \\
h = \begin{pmatrix} a' & b' \\ c' & d' \end{pmatrix} \in GL(2, \BB H).
\end{multline}
These two actions coincide on $SL(2, \BB H)$, which is defined
as the connected Lie subgroup of $GL(2, \BB H)$ whose Lie algebra is
$$
\mathfrak{sl}(2,\BB H) = \{ x \in \mathfrak{gl}(2,\BB H) ;\: \re (\tr x) =0 \}.
$$
This Lie algebra is isomorphic to $\mathfrak{so}(5,1)$.
\end{prop}

Differentiating these two actions we obtain two representations of
$\mathfrak{gl}(2, \BB H)$ on the space of harmonic functions which
we still denote by $\pi^0_l$ and $\pi^0_r$. These two representations
agree on $\mathfrak{sl}(2, \BB H)$.

\begin{df}  \label{reg-infinity-def}
Let $\phi$ be a $\BB C$-valued harmonic function on $\widehat{\BB H}$
with isolated singularities. We say that $\phi$ is {\em regular at infinity}
if the harmonic function
$$
\pi^0_l \begin{pmatrix} 0 & 1 \\ 1 & 0 \end{pmatrix} \phi =
\pi^0_r \begin{pmatrix} 0 & 1 \\ 1 & 0 \end{pmatrix} \phi =
\frac 1{N(X)} \cdot \phi(X^{-1})
$$
is regular at the origin.

Similarly, we say that a left-regular function $f$ or a right-regular
function $g$ on $\widehat{\BB H}$ with isolated singularities is 
{\em regular at infinity} if
$$
\pi_l \begin{pmatrix} 0 & 1 \\ 1 & 0 \end{pmatrix} f =
\frac {X^{-1}}{N(X)} \cdot f(X^{-1})
\qquad \text{or} \qquad
\pi_r \begin{pmatrix} 0 & 1 \\ 1 & 0 \end{pmatrix} g =
g(X^{-1}) \cdot \frac {X^{-1}}{N(X)}
$$
is regular at the origin.
\end{df}

In the following lemma we write $\mathfrak{gl}(2,\BB H)$ as
$$
\mathfrak{gl}(2,\BB H) = \biggl\{
\begin{pmatrix} A & B \\ C & D \end{pmatrix} ;\: A,B,C,D \in \BB H \biggr\}
$$
and think of $A,B,C,D$ as $2\times 2$ matrices. Also, let
$$
\partial = \begin{pmatrix} \partial_{11} & \partial_{21} \\
\partial_{12} & \partial_{22} \end{pmatrix} = \frac 12 \nabla.
$$

\begin{lem}  \label{gl(2,H)-action0}
The Lie algebra actions $\pi^0_l$ and $\pi^0_r$ of $\mathfrak{gl}(2,\BB H)$
on the space of harmonic functions are given by
\begin{align*}
\pi_l^0 \begin{pmatrix} A & 0 \\ 0 & 0 \end{pmatrix} &:
\phi \mapsto \tr \bigl( A \cdot (-X \cdot \partial \phi) \bigr)  \\
\pi_r^0 \begin{pmatrix} A & 0 \\ 0 & 0 \end{pmatrix} &:
\phi \mapsto \tr \bigl( A \cdot (-X \cdot \partial \phi - \phi) \bigr)  \\
\pi_l^0 \begin{pmatrix} 0 & B \\ 0 & 0 \end{pmatrix} =
\pi_r^0 \begin{pmatrix} 0 & B \\ 0 & 0 \end{pmatrix} &:
\phi \mapsto \tr \bigl( B \cdot (-\partial \phi ) \bigr)  \\
\pi_l^0 \begin{pmatrix} 0 & 0 \\ C & 0 \end{pmatrix} =
\pi_r^0 \begin{pmatrix} 0 & 0 \\ C & 0 \end{pmatrix} &:
\phi \mapsto \tr \Bigl( C \cdot \bigl(
X \cdot (\partial \phi) \cdot X + X\phi \bigr) \Bigr)  \\
&: \phi \mapsto \tr \Bigl( C \cdot \bigl(
X \cdot \partial (X\phi) \bigr) - X\phi \Bigr)  \\
\pi_l^0 \begin{pmatrix} 0 & 0 \\ 0 & D \end{pmatrix} &: \phi \mapsto
\tr \Bigl( D \cdot \bigl( (\partial \phi) \cdot X + \phi \bigr) \Bigr)  \\
&: \phi \mapsto \tr \Bigl( D \cdot \bigl( \partial (X\phi) - \phi \bigr)
\Bigr)  \\
\pi_r^0 \begin{pmatrix} 0 & 0 \\ 0 & D \end{pmatrix} &: \phi \mapsto
\tr \Bigl( D \cdot \bigl( (\partial \phi) \cdot X \bigr) \Bigr)  \\
&: \phi \mapsto \tr \Bigl( D \cdot \bigl( \partial (X\phi) - 2\phi \bigr)
\Bigr).
\end{align*}
\end{lem}

\subsection{Laurent Polynomials}

We use (\ref{SU(2)}) to identify $SU(2)$ with unit quaternions.
In particular, the $\BB H$-actions on $\BB S$ and $\BB S'$
restrict to $SU(2)$, and the irreducible representations of
$SU(2)$ are realized in the symmetric powers of $\BB S$.
Let $P(s'_1,s'_2)$ be a polynomial function on
$\BB S' = \{ (s'_1,s'_2) ;\: s'_1,s'_2 \in \BB C\}$, then an element
$X =\begin{pmatrix} x_{11} & x_{12} \\ x_{21} & x_{22} \end{pmatrix} \in SU(2)$
acts on $P(s'_1,s'_2)$ by
$$
P(s'_1,s'_2) \mapsto (\tau(X)P)(s'_1,s'_2)=
P(s'_1x_{11}+s'_2x_{21}, s'_1x_{12}+s'_2x_{22}).
$$
The irreducible representation $V_l$ of $SU(2)$ of level $l$,
$l=0 , \frac 12 , 1 , \frac 32, \dots$,
is realized in the space of homogeneous polynomials of degree $2l$.
We denote the restriction of $\tau$ to $V_l$ by $\tau_l$.
As a basis of $V_l$ we can choose the monomials
$$
P_m = (s'_1)^{l-m} (s'_2)^{l+m},
\qquad m=-l,-l+1,\dots,l;
$$
there are exactly $\dim_{\BB C} V_l = 2l+1$ such monomials.
Next we consider the matrix coefficients $t^l_{n\,\underline{m}}(X)$ of
the linear transformation $\tau_l(X)$ with respect to the basis
$\{P_m; \: m=-l,-l+1,\dots,l \}$.
That is $t^l_{n\,\underline{m}}(X)$'s are the $\BB C$-valued functions
on $SU(2)$ uniquely determined by the equation
\begin{multline*}
(\tau_l(X) P_m)(s'_1,s'_2) = 
(s'_1x_{11}+s'_2x_{21})^{l-m} (s'_1x_{12}+s'_2x_{22})^{l+m} \\
=\sum_{n=-l}^l t^l_{n\,\underline{m}}(X) \cdot (s'_1)^{l-n}(s'_2)^{l+n}
=\sum_{n=-l}^l t^l_{n\,\underline{m}}(X) \cdot P_n(s'_1,s'_2), \qquad
X=\begin{pmatrix} x_{11} & x_{12} \\ x_{21} & x_{22} \end{pmatrix}.
\end{multline*}
One can easily obtain the following formula for the matrix coefficients:
\begin{equation}  \label{int_m}
t^l_{n\,\underline{m}}(X) = \frac 1{2\pi i}
\oint (sx_{11}+x_{21})^{l-m} (sx_{12}+x_{22})^{l+m} s^{-l+n} \,\frac{ds}s,
\end{equation}
where the integral is taken over a loop in $\BB C$ going once around the origin
in the counterclockwise direction (cf. \cite{V}).
Although the matrix coefficients are defined as functions on $SU(2)$,
they naturally extend as polynomial functions to $\BB H$.
Differentiating (\ref{int_m}) under the integral sign we obtain:

\begin{lem}
The matrix coefficients $t^l_{n\,\underline{m}}(X)$ are harmonic:
$$
\square t^l_{n\,\underline{m}}(X) =0,
\qquad m,n =-l,-l+1,\dots,l.
$$
\end{lem}

We identify the complex-valued polynomial functions on $\BB H$ with
polynomial functions on $\BB H \otimes \BB C$ and denote both spaces by
$\BB C [z^0,z^1,z^2,z^3]$.
We define the spaces of harmonic polynomials on $\BB H$ and $\BB H^{\times}$:
$$
{\cal H}^+ = \{ \phi \in \BB C [z^0,z^1,z^2,z^3] ;\: \square \phi = 0 \},
\qquad
{\cal H} = \{ \phi \in \BB C [z^0,z^1,z^2,z^3, N(Z)^{-1}] ;\:
\square \phi = 0 \}.
$$
We denote by $\BB C [z^0,z^1,z^2,z^3]_d$ and
$\BB C [z^0,z^1,z^2,z^3, N(Z)^{-1}]_d$ the subspaces of degree $d$,
where the variables $z^0,z^1,z^2,z^3$ are assigned degree $+1$ and
$N(Z)^{-1}$ is assigned degree $-2$. Define
\begin{align*}
{\cal H}(n) &=
\{ \phi \in \BB C [z^0,z^1,z^2,z^3,N(Z)^{-1}]_n ;\: \square \phi=0 \},
\qquad n \in \BB Z,  \\
{\cal H} &= \bigoplus_{n \in \BB Z} {\cal H}(n)
= {\cal H}^+ \oplus {\cal H}^-, \qquad \text{where}  \\
{\cal H}^+ &= \bigoplus_{n \in \BB Z, \: n \ge 0} {\cal H}(n), \qquad
{\cal H}^- = \bigoplus_{n \in \BB Z, \: n<0} {\cal H}(n).
\end{align*}
It is natural to call $\BB C [z^0,z^1,z^2,z^3, N(Z)^{-1}]$ the space of
Laurent polynomials on $\BB H^{\times}$.

We use (\ref{SU(2)}) to identify $SU(2)$ with unit quaternions
and realize $SU(2) \times SU(2)$ as diagonal elements of $GL(2,\BB H)$.
Thus the actions (\ref{left_action}) and (\ref{right_action}) of
$GL(2,\BB H)$ on harmonic functions both restrict to
\begin{equation}  \label{SL2C-action}
\widetilde{\pi}^0(a \times b) : \phi(Z) \mapsto \phi(a^{-1} Z b),
\qquad a, b \in SU(2), \: \phi \in {\cal H}.
\end{equation}
Clearly, $\widetilde{\pi}^0$ preserves each ${\cal H}(n)$, $n \in \BB Z$.

\begin{prop}  \label{extension}
We have $\BB C [z^0,z^1,z^2,z^3] \simeq {\cal H}^+ \cdot \BB C[N(Z)]$ and
$$
\BB C [z^0,z^1,z^2,z^3,N(Z)^{-1}] \simeq
{\cal H}^+ \cdot \BB C[N(Z),N(Z)^{-1}].
$$
Moreover, the restriction map $\phi \mapsto \phi \bigl|_{S^3}$ induces
the following isomorphisms:
$$
{\cal H}^+ \simeq
\Bigl \{ \begin{matrix} \text{linear span of matrix}  \\
\text{coefficients of $V_l$, $l=0,\frac 12, 1, \frac 32,\dots$} \end{matrix}
\Bigr \}
$$
so that, for $d \ge 0$,
$$
{\cal H}(d) \simeq
\Bigl \{ \begin{matrix} \text{linear span of matrix}  \\
\text{coefficients of $V_{\frac d2}$} \end{matrix} \Bigr \}
$$
and ${\cal H}^- \simeq {\cal H}^+$ so that
${\cal H}(-1) = 0$ and ${\cal H}(d) \simeq {\cal H} (-d-2) \cdot N(z)^{d+1}$
if $d < -1$.

As a representation of $SU(2) \times SU(2)$,
$$
{\cal H}(d) =
\{ \phi \in \BB C [z^0,z^1,z^2,z^3,N(Z)^{-1}]_d ;\: \square \phi=0 \}
\simeq
\begin{cases}
V_{\frac d2} \boxtimes V_{\frac d2} & \text{if $d \ge 0$;}  \\
0 & \text{if $d = -1$;}  \\
V_{-\frac d2-1} \boxtimes V_{-\frac d2-1} & \text{if $d < 0$.}
\end{cases}
$$
\end{prop}

Next we turn to spinor-valued Laurent polynomials
$$
\BB S \otimes_{\BB C} \BB C[z^0,z^1,z^2,z^3, N(Z)^{-1}]
\qquad \text{and} \qquad
\BB C[z^0,z^1,z^2,z^3, N(Z)^{-1}] \otimes_{\BB C} \BB S'.
$$
We fix the spaces
\begin{align*}
{\cal V} &= \{ \text{left-regular functions } f \in
\BB S \otimes_{\BB C} \BB C[z^0,z^1,z^2,z^3, N(Z)^{-1}] \}, \\
{\cal V}' &= \{ \text{right-regular functions } g \in
\BB C[z^0,z^1,z^2,z^3, N(Z)^{-1}] \otimes_{\BB C} \BB S'\}.
\end{align*}
These spaces have a natural grading
$$
{\cal V} = \bigoplus_{n \in \BB Z} {\cal V}(n)
\qquad \text{and} \qquad
{\cal V}' = \bigoplus_{n \in \BB Z} {\cal V}'(n)
$$
by the polynomial degree $n$.

By the Cauchy-Fueter formula (Corollary \ref{zero})
the pairing (\ref{fgpairing}) for $f \in {\cal V}$ and $g \in {\cal V}'$
is zero unless the open set $U$ contains the origin.
And if $U$ does contain the origin, the value of 
$\langle g, f \rangle_U$ is the same for all such sets $U$.
For concreteness, we choose $U$ to be the open unit ball and write
$\langle g, f \rangle$ for this $\langle g, f \rangle_U$, so
\begin{equation}  \label{vpairing}
\langle g, f \rangle = \int_{S^3} \langle g(X), Dx \cdot f(X) \rangle
= \int_{S^3} \langle g(X) \cdot Dx, f(X) \rangle \quad \in \BB C,
\qquad f \in {\cal V}, \: g \in {\cal V}'.
\end{equation}

\begin{prop}  \label{vpairing_prop}
The pairing (\ref{vpairing}) between ${\cal V}(n)$ and ${\cal V}'(m)$
is non-degenerate when $m+n+3=0$ and zero for all other $m$ and $n$.
\end{prop}

Restricting the actions (\ref{spinor_left_action}) and
(\ref{spinor_right_action}) to the diagonal matrices we obtain
$SU(2) \times SU(2)$ actions
\begin{equation}  \label{vaction}
\pi_l \begin{pmatrix} a & 0 \\ 0 & d \end{pmatrix}
: \: f(Z) \mapsto d \cdot f(a^{-1}Zd)
\qquad \text{and} \qquad
\pi_r \begin{pmatrix} a & 0 \\ 0 & d \end{pmatrix}
: \: g(Z) \mapsto g(a^{-1}Zd) \cdot a^{-1},
\end{equation}
where $a,d \in \BB H$, $|a|=|d|=1$, $f \in {\cal V}(n)$, $g \in {\cal V}'(m)$.

\begin{prop}
As representations of $SU(2) \times SU(2)$,
\begin{align*}
{\cal V}(n) & \simeq
\begin{cases}
V_{\frac n2} \boxtimes V_{\frac {n+1}2} & \text{if $n \ge 0$;}  \\
0 & \text{if $n = -1$ or $-2$;}  \\
V_{-\frac n2-1} \boxtimes V_{\frac {-n-3}2} & \text{if $n < -2$,}
\end{cases}
\qquad \dim_{\BB C} {\cal V}(n) = (n+2)(n+1),
\\
{\cal V}'(m) & \simeq
\begin{cases}
V_{\frac {m+1}2} \boxtimes V_{\frac m2} & \text{if $m \ge 0$;}  \\
0 & \text{if $m = -1$ or $-2$;}  \\
V_{\frac {-m-3}2} \boxtimes V_{-\frac m2-1} & \text{if $m < -2$,}
\end{cases}
\qquad \dim_{\BB C} {\cal V}'(m) = (m+2)(m+1).
\end{align*}
\end{prop}

\subsection{Bases for Regular Functions}  \label{bases_section}

Recall that we denote by $\partial_{ij}$
the partial derivatives $\frac{\partial}{\partial x_{ij}}$, $1 \le i,j \le 2$.
By differentiating the formula for matrix coefficients (\ref{int_m})
under the integral sign we obtain:

\begin{lem}  \label{deriv_calc}
\begin{align*}
\partial_{11} t^l_{n\,\underline{m}}(X) &=
(l-m) t^{l- \frac 12}_{n+ \frac 12 \, \underline{m+ \frac 12}}(X),  \\
\partial_{12} t^l_{n\,\underline{m}}(X) &=
(l+m) t^{l- \frac 12}_{n+ \frac 12 \, \underline{m- \frac 12}}(X),  \\
\partial_{21} t^l_{n\,\underline{m}}(X) &=
(l-m) t^{l- \frac 12}_{n- \frac 12 \, \underline{m+ \frac 12}}(X),  \\
\partial_{22} t^l_{n\,\underline{m}}(X) &=
(l+m) t^{l- \frac 12}_{n- \frac 12 \, \underline{m- \frac 12}}(X).
\end{align*}
\end{lem}
In the above formulas the parameters $l'$, $m'$, $n'$
in some matrix coefficients $t^{l'}_{m'\,\underline{n'}}(X)$ may be out of
allowed range $l'=0,\frac 12, 1, \frac 32,\dots$ and $|m'|,|n'| \le l'$.
If that happens we declare such coefficients equal zero. 

We also have the following multiplication identities for matrix coefficients.

\begin{lem}  \label{mult-ids}
$$
\Bigl( t^{l-\frac 12}_{m + \frac 12 \, \underline{n}}(X),
t^{l-\frac 12}_{m - \frac 12 \, \underline{n}}(X) \Bigr)
\begin{pmatrix} x_{11} & x_{12} \\ x_{21} & x_{22} \end{pmatrix}
=
\Bigl( t^{l}_{m \, \underline{n- \frac 12}}(X),
t^{l}_{m \, \underline{n+ \frac 12}}(X) \Bigr)
$$
and
$$
\begin{pmatrix} x_{11} & x_{12} \\ x_{21} & x_{22} \end{pmatrix}
\begin{pmatrix}
(l-m+ \frac 12) t^l_{n \, \underline{m+ \frac 12}}(X)  \\
(l+m+ \frac 12) t^l_{n \, \underline{m- \frac 12}}(X)
\end{pmatrix}
=
\begin{pmatrix}
(l-n+1) t^{l + \frac 12}_{n - \frac 12 \, \underline{m}}(X)  \\
(l+n+1) t^{l + \frac 12}_{n + \frac 12 \, \underline{m}}(X)
\end{pmatrix}
$$
(in the sense of matrix multiplication).
\end{lem}

\begin{prop}  \label{dual_bases}
The following $\BB S$-valued functions form a basis for ${\cal V}(2l)$ --
the left-regular polynomial functions on $\BB H$ of degree $2l$,
$l=0,\frac 12, 1, \frac 32,\dots$:
$$
\begin{pmatrix}
(l-m+ \frac 12) t^l_{n \, \underline{m+ \frac 12}}(X)  \\
(l+m+ \frac 12) t^l_{n \, \underline{m- \frac 12}}(X)
\end{pmatrix},
\qquad
\begin{matrix}
m =-l-\frac 12 ,-l+\frac 32,\dots,l+\frac 12,  \\
n =-l,-l+1,\dots,l.
\end{matrix}
$$
The functions
$$
\begin{pmatrix}
(l-n+ \frac 12) t^l_{n- \frac 12 \, \underline{m}}(X^{-1}) \cdot N(X)^{-1}  \\
(l+n+ \frac 12) t^l_{n+ \frac 12 \, \underline{m}}(X^{-1}) \cdot N(X)^{-1}
\end{pmatrix},
\qquad
\begin{matrix}
m=-l,-l+1,\dots,l,  \\
n =-l+\frac 12 ,-l+\frac 32,\dots,l-\frac 12,
\end{matrix}
$$
form a basis of ${\cal V}(-2l-2)$ --
the left-regular polynomial functions on $\BB H^{\times}$ of degree $(-2l-2)$.
Altogether these functions form a basis of
${\cal V} = \bigoplus_{k \in \BB Z} {\cal V}(k)$.

The basis dual with respect to the pairing (\ref{vpairing}) consists of
$\BB S'$-valued right-regular functions
$$
\Bigl(
t^{l+\frac 12}_{m \, \underline{n- \frac 12}}(X^{-1}) \cdot N(X)^{-1},
t^{l+\frac 12}_{m \, \underline{n+ \frac 12}}(X^{-1}) \cdot N(X)^{-1}
\Bigr),
\qquad
\begin{matrix}
m =-l-\frac 12 ,-l+\frac 32,\dots,l+\frac 12,  \\
n =-l,-l+1,\dots,l,
\end{matrix}
$$
which, for $l$ fixed, span ${\cal V}'(-2l-3)$, and
$$
\Bigl(
t^{l-\frac 12}_{m + \frac 12 \, \underline{n}}(X),
t^{l-\frac 12}_{m - \frac 12 \, \underline{n}}(X)
\Bigr),
\qquad
\begin{matrix}
m=-l,-l+1,\dots,l,  \\
n =-l+\frac 12 ,-l+\frac 32,\dots,l-\frac 12,
\end{matrix}
$$
which, for $l$ fixed, span ${\cal V}'(2l-1)$.
Altogether these functions form a basis of
${\cal V}' = \bigoplus_{k \in \BB Z} {\cal V}'(k)$.
\end{prop}

Next we prove a quaternionic analogue of the complex variable formula
$$
\frac 1{x-y} =
\begin{cases}
\frac 1x \sum_{k=0}^{\infty} x^{-k}y^k,
& \text{valid for $x,y \in \BB C$, $|y| < |x|$;}  \\
\quad  \\
-\frac 1y \sum_{k=0}^{\infty} x^ky^{-k},
& \text{valid for $x,y \in \BB C$, $|y| > |x|$.}
\end{cases}
$$

\begin{prop}  \label{k0_expansion}
We have the following matrix coefficient expansions
$$
k_0(X-Y) = \frac 1{N(X-Y)} = N(X)^{-1} \cdot \sum_{l,m,n}
t^l_{n \, \underline{m}}(X^{-1}) \cdot t^l_{m \, \underline{n}}(Y)
$$
which converges pointwise absolutely in the region
$\{ (X,Y) \in \BB H^{\times} \times \BB H^{\times}; |Y|<|X| \}$, and
$$
k_0(X-Y) = N(Y)^{-1} \cdot \sum_{l,m,n} 
t^l_{n \, \underline{m}}(X) \cdot t^l_{m \, \underline{n}}(Y^{-1})
$$
which converges pointwise absolutely in the region
$\{ (X,Y) \in \BB H^{\times} \times \BB H^{\times}; |X|<|Y| \}$.
The sums are taken first over all $m,n = -l, -l+1, \dots, l$, then
over $l=0,\frac 12, 1, \frac 32,\dots$.
\end{prop}

\pf
By the multiplicativity property of matrix coefficients,
$$
\sum_{m,n=-l, -l+1, \dots, l} 
t^l_{n \, \underline{m}}(X^{-1}) \cdot t^l_{m \, \underline{n}}(Y)
= \sum_{n=-l, -l+1, \dots, l} t^l_{n \, \underline{n}}(X^{-1}Y)
= \tr \bigl( \tau_l (X^{-1}Y) \bigr).
$$
Now, any element $X^{-1}Y \in \BB H$ can be diagonalized and it is known that
\begin{equation}  \label{trace}
\tr \biggl(
\tau_l \begin{pmatrix} \lambda_1 & 0 \\ 0 & \lambda_2 \end{pmatrix} \biggr) =
\frac {\lambda_1^{2l+1} - \lambda_2^{2l+1}}{\lambda_1 - \lambda_2}.
\end{equation}
Although this expression is valid only on the region where
$\lambda_1 \ne \lambda_2$, the right hand side clearly continues
analytically across the set of $X^{-1}Y$ for which $\lambda_1=\lambda_2$.
Now, an element $h \in \BB H \subset \mathfrak{gl}(2,\BB C)$ has both
eigenvalues of length less than one if and only if $N(h)<1$.
Thus, in the region $\{|Y|<|X|\}$, we have $|\lambda_1|,|\lambda_2|<1$
and summing the terms in (\ref{trace}) over $l=0,\frac 12, 1, \frac 32,\dots$
we get an absolutely convergent geometric series:
\begin{multline*}
\sum_{l=0,\frac 12, 1, \frac 32,\dots} 
\tr \biggl(
\tau_l \begin{pmatrix} \lambda_1 & 0 \\ 0 & \lambda_2 \end{pmatrix} \biggr)
= \frac 1{\lambda_1 - \lambda_2} \sum_{l=0,\frac 12, 1, \frac 32,\dots} 
(\lambda_1^{2l+1} - \lambda_2^{2l+1})  \\
= \frac 1{(1-\lambda_1)(1-\lambda_2)}
= \frac {N(X)}{N(X-Y)}.
\end{multline*}
The other expansion follows from the first one by switching the
variables $X$ and $Y$ and the indices $m$ and $n$.
\qed

Applying the differential operator $\nabla$ to both sides of the two
expansions we obtain matrix coefficient expansions for
$\frac {(X-Y)^{-1}}{N(X-Y)}$:

\begin{prop}  \label{k-expansion}
We have the following matrix coefficient expansions
\begin{multline*}
k(X-Y) = \frac {(X-Y)^{-1}}{N(X-Y)}  \\
= \sum_{l,m,n}
\begin{pmatrix}
(l-m+ \frac 12) t^l_{n \, \underline{m+ \frac 12}}(Y)  \\
(l+m+ \frac 12) t^l_{n \, \underline{m- \frac 12}}(Y)
\end{pmatrix}
\cdot
\Bigl(
t^{l+\frac 12}_{m \, \underline{n- \frac 12}}(X^{-1}) \cdot N(X)^{-1},
t^{l+\frac 12}_{m \, \underline{n+ \frac 12}}(X^{-1}) \cdot N(X)^{-1}
\Bigr),
\end{multline*}
which converges pointwise absolutely in the region
$\{ (X,Y) \in \BB H^{\times} \times \BB H^{\times}; |Y|<|X| \}$, and
the sum is taken first over all
$m =-l-\frac 12 ,-l+\frac 32,\dots,l+\frac 12$ and $n =-l,-l+1,\dots,l$,
then over $l=0,\frac 12, 1, \frac 32,\dots$.
In other words, we sum over the basis elements of
${\cal V}^+ = \bigoplus_{k \ge 0} {\cal V}(k)$ given in
Proposition \ref{dual_bases} multiplied by their respective dual basis
elements of ${\cal V}'^- = \bigoplus_{k < 0} {\cal V}'(k)$.
Similarly,
$$
k(X-Y) = - \sum_{l,m,n}
\begin{pmatrix}
(l-n+ \frac 12) t^l_{n- \frac 12 \, \underline{m}}(Y^{-1}) \cdot N(Y)^{-1}  \\
(l+n+ \frac 12) t^l_{n+ \frac 12 \, \underline{m}}(Y^{-1}) \cdot N(Y)^{-1}
\end{pmatrix}
\cdot
\Bigl(
t^{l-\frac 12}_{m + \frac 12 \, \underline{n}}(X),
t^{l-\frac 12}_{m - \frac 12 \, \underline{n}}(X)
\Bigr),
$$
which converges pointwise absolutely in the region
$\{ (X,Y) \in \BB H^{\times} \times \BB H^{\times}; |X|<|Y| \}$,
and the sum is taken first over all
$m =-l,-l+1,\dots,l$ and $n =-l+\frac 12 ,-l+\frac 32,\dots,l-\frac 12$,
then over $l =0,\frac 12, 1, \frac 32,\dots$.
In other words, we sum over the basis elements of
${\cal V}^- = \bigoplus_{k < 0} {\cal V}(k)$ given in
Proposition \ref{dual_bases} multiplied by their respective dual basis
elements of ${\cal V}'^+ = \bigoplus_{k \ge 0} {\cal V}'(k)$.
\end{prop}

This matrix coefficient expansion can be used to give an alternative proof
of the Cauchy-Fueter formulas (Theorem \ref{Fueter})
using the duality of the bases.

We also have the following matrix coefficient expansions. Their proofs are
similar to those for $\frac 1{N(X-Y)}$.

\begin{prop}  \label{N-square-expansion}
We have the following matrix coefficient expansions
$$
\frac 1{N(X-Y)^2} = \sum_{k,l,m,n}
(2l+1) t^l_{n \, \underline{m}}(X^{-1}) \cdot N(X)^{-k-2} \cdot
t^l_{m \, \underline{n}}(Y) \cdot N(Y)^k
$$
which converges pointwise absolutely in the region
$\{ (X,Y) \in \BB H^{\times} \times \BB H^{\times}; |Y|<|X| \}$, and
$$
\frac 1{N(X-Y)^2} = \sum_{k,l,m,n} 
(2l+1) t^l_{n \, \underline{m}}(X) \cdot N(X)^k \cdot
t^l_{m \, \underline{n}}(Y^{-1}) \cdot N(Y)^{-k-2}
$$
which converges pointwise absolutely in the region
$\{ (X,Y) \in \BB H^{\times} \times \BB H^{\times}; |X|<|Y| \}$.
The sums are taken first over all $m,n = -l, -l+1, \dots, l$,
then over $k=0,1,2,3,\dots$ and $l=0,\frac 12, 1, \frac 32,\dots$.
\end{prop}

\subsection{The Spaces ${\cal H}$, ${\cal V}$ and ${\cal V}'$ as
Representations of $\mathfrak{gl}(2,\BB H)$}

We decompose the spaces ${\cal H}$, ${\cal V}$ and ${\cal V}'$ into their
``positive'' and ``negative'' parts:
\begin{center}
\begin{tabular}{ccccc}
${\cal H} = {\cal H}^+ \oplus {\cal H}^-$, & \qquad &
${\cal H}^+ = \bigoplus_{n \in \BB Z,\: n \ge 0} {\cal H}(n)$, & \qquad &
${\cal H}^- = \bigoplus_{n \in \BB Z,\: n < 0} {\cal H}(n)$,  \\
${\cal V} = {\cal V}^+ \oplus {\cal V}^-$, & \qquad &
${\cal V}^+ = \bigoplus_{n \in \BB Z,\: n \ge 0} {\cal V}(n)$, & \qquad &
${\cal V}^- = \bigoplus_{n \in \BB Z,\: n < 0} {\cal V}(n)$,  \\
${\cal V}' = {\cal V}'^+ \oplus {\cal V}'^-$, & \qquad &
${\cal V}'^+ = \bigoplus_{n \in \BB Z,\: n \ge 0} {\cal V}'(n)$, & \qquad &
${\cal V}'^- = \bigoplus_{n \in \BB Z,\: n < 0} {\cal V}'(n)$.
\end{tabular}
\end{center}

We write $\partial_i$ for $\frac{\partial}{\partial x^i}$, $0 \le i \le 3$,
and introduce two differential operators acting on functions on $\BB H$ --
the degree operator and the degree operator plus identity:
$$
\deg = x^0\partial_0+x^1\partial_1+x^2\partial_2+x^3\partial_3,
\qquad
\widetilde{\deg}=1+\deg.
$$

\begin{thm}  \label{unitary-u(2,2)-action}
\begin{enumerate}
\item
The spaces ${\cal H}^+$, ${\cal H}^-$, ${\cal V}^+$, ${\cal V}^-$,
${\cal V}'^+$ and ${\cal V}'^-$ are irreducible representations of
$\BB C \otimes \mathfrak{gl}(2,\BB H) \simeq \mathfrak{gl}(4,\BB C)$
(and $\BB C \otimes \mathfrak{sl}(2,\BB H) \simeq \mathfrak{sl}(4,\BB C)$).
\item
Each of these spaces has a unitary structure such that the real form
$\mathfrak u(2,2)$ of $\mathfrak{gl}(4,\BB C)$ acts unitarily.
Explicitly, the unitary structure on ${\cal V}^{\pm}$ is
$$
(f_1,f_2) = \int_{S^3} \langle f_2^+(X), f_1(X) \rangle \,dS,
$$
similarly for ${\cal V}'^{\pm}$,
and the unitary structure on ${\cal H}^{\pm}$ is
$$
(\phi_1,\phi_2) =
\int_{S^3} (\widetilde{\deg} \phi_1)(X) \cdot \B{\phi_2}(X) \,dS.
$$
\item
There is a non-degenerate $\mathfrak{u}(2,2)$-invariant bilinear
pairing between $(\pi^0_l, {\cal H}^+)$ and $(\pi^0_r, {\cal H}^-)$
given by
\begin{equation}  \label{H-pairing}
\langle \phi_1,\phi_2 \rangle =
\frac 1{2\pi^2} \int_{S_R} (\widetilde{\deg} \phi_1)(X) \cdot \phi_2(X)
\,\frac{dS}R,
\qquad \phi_1 \in {\cal H}^+, \: \phi_2 \in {\cal H}^-,
\end{equation}
which is independent of the choice of $R>0$.
\end{enumerate}
\end{thm}

\begin{rem}
Note the appearance in this theorem of
$\mathfrak{u}(2,2) \simeq \mathfrak{so}(4,2)$ -- the Lie algebra
of the conformal group of the Minkowski space $\BB M$.
This demonstrates the necessity of complexification of 
$\mathfrak{gl}(2,\BB H)$ and $\BB H$ and passing to different real forms
such as $\mathfrak{u}(2,2)$ and $\BB M$ which will be done in the next chapter.
\end{rem}

The representations ${\cal H}^+$, ${\cal V}^+$ and ${\cal V}'^+$ belong to
the most degenerate series of unitary representations of $\mathfrak{su}(2,2)$
with highest weights. This series is parameterized by spin, which runs over
all semi-integers. In particular, ${\cal H}^+$, ${\cal V}^+$ and ${\cal V}'^+$
have spins $0$, $1/2$ and $1/2$ respectively.

We denote by $\operatorname{Map} (S^3, \BB C)$, 
$\operatorname{Map} (S^3, \BB S)$ and
$\operatorname{Map} (S^3, \BB S')$
the spaces of smooth functions on $S^3$ with values in $\BB C$,
$\BB S$ and $\BB S'$ respectively.
The group $SU(2) \times SU(2)$ acts naturally on each of these spaces
as diagonal elements of $GL(2,\BB H)$.
We define the spaces of ``polynomial functions'' on $S^3$ with values in
$\BB C$, $\BB S$ and $\BB S'$ as the functions in respectively
$\operatorname{Map} (S^3, \BB C)$, $\operatorname{Map} (S^3, \BB S)$
and $\operatorname{Map} (S^3, \BB S')$ which lie in a finite-dimensional
subspace invariant under the respective $SU(2) \times SU(2)$-action.
We denote these spaces respectively by
$$
\pol(S^3), \qquad \pol(S^3, \BB S), \qquad \pol(S^3, \BB S').
$$

\begin{prop}
Restricting functions defined on $\BB H^{\times}$ to the sphere $S^3$
induces the following isomorphisms:
$$
{\cal H}^+ \simeq \pol(S^3) \simeq {\cal H}^-,
\qquad
{\cal V} \simeq \pol(S^3, \BB S)
\qquad \text{and} \qquad
{\cal V}' \simeq \pol(S^3, \BB S').
$$
Moreover, ${\cal H}^+$ is dense in $L^2(S^3)$, and
${\cal V}$ (respectively ${\cal V}'$) is dense in $L^2(S^3, \BB S)$
(respectively $L^2(S^3, \BB S')$).
\end{prop}

The first isomorphism in the proposition can be interpreted as
``every function on $S^3$ can be extended to a harmonic function on
$\BB H^{\times}$ in exactly two ways: so that it is regular at $0$
and so that it is regular at $\infty$.''
Similarly, the other two isomorphisms can be interpreted as
``every function on $S^3$ with values in $\BB S$ (respectively $\BB S'$)
can be extended to a left-regular (respectively right-regular) function on
$\BB H^{\times}$ in a unique way.''

\subsection{Poisson Formula}

In this subsection we prove a Poisson-type formula for harmonic
functions defined on a ball in $\BB H$.
This formula is a special case $n=4$ of a general formula that expresses
harmonic functions in the interior of a ball in $\BB R^n$
as an integral over the boundary sphere (see, for instance, \cite{SW}).
However, our method uses representation theory and can be applied to
prove other results such as Theorems \ref{Poisson-M-thm}, \ref{2pole1},
\ref{2pole1-M}.
As before, $S^3_R$ denotes the sphere of radius $R$,
i.e. $S^3_R = \{X \in \BB H ;\: |X|=R \}$.
Recall the conformal action of $GL(2,\BB H)$ on $\widehat{\BB H}$:
$$
X \mapsto \tilde X = (aX+b)(cX+d)^{-1}, \qquad
h^{-1} = \begin{pmatrix} a & b \\ c & d \end{pmatrix} \in GL(2, \BB H).
$$

\begin{lem}  \label{N=1-preserving}
Let $G \subset GL(2, \BB H)$ be the subgroup consisting of all elements
of $GL(2,\BB H)$ preserving the unit sphere. Then
$$
Lie(G) = \biggl\{ \begin{pmatrix} A & B \\ B^+ & D \end{pmatrix} ;\:
A,B,D \in \BB H, \: \re A = \re D \biggr\}.
$$
\end{lem}

\pf
The Lie algebra of $G$ consists of all matrices
$\begin{pmatrix} A & B \\ C & D \end{pmatrix}$, $A,B,C,D \in \BB H$,
which generate vector fields tangent to $S^3_1$.
Such a matrix $\begin{pmatrix} A & B \\ C & D \end{pmatrix}$ generates a vector
field
$$
\frac d{dt} \bigl( (1+tA)X+tB \bigr) (tCX+1+tD)^{-1} \Bigr|_{t=0}
= AX+B-XCX-XD.
$$
A vector field is tangent to $S^3_1$ if and only if it is orthogonal
with respect to (\ref{bilinear}) to the vector field $X$ for $|X|=1$:
$$
0=\re \bigl( (AX+B-XCX-XD)X^+ \bigr) = \re (A-D+BX^+-XC),
\qquad \forall X \in S^3_1.
$$
It follows that $\re A = \re D$ and $C=B^+$, hence the result.
\qed

\begin{cor}
Let $G_0 \subset GL(2, \BB H)$ be the connected subgroup with Lie algebra
$$
\mathfrak g_0 = \{ x \in Lie(G);\: \re(\tr x)=0 \}
= \biggl\{ \begin{pmatrix} A & B \\ B^+ & D \end{pmatrix} ;\:
A,B,D \in \BB H, \: \re A = \re D =0 \biggr\}.
$$
Then $G_0$ preserves the unit sphere $S^3_1 = \{ X \in \BB H ;\: |X|=1 \}$,
the open ball $\{ X \in \BB H ;\: |X|<1 \}$ and the open set
$\{ X \in \BB H ;\: |X|>1 \}$.
\end{cor}

The Lie algebra $\mathfrak g_0$ and the Lie group $G_0$ are isomorphic
to $\mathfrak{sp}(1,1)$ and $Sp(1,1)$ respectively
(see, for example, \cite{H}).

\begin{lem}
The group $G_0$ is generated by $SU(2) \times SU(2)$ realized as
the subgroup of diagonal matrices 
$\begin{pmatrix} a & 0 \\ 0 & d \end{pmatrix}$, $a,d \in \BB H$,
$|a|=|d|=1$, and the one-parameter group
$$
G_0' = \biggl\{ \begin{pmatrix}
\cosh t & \sinh t \\ \sinh t & \cosh t \end{pmatrix} ;\:
t \in \BB R \biggr\}.
$$
\end{lem}

\pf
The subgroup of diagonal matrices
$\begin{pmatrix} a & 0 \\ 0 & d \end{pmatrix}$
with $|a|=|d|=1$ has Lie algebra 
$$
\mathfrak{su}(2) \oplus \mathfrak{su}(2) =
\biggl\{ \begin{pmatrix} A & 0 \\ 0 & D \end{pmatrix} ;\:
A,D \in \BB H, \: \re A = \re D =0 \biggr\}.
$$
The one-parameter group $G_0'$ has Lie algebra
$$
\mathfrak{g}_0' = \biggl\{ \begin{pmatrix} 0 & t \\ t & 0 \end{pmatrix} ;\:
t \in \BB R \biggr\}.
$$
An easy computation shows that $\mathfrak{su}(2) \oplus \mathfrak{su}(2)$
and $\mathfrak{g}_0'$ generate $\mathfrak{g}_0$. Since the groups
$SU(2) \times SU(2)$, $G_0'$ and $G_0$ are connected, the result follows.
\qed

We are now ready to state the main result of this subsection.

\begin{thm} \label{Poisson}
Let $\phi$ be a real analytic harmonic function defined on a closed ball
$\{X \in \BB H ;\: |X| \le R \}$, for some $R>0$. Then
\begin{align*}
\phi(Y) &= \frac 1{2\pi^2} \int_{X \in S^3_R}
\frac {R^2 - N(Y)}{N(X-Y)^2} \cdot \frac{dS}R \cdot \phi(X),  \\
&= - \frac 1{2\pi^2} \int_{X \in S^3_R}
\Bigl( \widetilde{\deg}_X \frac 1{N(X-Y)} \Bigr)
\cdot \frac{dS}R \cdot \phi(X),  \\
&= \frac 1{2\pi^2} \int_{X \in S^3_R} \frac 1{N(X-Y)} \cdot \frac{dS}R
\cdot (\widetilde{\deg}\phi)(X),
\qquad \forall Y \in \BB H, \: |Y|<R.
\end{align*}
\end{thm}

By direct computation we have
$$
\widetilde{\deg}_X \frac 1{N(X-Y)} = - \frac {N(X) - N(Y)}{N(X-Y)^2}.
$$
Since $N(X)=R^2$, we see that the first two formulas are equivalent, and
the function $\frac {R^2 - N(Y)}{N(X-Y)^2}$ is harmonic with respect to
the $Y$ variable.
First we consider the special case $R=1$.

\begin{prop}  \label{Poisson_R=1}
Let $\phi$ be a harmonic function defined on a closed ball
$\{X \in \BB H ;\: |X| \le 1 \}$. Then
\begin{equation*}
\phi(Y) = \frac 1{2\pi^2} \int_{X \in S^3_1}
\frac {1-N(Y)}{N(X-Y)^2} \cdot \phi(X) \,dS,
\qquad \forall Y \in \BB H, \: |Y|<1.
\end{equation*}
\end{prop}

\pf
We define a map
$$
\phi \mapsto \operatorname{Harm}(\phi), \qquad
(\operatorname{Harm}\phi)(Y) = \frac 1{2\pi^2} \int_{X \in S^3_1}
\frac {1 - N(Y)}{N(X-Y)^2} \cdot \phi(X) \,dS.
$$
We want to show that $\operatorname{Harm}$ is the identity mapping.
Let $\widehat{\cal H}$ denote the space of real analytic harmonic functions
on the closed ball $\{ X \in \BB H ;\: |X| \le 1 \}$.

\begin{lem}  \label{G_0-equivariant}
The map $\operatorname{Harm}: \widehat{\cal H} \to \widehat{\cal H}$
is equivariant with respect to the $\pi_l^0$-action of $G_0$ on
$\widehat{\cal H}$.
\end{lem}

\pf
The map
$\operatorname{Harm}: \widehat{\cal H} \to \widehat{\cal H}$
is equivariant with respect to the $\pi_l^0$ action of $SU(2) \times SU(2)$
on $\widehat{\cal H}$.
Since the group $G_0$ is generated by $SU(2) \times SU(2)$ and
the one-parameter group
$$
G_0' = \biggl\{ \begin{pmatrix}
\cosh t & \sinh t \\ \sinh t & \cosh t \end{pmatrix} ;\:
t \in \BB R \biggr\},
$$
it is sufficient to show that $\operatorname{Harm}$ is $G_0'$-equivariant.
We want to compare $(\operatorname{Harm}\tilde\phi)(Y)$ with
$(\widetilde {\operatorname{Harm}\phi})(Y)$,
where $\tilde\phi = \pi^0_l(g) \phi$ and
$(\widetilde {\operatorname{Harm}\phi})
= \pi^0_l(g) (\operatorname{Harm} (\phi))$,
$g \in G_0$. A straightforward computation implies

\begin{lem}  \label{Jacobian-lemma}
Fix an element
$g = \begin{pmatrix} \cosh t & \sinh t \\ \sinh t & \cosh t \end{pmatrix}
\in G_0'$ and consider its conformal action on $S^3_1$:
$$
\pi_l(g):\: X \mapsto \tilde X =
(\cosh t X - \sinh t)(-\sinh t X + \cosh t)^{-1}.
$$
Then the Jacobian $J$ of this map is
$$
\pi_l(g)^* (dS) = J\,dS =\frac 1{N(-\sinh t X + \cosh t)} \cdot
\frac{1 - (\re \tilde X)^2}{1 - (\re X)^2} \,dS.
$$
\end{lem}

Since $G_0'$ is connected, it is sufficient to verify the $G_0'$-equivariance
of $\operatorname{Harm}$ on the Lie algebra level.
For $t \to 0$ and modulo terms of order $t^2$, we have:
$$
\tilde X = (X -t)(1-tX)^{-1}=X+t(X^2-1),
\qquad \tilde Y = Y+t(Y^2-1),
$$
$$
N(-\sinh t X + \cosh t) = 1 - 2t \re X, \qquad
N(-\sinh t Y + \cosh t) = 1 - 2t \re Y,
$$
$$
N(\tilde Y) = N(Y) + 2t \re (Y^+(Y^2-1)) = N(Y) - 2t(1-N(Y)) \re Y.
$$
Since $\re(X^2)=2(\re X)^2-N(X)=2(\re X)^2-1$,
$$
1 - (\re \tilde X)^2 = 1 - (\re X + t\re(X^2) -t)^2
= (1 -(\re X)^2)(1+4t \re X).
$$
Continuing to work modulo terms of order $t^2$, we get
\begin{multline*}
(\operatorname{Harm}\tilde\phi)(Y)
= \frac 1{2\pi^2} \int_{X \in S^3_1} \frac {1 - N(Y)}{N(X-Y)^2} \cdot
\frac 1{N(-\sinh t X + \cosh t)} \cdot \phi(\tilde X) \,dS  \\
= \frac 1{2\pi^2} \int_{X \in S^3_1} \frac {1 - N(Y)}{N(X-Y)^2}
\cdot (1 + 2t \re X) \cdot \phi(\tilde X) \,dS.
\end{multline*}
On the other hand, using Lemmas \ref{X-Y_difference} and \ref{Jacobian-lemma},
\begin{multline*}
(\widetilde {\operatorname{Harm}\phi})(Y) =
\frac 1{N(-\sinh t Y + \cosh t)} \cdot (\operatorname{Harm}\phi)(\tilde Y)  \\
= \frac 1{N(-\sinh t Y + \cosh t)} \cdot
\frac 1{2\pi^2} \int_{\tilde X \in S^3_1}
\frac {1 - N(\tilde Y)}{N(\tilde X - \tilde Y)^2} \cdot \phi( \tilde X) \,dS \\
= \frac 1{2\pi^2} \int_{X \in S^3_1} \frac {1 - N(Y)}{N(X-Y)^2}
\cdot (1 + 2t \re X) \cdot \phi(\tilde X) \,dS.
\end{multline*}
This proves that the map
$\operatorname{Harm}: \widehat{\cal H} \to \widehat{\cal H}$
is $G_0$-equivariant.
\qed

We can now finish the proof of Proposition \ref{Poisson_R=1}.
It is easy to show that the space $\widehat{\cal H}$ viewed as a
representation of $G_0$ is irreducible.
The group $SU(2)\times SU(2)$ is a maximal compact
subgroup of $G_0$, and the map $\operatorname{Harm}$ must preserve
the space of $SU(2)\times SU(2)$-finite vectors ${\cal H}$.
Hence by Schur's Lemma applied to the irreducible
$(SU(2)\times SU(2), \mathfrak{g}_0)$-module ${\cal H}$,
there exists a $\lambda \in \BB C$ such that the map 
$\operatorname{Harm}: \cal H \to \cal H$ is given by multiplication
by $\lambda$.
To pin down the value of $\lambda$ we substitute $\phi(X)=1$ and $Y=0$,
and we immediately see that $\lambda=1$. Since ${\cal H}$ is dense in
$\widehat{\cal H}$, this finishes our proof of Proposition \ref{Poisson_R=1}.
\qed

\noindent {\it Proof of Theorem \ref{Poisson}.}
Given a harmonic function $\phi$ on the closed ball
$\{X \in \BB H ;\: |X| \le R \}$, we can consider a function
$\phi'(X)=\phi(RX)$, then by Proposition \ref{Poisson_R=1},
$\operatorname{Harm}(\phi') = \phi'$, and the first integral
formula follows.

To prove the last integral formula we consider the following
integral operator on $\widehat{\cal H}$:
$$
\phi \mapsto \operatorname{Harm}'(\phi), \qquad
(\operatorname{Harm}'\phi)(Y) =
\int_{X \in S^3_R} \frac 1{N(X-Y)} \cdot \frac{dS}R \cdot \phi(X).
$$
Since
$$
\widetilde{\deg}_X \frac 1{N(X-Y)} = - \widetilde{\deg}_Y \frac 1{N(X-Y)},
$$
we see that
$$
\operatorname{Harm} = \widetilde{\deg}_Y \circ \operatorname{Harm}'.
$$
But the differential operator $\widetilde{\deg}$ is injective on
$\widehat{\cal H}$, hence
$\operatorname{Harm}' \circ \widetilde{\deg}_X$
is the identity operator on $\widehat{\cal H}$.
This proves the last integral formula.
\qed

\begin{rem}
One can give at least two alternative proofs of Theorem \ref{Poisson}:
One proof -- using the matrix coefficient expansion of $\frac 1{N(X-Y)}$
given in Proposition \ref{k0_expansion} and another --
using the equivariance of the map
$\operatorname{Harm}: \widehat{\cal H} \to \widehat{\cal H}$
with respect to $SU(2) \times SU(2)$ only, thus avoiding the computational
part of showing the $G_0$-equivariance.
However, these proofs do not generalize so well to other cases like
the case of the Minkowski space $\BB M$ that will be discussed in Subsection
\ref{Poisson-M}.
\end{rem}

Changing the variables $X \mapsto X^{-1}$ we obtain

\begin{cor}  \label{Poisson 1/X}
Let $\phi$ be a harmonic function defined on a closed set
$\{X \in \BB H ;\: |X| \ge R \}$, for some $R>0$, and regular at infinity.
Then
$$
\phi(Y) = -\frac 1{2\pi^2} \int_{X \in S^3_R} \frac 1{N(X-Y)} \cdot \frac{dS}R
\cdot (\widetilde{\deg}\phi)(X),
\qquad \forall Y \in \BB H, \: |Y| >R.
$$
\end{cor}

\subsection{Hydrogen Atom: the Discrete Spectrum}

In this subsection we give an example of an application of
quaternionic analysis to physics. We recast the construction of \cite{BI}(I)
of the discrete part of the spectral decomposition of the three-dimensional
Laplacian with the Coulomb potential, starting with the Poisson formula.

Let $\phi$ be a harmonic function which is homogeneous of degree $2l$, then
$$
(\widetilde{\deg}\phi)(X) = (2l+1) \phi(X).
$$
In this case the Poisson formula (Theorem \ref{Poisson}) for $R=1$ yields
\begin{equation}  \label{Poisson-hydro}
\phi(Y) = \frac {2l+1}{2\pi^2} \int_{X \in S^3} \frac {\phi(X)}{N(X-Y)} \,{dS},
\qquad \forall Y \in \BB H, \: |Y|<1,
\end{equation}
and $\phi$ is a linear combination of the matrix coefficients
$t^l_{n\, \underline{m}}$'s which span a vector space of polynomials
of dimension $(2l+1)^2$.

Fix a $\rho >0$ and apply the Cayley transform
\begin{equation}  \label{Cayley}
\pi_l \begin{pmatrix} \rho & \rho \\ -1 & 1 \end{pmatrix}: \BB H \to \BB H,
\qquad Z \mapsto X= \frac {Z-\rho}{Z+\rho}.
\end{equation}
It has the following properties:
$$
N(X)=1 \Longleftrightarrow \re Z =0, \qquad
N(X)<1 \Longleftrightarrow \re Z >0.
$$
Thus we can allow in (\ref{Poisson-hydro})
$$
X=\frac {Z-\rho}{Z+\rho}, \quad \re Z =0
\qquad \text{and} \qquad 
Y=\frac {W-\rho}{W+\rho}, \quad \re W >0.
$$
From Lemma \ref{restriction} we know that $dS = \frac {Dx}X$,
hence from Proposition \ref{Dx-pullback} and Lemma \ref{X-Y_difference}
we obtain:

\begin{lem}
The pull-back of the measure $dS$ on $S^3$ by the Cayley transform
(\ref{Cayley}) is
$$
\frac {8\rho^3}{(N(Z)+\rho^2)^3} \,dS(Z),
$$
where $dS(Z)$ denotes the Euclidean measure on $\{ \re Z =0 \}$. Also
$$
N(X-Y) = \frac {4\rho^2 \cdot N(Z-W)}{(N(Z)+\rho^2)(N(W)+\rho^2)}.
$$
\end{lem}

Thus we can rewrite (\ref{Poisson-hydro}) in terms of $Z$ and $W$:
\begin{equation}  \label{Poisson-intermediate1}
\phi \biggl( \frac {W-\rho}{W+\rho} \biggr) = \frac {(2l+1)\rho}{\pi^2}
\int_{\re Z=0}\phi \biggl( \frac {Z-\rho}{Z+\rho} \biggr) \cdot 
\frac {N(W)+\rho^2}{(N(Z)+\rho^2)^2} \cdot \frac {dS(Z)}{N(Z-W)},
\end{equation}
for all $W \in \BB H$ with $\re W>0$.
Now we introduce
$$
\psi(Z) = \frac 1{(N(Z)+\rho^2)^2} \cdot
\phi \biggl( \frac {Z-\rho}{Z+\rho} \biggr),
\qquad
\psi(W) = \frac 1{(N(W)+\rho^2)^2} \cdot
\phi \biggl( \frac {W-\rho}{W+\rho} \biggr).
$$
Then we can rewrite (\ref{Poisson-intermediate1}) as
\begin{equation}  \label{Poisson-intermediate2}
(N(W)+\rho^2) \cdot \psi(W) = \frac {(2l+1)\rho}{\pi^2} \int_{\re Z=0}
\psi(Z) \cdot \frac{dS(Z)}{N(Z-W)},
\qquad \forall W \in \BB H, \: \re W>0.
\end{equation}

Next we consider the Fourier transforms over the three-dimensional planes
in $\BB H$ parallel to $\{\re Z=0\}$:
$$
\hat\psi_t(\xi) = \frac 1{(2\pi)^3}\int_{\re Z = t}
\psi(Z) \cdot e^{-i \xi \cdot Z} \,dS(Z),
\qquad t \in \BB R,\: t>0,
$$
where  $\xi \cdot Z$ denotes $\re (\xi Z^+)$.
Since $\phi(X)$ is a polynomial, 
$\psi(Z) = \frac 1{(N(Z)+\rho^2)^2} \cdot
\phi \bigl( \frac {Z-\rho}{Z+\rho} \bigr)$
is non-singular for $\re Z>0$ and decays sufficiently fast for the integral to
be convergent. Applying the Fourier transform, we can write
$$
\psi(Z) = \int_{\BB R^3} \hat\psi_0(\xi) \cdot e^{i \xi \cdot Z} \,dS(\xi),
\qquad
\psi(W) = \int_{\BB R^3} \hat\psi_{\re W}(\xi) \cdot e^{i \xi \cdot \im W}
\,dS(\xi),
$$
where we identify $\BB R^3$ with $\{ Z \in \BB H;\: \re Z =0 \}$.
It follows from (\ref{Poisson-intermediate2}) that
\begin{equation}  \label{Poisson-intermediate3}
(N(W)+\rho^2) \int_{\BB R^3} \hat \psi_{\re W}(\xi) \cdot e^{i\xi\cdot \im W}
\,dS(\xi)
= \frac {(2l+1)\rho}{\pi^2} \int_{\BB R^3} \int_{\BB R^3}
\frac {\hat\psi_0(\xi) \cdot e^{i \xi \cdot Z}}{N(Z-W)} \,dS(\xi)dS(Z).
\end{equation}
The last integral is not absolutely convergent and should be understood as
the limit
$$
\lim_{R \to \infty}
\int_{Z \in \BB R^3, \: |Z|<R} \int_{\BB R^3}
\frac {\hat\psi_0(\xi) \cdot e^{i \xi \cdot Z}}{N(Z-W)} \,dS(\xi)dS(Z).
$$

\begin{lem}
$$
\lim_{R \to \infty} \int_{Z \in \BB R^3, \: |Z|<R}
\frac {e^{i \xi \cdot Z}}{N(Z-W)} \,dS(Z) =
2\pi^2 \cdot e^{i \xi \cdot \im W} \cdot \frac {e^{-|\xi| \re W}}{|\xi|},
$$
where $|\xi| = \bigl( (\xi^1)^2 + (\xi^2)^2 + (\xi^3)^2 \bigr)^{1/2}$.
\end{lem}

\pf
First we make a shift $Z=\im W+ Z'$, then the integral becomes
$$
e^{i \xi \cdot \im W} \lim_{R \to \infty} \int_{Z' \in \BB R^3, \: |Z|<R}
\frac {e^{i \xi \cdot Z'}}{(\re W)^2 + N(Z')} \,dS(Z').
$$
Clearly, the latter integral is invariant with respect to rotations of $\xi$.
Thus we may assume that $\xi = |\xi|e_1$, and the integral reduces to
$$
\lim_{R \to \infty} \int_{Z' \in \BB R^3, \: |Z|<R}
\frac {e^{i|\xi|z^1}}{(\re W)^2 + (z^1)^2 + (z^2)^2 + (z^3)^2} \,dS(Z).
$$
We set $\alpha = \bigl( (\re W)^2 + (z^2)^2 + (z^3)^2 \bigr)^{1/2}$ and
integrate with respect to $z^1$:
$$
\int_{-\infty}^{\infty} \frac {e^{i|\xi|z^1}}{(z^1)^2 + \alpha^2} \,dz^1
= \frac {\pi}{\alpha} e^{-|\xi|\alpha}.
$$
Thus we obtain
$$
\int_{\BB R^2} \frac {\pi}{\alpha} e^{-|\xi|\alpha} \,dz^2dz^3.
$$
The last integral can be computed in polar coordinates.
Writing $z^2 = r \cos \theta$, $z^3 = r \sin \theta$ and substituting
$s = ((\re W)^2 + r^2)^{1/2}$ we get:
$$
2 \pi^2 \int_0^{\infty} \frac {e^{-|\xi| \cdot ((\re W)^2 + r^2)^{1/2}}}
{((\re W)^2 + r^2)^{1/2}} r \,dr
= 2 \pi^2 \int_{\re W}^{\infty} e^{-|\xi|s} \,ds
= 2\pi^2 \cdot \frac {e^{-|\xi| \re W}}{|\xi|}.
$$
\qed

Let
$$
\Delta = \frac {\partial^2}{(\partial \xi^1)^2} +
\frac {\partial^2}{(\partial \xi^2)^2} +
\frac {\partial^2}{(\partial \xi^3)^2},
$$
then (\ref{Poisson-intermediate3}) can be rewritten as
\begin{multline*}
\int_{\BB R^3} \bigl(((\re W)^2+\rho^2 - \Delta)\hat\psi_{\re W} \bigr)
(\xi) \cdot e^{i \xi \cdot \im W} \,dS(\xi)  \\
= (N(W)+\rho^2) \int_{\BB R^3}
\hat\psi_{\re W}(\xi) \cdot e^{i \xi \cdot \im W} \,dS(\xi) \\
= 2(2l+1)\rho \int_{\BB R^3} \hat\psi_0(\xi) \cdot
e^{i \xi \cdot \im W} \cdot \frac {e^{-|\xi| \re W}}{|\xi|} \,dS(\xi).
\end{multline*}
Therefore,
$$
\biggl( -\frac 12 \Delta + \frac {(\re W)^2+\rho^2}2 \biggr)
\hat\psi_{\re W} (\xi)
= \frac {(2l+1)\rho}{|\xi|} \cdot e^{-|\xi| \re W} \cdot \hat\psi_0(\xi).
$$
Letting $\re W \to 0^+$, we obtain
\begin{equation}  \label{Schrod}
\biggl( -\frac 12 \Delta + \frac {\rho^2}2 \biggr) \hat\psi_0 (\xi)
= \frac {(2l+1)\rho}{|\xi|} \cdot \hat\psi_0(\xi).
\end{equation}
Let $\kappa = (2l+1)\rho$, then the equation becomes
\begin{equation}  \label{Schrodinger}
- \biggl( \frac 12 \Delta + \frac {\kappa}{|\xi|} \biggr) \hat\psi_0 (\xi)
= E \cdot \hat\psi_0(\xi),
\end{equation}
where
\begin{equation}  \label{E}
E = - \frac {\rho^2}2 = -\frac {\kappa^2}{2(2l+1)^2}.
\end{equation}
Thus we obtain the spectrum of (\ref{Schrodinger}) with eigenvalues given by
(\ref{E}). Note that the eigenfunctions can be found explicitly by applying
the Cayley and Fourier transforms to harmonic functions of degree $2l$.

\begin{rem}
One can also show that (\ref{E}) gives {\em all} the negative eigenvalues
and $\hat \phi_0(\xi)$ are all the corresponding eigenfunctions.
It is possible to reverse all the steps in the argument by presenting
$$
\frac {2\pi^2}{|\xi|} = \lim_{w^0 \to 0^+} \lim_{R \to \infty}
\int_{Z \in \BB R^3, \: |Z|<R} \frac {e^{i \xi \cdot Z}}{(w^0)^2+N(Z)} \,dS(Z),
$$
where $|\xi| = \bigl( (\xi^1)^2 + (\xi^2)^2 + (\xi^3)^2 \bigr)^{1/2}$.
Then one arrives at a version of (\ref{Poisson-hydro}) with $N(Y)<1$ and
$Y$ approaching the unit sphere.
Then the harmonic extension yields (\ref{Poisson-hydro}).
\end{rem}

Finally, the spectral decomposition of (\ref{Schrodinger}) with positive
eigenvalues $E$ comes from a Minkowski analogue of the Poisson formula and
will be considered in the next chapter.

\section{Quaternionic Analysis in $\HC$ and the Minkowski Space}
\label{M-section}

\subsection{The Space of Quaternions $\HC$ and the Minkowski Space $\BB M$}

Let $\HC = \BB C \otimes \BB H$ be the algebra of complexified quaternions.
We define a complex conjugation on $\HC$ with respect to $\BB H$:
$$
Z = z^0e_0 + z^1e_1 + z^2e_2 + z^3e_3 \quad \mapsto \quad
Z^c = \B{z^0}e_0 + \B{z^1}e_1 + \B{z^2}e_2 + \B{z^3}e_3,
\qquad z^0,z^1,z^2,z^3 \in \BB C,
$$
so that $\BB H = \{ Z \in \HC ;\: Z^c=Z \}$.
Then we can realize the Minkowski space $\BB M$ as a real form of $\HC$:
$$
\BB M = \{Z \in \HC ; \:Z^{c+} = - Z \}.
$$
Note that
$$
\HC \ni Z \mapsto N(Z) = \det Z
$$
is a quadratic form over $\BB C$.
The signature of this quadratic form restricted to $\BB M$ is $(3,1)$.
The corresponding symmetric bilinear form on $\HC$ is
\begin{equation}  \label{bilinear_form}
\langle Z, W \rangle = \frac 12 \tr(Z^+ W) = \frac 12 \tr(Z W^+).
\end{equation}

The algebra $\HC$ can be realized as complex $2 \times 2$ matrices:
$$
\HC = \biggl\{
Z= \begin{pmatrix} z_{11} & z_{12} \\ z_{21} & z_{22} \end{pmatrix}
; \: z_{ij} \in \BB C \biggr\}.
$$
Then the operations of quaternionic conjugation and complex conjugation
on $\HC$ become
$$
Z^+= \begin{pmatrix} z_{22} & -z_{12} \\ -z_{21} & z_{11} \end{pmatrix},
\qquad
Z^c = \begin{pmatrix}
\B{z_{22}} & -\B{z_{21}} \\ -\B{z_{12}} & \B{z_{11}} \end{pmatrix}.
$$
Observe that $Z^{c+}$ is the matrix adjoint $Z^*$.
Then
$$
\BB M = \{Z \in \HC ; \: Z^* = - Z \}
= \biggl\{
Z= \begin{pmatrix} z_{11} & z_{12} \\ z_{21} & z_{22} \end{pmatrix} \in \HC
; \: z_{11}, z_{22} \in i \BB R, \: z_{21} = -\overline{z_{12}} \biggr\}.
$$
The Minkowski space $\BB M$ is spanned over $\BB R$ by the four matrices
$$
\tilde e_0 = -ie_0 = \begin{pmatrix} -i & 0 \\ 0 & -i \end{pmatrix}, \qquad
e_1 = \begin{pmatrix} 0 & -i \\ -i & 0 \end{pmatrix}, \qquad
e_2 = \begin{pmatrix} 0 & -1 \\ 1 & 0 \end{pmatrix}, \qquad
e_3 = \begin{pmatrix} -i & 0 \\ 0 & i \end{pmatrix},
$$
so
$$
\BB M = \biggl\{ y^0 \tilde e_0 + y^1 e_1 + y^2 e_2 + y^3 e_3 =
\begin{pmatrix} -iy^0-iy^3 & -iy^1-y^2 \\ -iy^1+y^2 & -iy^0+iy^3 \end{pmatrix}
;\: y^0, y^1, y^2, y^3 \in \BB R \biggr\}.
$$
The quaternionic conjugation in this basis becomes
$$
\tilde e_0^+ = \tilde e_0,
\quad e_1^+ = -e_1, \quad e_2^+ = -e_2, \quad e_3^+ = -e_3.
$$
The elements $\tilde e_0$, $e_1$, $e_2$, $e_3$
are orthogonal with respect to the bilinear form (\ref{bilinear_form}) and
$\langle \tilde e_0, \tilde e_0 \rangle = -1$,
$\langle e_1, e_1 \rangle = \langle e_2, e_2 \rangle
= \langle e_3, e_3 \rangle = 1$.
We select the orientation of $\BB M$ so that $\{\tilde e_0, e_1, e_2, e_3 \}$
is a positively oriented basis. Set
$$
\|Z\| = \frac 1{\sqrt 2}
\sqrt{|z_{11}|^2 + |z_{12}|^2 + |z_{21}|^2 + |z_{22}|^2}, 
$$
so that $\|e_i\|=1$, $0 \le i \le 3$.
The corresponding Euclidean volume form on $\BB M$ is
$dV = dy^0 \wedge dy^1 \wedge dy^2 \wedge dy^3$.

We extend $Dx$ (defined on $\BB H$) to a holomorphic $\HC$-valued 3-form
on $\HC$ as
\begin{equation}  \label{Dz-explicit}
Dz = e_0 dz^1 \wedge dz^2 \wedge dz^3 - e_1 dz^0 \wedge dz^2 \wedge dz^3
+ e_2 dz^0 \wedge dz^1 \wedge dz^3 - e_3 dz^0 \wedge dz^1 \wedge dz^2,
\end{equation}
where we write $z_j = x_j + i \tilde x_j$, $x_j, \tilde x_j \in \BB R$,
and $dz^j = dx^j + i d\tilde x^j$, $0 \le j \le 3$.
Then we define a 3-form on $\BB M$ by $Dy = Dz \bigl |_{\BB M}$.

\begin{prop}
The 3-form $Dy$ takes values in $i \BB M \subset \HC$, is given explicitly by
\begin{equation}  \label{Dy-explicit}
Dy = i \tilde e_0 dy^1 \wedge dy^2 \wedge dy^3
+i e_1 dy^0 \wedge dy^2 \wedge dy^3
-i e_2 dy^0 \wedge dy^1 \wedge dy^3
+i e_3 dy^0 \wedge dy^1 \wedge dy^2,
\end{equation}
and satisfies the following property:
\begin{multline*}
\langle Y_1, Dy(Y_2,Y_3,Y_4) \rangle =
\frac 12 \tr (Y_1^+ \cdot Dy (Y_2,Y_3,Y_4)) =
-i \cdot dV(Y_1, Y_2, Y_3, Y_4)  \\
\qquad \forall Y_1, Y_2, Y_3, Y_4 \in \BB M.
\end{multline*}
\end{prop}

\subsection{Regular functions on $\BB M$ and $\HC$}

We introduce linear differential operators on $\BB M$
\begin{align*}
\nabla^+_{\BB M} &= - \tilde e_0 \frac{\partial}{\partial y^0} +
e_1 \frac{\partial}{\partial y^1} + e_2 \frac{\partial}{\partial y^2}
+ e_3 \frac{\partial}{\partial y^3}
\qquad \text{and}  \\
\nabla_{\BB M} &= - \tilde e_0 \frac{\partial}{\partial y^0} -
e_1 \frac{\partial}{\partial y^1} - e_2 \frac{\partial}{\partial y^2} -
e_3 \frac{\partial}{\partial y^3}
\end{align*}
which may be applied to functions on the left and on the right.

Fix an open subset $U \subset \BB M$, and let $f$ be a
differentiable function on $U$ with values in $\HC$.

\begin{df}
The function $f$ is {\em left-regular} if it satisfies
$$
(\nabla^+_{\BB M} f)(Y) = - \tilde e_0 \frac{\partial f}{\partial y^0}(Y) +
e_1 \frac{\partial f}{\partial y^1}(Y) +
e_2 \frac{\partial f}{\partial y^2}(Y) +
e_3 \frac{\partial f}{\partial y^3}(Y) =0,
\qquad \forall Y \in U.
$$
Similarly, $f$ is {\em right-regular} if
$$
(f\nabla^+_{\BB M} )(Y) = - \frac{\partial f}{\partial y^0}(Y) \tilde e_0 +
\frac{\partial f}{\partial y^1}(Y)e_1 +
\frac{\partial f}{\partial y^2}(Y)e_2 +
\frac{\partial f}{\partial y^3}(Y)e_3 =0,
\qquad \forall Y \in U.
$$
\end{df}

We define a second-order differential operator on $\BB M$
$$
\square_{3,1} = - \frac {\partial^2}{(\partial y^0)^2} +
\frac {\partial^2}{(\partial y^1)^2} +
\frac {\partial^2}{(\partial y^2)^2} +
\frac {\partial^2}{(\partial y^3)^2}
=\nabla_{\BB M} \nabla^+_{\BB M} = \nabla^+_{\BB M} \nabla_{\BB M}.
$$ 
Thus we have factored the wave operator on $\BB M$ into
two first order differential operators.

\begin{prop}
For any ${\cal C}^1$-function $f: U \to \HC$,
$$
d(f \cdot Dy) = df \wedge Dy = -i(f \nabla^+_{\BB M}) dV,
\qquad
d(Dy \cdot f) = - Dy \wedge df = -i(\nabla^+_{\BB M} f) dV.
$$
\end{prop}

\begin{cor}
Let $f: U \to \HC$ be of class ${\cal C}^1$, then
$$
\nabla^+_{\BB M} f = 0 \quad \Longleftrightarrow \quad Dy \wedge df =0,
\qquad
f \nabla^+_{\BB M} = 0 \quad \Longleftrightarrow \quad df \wedge Dy =0.
$$
\end{cor}

Thus the definitions of regular functions are analogous to the ones in
$\BB H$. However, as we will see in the next subsection, the important
difference is in the type of singularity of the kernels $k_0(X-Y)$, $k(X-Y)$.
This leads us to a more general view of quaternionic analysis that involves
functions on $\HC$, real forms of $\HC$ -- such as $\BB H$, $\BB M$ --
and relations between such functions given by complex continuation.

Let $U^{\BB C} \subset \HC$ be an open set and $f^{\BB C}: U^{\BB C} \to \HC$ a
differentiable function. We say that $f^{\BB C}$ is {\em holomorphic}
if it is holomorphic with respect to the complex variables
$z_{ij}$, $1 \le i,j \le 2$.
For such holomorphic function $f^{\BB C}$ the following derivatives are equal:
$$
\nabla^+_{\BB M} f^{\BB C} = \nabla^+ f^{\BB C} =
2 \begin{pmatrix} \frac {\partial}{\partial z_{22}} &
- \frac {\partial}{\partial z_{21}}  \\
- \frac {\partial}{\partial z_{12}} &
\frac {\partial}{\partial z_{11}} \end{pmatrix} f^{\BB C}, \qquad
f^{\BB C} \nabla^+_{\BB M} = f^{\BB C} \nabla^+ =
2 f^{\BB C} \begin{pmatrix} \frac {\partial}{\partial z_{22}} &
- \frac {\partial}{\partial z_{21}}  \\
- \frac {\partial}{\partial z_{12}} &
\frac {\partial}{\partial z_{11}} \end{pmatrix},
$$
$$
\nabla_{\BB M} f^{\BB C} = \nabla f^{\BB C} =
2 \begin{pmatrix} \frac {\partial}{\partial z_{11}} &
\frac {\partial}{\partial z_{21}}  \\
\frac {\partial}{\partial z_{12}} &
\frac {\partial}{\partial z_{22}} \end{pmatrix} f^{\BB C}, \qquad
f^{\BB C} \nabla_{\BB M} = f^{\BB C} \nabla =
2 f^{\BB C} \begin{pmatrix} \frac {\partial}{\partial z_{11}} &
\frac {\partial}{\partial z_{21}}  \\
\frac {\partial}{\partial z_{12}} &
\frac {\partial}{\partial z_{22}} \end{pmatrix}.
$$
Since we are interested in holomorphic functions only, we will abuse the
notation and denote by $\nabla$ and $\nabla^+$ the holomorphic differential
operators
$2 \begin{pmatrix} \frac {\partial}{\partial z_{11}} &
\frac {\partial}{\partial z_{21}}  \\
\frac {\partial}{\partial z_{12}} &
\frac {\partial}{\partial z_{22}} \end{pmatrix}$
and
$2 \begin{pmatrix} \frac {\partial}{\partial z_{22}} &
- \frac {\partial}{\partial z_{21}}  \\
- \frac {\partial}{\partial z_{12}} &
\frac {\partial}{\partial z_{11}} \end{pmatrix}$ respectively.

\begin{df}
Let $f^{\BB C}: U^{\BB C} \to \HC$ be a function.
We say that $f^{\BB C}$ is {\em holomorphic left-regular} if it is holomorphic
and $\nabla^+ f^{\BB C} =0$.
Similarly, $f^{\BB C}$ is {\em holomorphic right-regular} if it is holomorphic
and $f^{\BB C} \nabla^+ =0$.
\end{df}

As in the real case, we have:
\begin{lem}
For any holomorphic function $f^{\BB C}: U^{\BB C} \to \HC$,
$$
\nabla^+_{\BB M} f^{\BB C} = 0 \quad \Longleftrightarrow \quad
Dz \wedge df^{\BB C} =0 \quad \Longleftrightarrow \quad
\nabla^+ f^{\BB C} = 0,
$$
$$
f^{\BB C} \nabla^+_{\BB M} = 0 \quad \Longleftrightarrow \quad
df^{\BB C} \wedge Dz =0 \quad \Longleftrightarrow \quad
f^{\BB C} \nabla^+= 0.
$$
\end{lem}

The restriction relations
$$
Dz \bigl |_{\BB H} = Dx,
\qquad
Dz \bigl |_{\BB M} = Dy
$$
imply that the restriction of a holomorphic left- or right-regular function
to $U_{\BB H} = U^{\BB C} \cap \BB H$ produces a left- or right-regular
function on $U_{\BB H}$ respectively.
And the restriction of a holomorphic left- or right-regular function
to $U_{\BB M} = U^{\BB C} \cap \BB M$ also produces a left- or right-regular
function on $U_{\BB M}$ respectively.
Conversely, if one starts with, say, a left-regular function on $\BB M$,
extends it holomorphically to a left-regular function on $\HC$
and then restricts the extension to $\BB H$, the resulting function
is left-regular on $\BB H$.

\begin{lem}  \label{closed-M}
We have:
\begin{enumerate}
\item
$\square_{3,1} \frac 1{N(Z)} = 0$;
\item
$\nabla \frac 1{N(Z)}
= \frac 1{N(Z)} \nabla
= -2 \frac {Z^{-1}}{N(Z)} = -2 \frac {Z^+} {N(Z)^2}$;
\item
$\frac {Z^{-1}}{N(Z)} = \frac {Z^+} {N(Z)^2}$
is a holomorphic left- and right-regular function defined wherever
$N(Z) \ne 0$;
\item
The form
$\frac {Z^{-1}}{N(Z)} \cdot Dz = \frac {Z^+} {N(Z)^2} \cdot Dz$ is a closed
holomorphic $\HC$-valued 3-form defined wherever $N(Z) \ne 0$.
\end{enumerate}
\end{lem}

\subsection{Fueter Formula for Holomorphic Regular Functions on $\BB M$}

We are interested in extensions of the Cauchy-Fueter formula
(Theorem \ref{Fueter}) to functions on $\BB M$.
Let $U \subset \BB M$ be an open subset, and let $f$ be a ${\cal C}^1$-function
defined on a neighborhood of $\overline{U}$ such that $\nabla^+_{\BB M} f =0$.
In this section we extend the Cauchy-Fueter integral formula
to left-regular functions which can be extended holomorphically to a
neighborhood of $\overline{U}$ in $\HC$.
Observe that the expression in the integral formula
$\frac {(X - Y)^{-1}}{N(X - Y)} \cdot Dx$
is nothing else but the restriction to $\BB H$ of the
holomorphic 3-form $\frac {(Z-Y)^{-1}} {N(Z-Y)} \cdot Dz$
which is the 3-form from Lemma \ref{closed-M} translated by $Y$.
For this reason we expect an integral formula of the kind
$$
f(Y_0) = \frac 1{2\pi^2} \int_{\partial U}
\biggl( \frac {(Z-Y_0)^{-1}} {N(Z-Y_0)} \cdot Dz \biggr)
\biggl|_{\BB M} \cdot f(Y),
\qquad \forall Y_0 \in U.
$$
However, the integrand is singular wherever $N(Z-Y_0)=0$.
We resolve this difficulty by deforming the contour of integration $\partial U$
in such a way that the integral is no longer singular.

Fix an $\epsilon \in \BB R$ 
and define the $\epsilon$-deformation
$h_{\epsilon}: \HC \to \HC$, $Z \mapsto Z_{\epsilon}$, by:
\begin{align*}
z_{11} \quad &\mapsto \quad z_{11} - i\epsilon z_{22}  \\
z_{22} \quad &\mapsto \quad z_{22} - i\epsilon z_{11}  \\
z_{12} \quad &\mapsto \quad z_{12} + i\epsilon z_{12}  \\
z_{21} \quad &\mapsto \quad z_{21} + i\epsilon z_{21}.
\end{align*}

Define a quadratic form on $\HC$
$$
S(Z)= - (z_{11}^2 + z_{22}^2 + 2 z_{12}z_{21}).
$$

\begin{lem}
We have the following identities:
$$
Z_{\epsilon} = Z - i \epsilon Z^+, \qquad
(Z_{\epsilon})^+ = Z^+ - i \epsilon Z,
$$
$$
N(Z_{\epsilon}) = (1-\epsilon^2) N(Z) + i \epsilon S(Z),
$$
$$
S(Y) = 2\|Y\|^2, \qquad \forall Y \in \BB M.
$$
\end{lem}

For $Z_0 \in \HC$ fixed, we use the notation
$$
h_{\epsilon, Z_0} (Z) = Z_0 + h_{\epsilon}(Z-Z_0) = Z - i\epsilon (Z-Z_0)^+.
$$

\begin{thm}  \label{holomorphic_Fueter}
Let $U \subset \BB M$ be an open bounded subset with piecewise
${\cal C}^1$ boundary $\partial U$,
and let $f(Y)$ be a ${\cal C}^1$-function defined on a
neighborhood of the closure $\overline{U}$ such that $\nabla_{\BB M}^+ f =0$.
Suppose that $f$ extends to a holomorphic left-regular function
$f^{\BB C} : V^{\BB C} \to \HC$ with $V^{\BB C} \subset \HC$ an open 
subset containing $\overline{U}$, then
$$
-\frac 1{2\pi^2} \int_{(h_{\epsilon, Y_0})_*(\partial U)}
\frac {(Z-Y_0)^{-1}} {N(Z-Y_0)} \cdot Dz \cdot f^{\BB C}(Z) =
\begin{cases}
f(Y_0) & \text{if $Y_0 \in U$;} \\
0 & \text{if $Y_0 \notin \overline{U}$.}
\end{cases}
$$
for all $\epsilon \ne 0$ sufficiently close to 0.
\end{thm}

\begin{rem}
For all $\epsilon \ne 0$ sufficiently close to 0 the contour of integration
$(h_{\epsilon, Y_0})_*(\partial U)$ lies inside $V^{\BB C}$ and the integrand
is non-singular, thus the integrals are well defined. Moreover, we will see
that the value of the integral becomes constant when the parameter $\epsilon$
is sufficiently close to 0.
\end{rem}

\pf
Let $M = \sup_{Y \in \partial U} \|Y-Y_0\|$.
Without loss of generality we may assume that
$V^{\BB C}$ is the $\delta$-neighborhood of $\overline{U}$
for some $\delta >0$.
We will show that the integral formula holds for
$0 < |\epsilon| < \delta/2M$.
Clearly, for this choice of $\epsilon$ the contour of integration
$(h_{\epsilon,Y_0})_*(\partial U)$ lies inside $V^{\BB C}$ and, since
the integrand is a closed form, the integral stays constant
for $-\delta/2M < \epsilon < 0$ and $0 < \epsilon < \delta/2M$
(a priori the values of the integral may be different on these two intervals).

Since the case $Y_0 \notin \overline{U}$ is trivial, we assume $Y_0 \in U$.
Let $S_r = \{ Y \in \BB M ;\: \|Y-Y_0\| = r \}$
and $B_r = \{ Y \in \BB M ;\: \|Y-Y_0\| < r \}$
be the sphere and the ball of radius $r$ centered at $Y_0$,
and choose $r>0$ sufficiently small so that $B_r \subset U$ and
$r < \delta/2$. By Stokes'
\begin{multline*}
\int_{(h_{\epsilon, Y_0})_*(\partial U)}
\frac {(Z-Y_0)^{-1}} {N(Z-Y_0)} \cdot Dz \cdot f^{\BB C}(Z)
= \int_{(h_{\epsilon, Y_0})_*(S_r)}
\frac {(Z-Y_0)^{-1}} {N(Z-Y_0)} \cdot Dz \cdot f^{\BB C}(Z)  \\
= \int_{(h_{1, Y_0})_*(S_r)}
\frac {(Z-Y_0)^{-1}} {N(Z-Y_0)} \cdot Dz \cdot f^{\BB C}(Z),
\end{multline*}
where $S_r$ is oriented as the boundary of $B_r$.

Let $\operatorname{P}_{\BB H}$ be the projection
$\HC \twoheadrightarrow \BB H$ defined
by 
\begin{multline*}
Z = (x^0 + i \tilde x^0) e_0
+ (x^1 + i \tilde x^1) e_1
+ (x^2 + i \tilde x^2) e_2
+ (x^3 + i \tilde x^3) e_3  \\
\mapsto
X = x^0 e_0 + x^1 e_1 + x^2 e_2 + x^3 e_3,
\qquad
x^0,  x^1,  x^2,  x^3,
\tilde x^0,  \tilde x^1,  \tilde x^2,  \tilde x^3 \in \BB R,
\end{multline*}
and let $\operatorname{P}_{\BB H+Y_0} : \HC \twoheadrightarrow \BB H + Y_0$,
be the projection
$\operatorname{P}_{\BB H+Y_0} (Z) = \operatorname{P}_{\BB H} (Z-Y_0) + Y_0$.
We describe the supports of the cycles involved in integration
\begin{align*}
|S_r| &=
\{ Y_0 + a \tilde e_0 + be_1 + ce_2 + de_3 ;\: a^2 + b^2 + c^2 + d^2 = r^2 \},
\\
|(h_{1, Y_0})_*(S_r)| &= \{ Y_0 + (1-i)a \tilde e_0
+ (1+i)be_1 + (1+i)ce_2 + (1+i)de_3 ;  \\
& \hskip2.7in a^2 + b^2 + c^2 + d^2 = r^2 \},  \\
|(\operatorname{P}_{\BB H + Y_0} \circ h_{1,Y_0})_*(S_r)| &=
\{ Y_0 - ae_0 + be_1 + ce_2 + de_3 ;\: a^2 + b^2 + c^2 + d^2 = r^2 \}  \\
&= \{ X \in \BB H+Y_0 ;\: \|X-Y_0\| = r \}.
\end{align*}
Let $\tilde S_r = \{ X \in \BB H+Y_0 ;\: \|X-Y_0\| = r \}$
be the sphere oriented as the boundary of the open ball, then
$$
(\operatorname{P}_{\BB H + Y_0} \circ h_{1,Y_0})_*(S_r) = - \tilde S_r
$$
as 3-cycles.
The cycles $(h_{1, Y_0})_*(S_r)$ and $- \tilde S_r$ are homologous to each
other inside $\HC \setminus \{ N (Z - Y_0) =0 \}$.
Then by Stokes' again,
$$
\int_{(h_{1,Y_0})_*(S_r)}
\frac {(Z-Y_0)^{-1}} {N(Z-Y_0)} \cdot Dz \cdot f^{\BB C}(Z)  \\
= - \int_{\tilde S_r}
\frac {(Z-Y_0)^{-1}} {N(Z-Y_0)} \cdot Dz \cdot f^{\BB C}(Z).
$$
Finally, by the Fueter formula for the regular quaternions
(Theorem \ref{Fueter}), the last integral is $-2\pi^2 f(Y_0)$.
\qed

One can drop the assumption that $f(Y)$ extends to an open neighborhood
of $\overline{U}$ in $\HC$ and prove the following version of the Fueter
formula on $\BB M$ involving generalized functions:

\begin{thm}
Let $U \subset \BB M$ be a bounded open region with smooth boundary
$\partial U$.
Let $f: U \to \HC$ be a function which extends to a real-differentiable
function on an open neighborhood $V \subset \BB M$ of the closure
$\overline{U}$ such that $\nabla^+_{\BB M} f = 0$.
Then, for any point $Y_0 \in \BB M$ such that $\partial U$ intersects the
cone $\{ Y \in \BB M ;\: N (Y-Y_0) =0 \}$ transversally, we have:
$$
\lim_{\epsilon \to 0} \frac {-1}{2\pi^2} \int_{\partial U}
\frac {(Y-Y_0)^+} {\bigl( N(Y-Y_0) + i\epsilon \|Y-Y_0\|^2 \bigr)^2}
\cdot Dy \cdot f(Y)
=
\begin{cases}
f(Y_0) & \text{if $Y_0 \in U$;} \\
0 & \text{if $Y_0 \notin \overline{U}$.}
\end{cases}
$$
\end{thm}

\subsection{Fueter Formula for Hyperboloids in $\BB M$}

Now we would like to extend the Fueter formula in $\BB M$ to certain
non-compact cycles. Then we need to require additional regularity
conditions at infinity, which can be defined using conformal transformations
in $\HC$. We illustrate the general idea with an important example of
two-sheeted hyperboloids in $\BB M$ that has a representation theoretic
significance and an application to physics.
We realize the group $U(2)$ as
$$
U(2)=\{Z \in \HC ;\: Z^* = Z^{-1} \}.
$$

\begin{lem}
Consider an element
$\gamma = \begin{pmatrix} i & 1 \\ i & -1 \end{pmatrix} \in GL(2,\HC)$
with $\gamma^{-1} = \begin{pmatrix} -i & -i \\ 1 & -1 \end{pmatrix}$.
The fractional linear map on $\HC$
$$
\pi_l(\gamma): \: Z \mapsto -i(Z+1)(Z-1)^{-1}
$$
maps $\BB M \to U(2)$, has no singularities on $\BB M$, and sends
the two-sheeted hyperboloid $\{ Y \in \BB M;\: N(Y) = -1 \}$ into
the sphere $\{ Z \in U(2); N(Z)=1 \} = SU(2)$.

Conversely, the fractional linear map on $\HC$
$$
\pi_l(\gamma^{-1}): \: Z \mapsto (Z-i)(Z+i)^{-1}
$$
maps $U(2) \to \BB M$ (with singularities), and sends
the sphere $\{ Z \in U(2); N(Z)=1 \} = SU(2)$ into 
the two-sheeted hyperboloid $\{ Y \in \BB M;\: N(Y) = -1 \}$.
The singularities of $\pi_l(\gamma^{-1})$ on $SU(2)$ lie along the sphere
$\{ X \in SU(2) ; \re X =0 \}$.
\end{lem}

This way $U(2)$ can be regarded as a compactification of $\BB M$.

\begin{df}
A left- (respectively right-) regular function $f(Y)$ (respectively $g(Y)$)
on $\BB M$ is {\em regular at infinity} if $\pi_l(\gamma^{-1})f$
(respectively $\pi_r(\gamma^{-1})g$) extends to a regular function on $U(2)$.

We denote by ${\cal V}(\BB M)$ (respectively ${\cal V}'(\BB M)$) the spaces
of left- (respectively right-) regular function $f(Y)$ (respectively $g(Y)$)
on $\BB M$ such that $\pi_l(\gamma^{-1})f$ (respectively $\pi_r(\gamma^{-1})g$)
extends holomorphically to an open neighborhood of $U(2)$ in $\HC$.
\end{df}

Next we observe

\begin{lem}
Let $Y,Y' \in \BB M$ be such that $N(Y) \ne N(Y')$.
Then, for all sufficiently small $\epsilon \in \BB R \setminus \{0\}$,
$$
N \bigl( Y - (1+i\epsilon)Y' \bigr) \ne 0.
$$
\end{lem}

For $Y \in \BB M$ with $N(Y)<0$ and $\epsilon \in \BB R$, we define
$$
Y^{\epsilon} =
\begin{cases}
(1+i\epsilon)Y & \text{if $i\tr Y>0$;} \\
(1-i\epsilon)Y & \text{if $i\tr Y<0$.}
\end{cases}
$$
(Note that elements in $\BB M$ always have purely imaginary traces.)

\begin{lem}
For $Y \in \BB M$ we have:
\begin{align*}
\pi_l(\gamma)(Y) & \text{ has both eigenvalues of length $1$;}  \\
\text{if $\epsilon>0$, } \pi_l(\gamma) (Y^{\epsilon})
&\text{ has both eigenvalues of length $>1$;}  \\
\text{if $\epsilon<0$, } \pi_l(\gamma) (Y^{\epsilon})
&\text{ has both eigenvalues of length $<1$.}
\end{align*}
\end{lem}

Let $R>0$, and consider a two-sheeted hyperboloid
$H_R = \{ Y \in \BB M ;\: N(Y)=-R^2 \}$.
We orient it so that $\{e_1,e_2,e_3\}$ form positively oriented bases
of the tangent spaces of $H_R$ at $\pm iR$.
This way $\pi_l(\gamma): H_1 \to SU(2)$
preserves the orientations.
For an $\epsilon \in \BB R$, we define a deformed hyperboloid in $\HC$:
$$
H_R^{\epsilon} = \{ Y^{\epsilon} ;\: Y \in \BB M,\: N(Y)=-R^2 \}.
$$
The orientation of $H_R^{\epsilon}$ is induced from that of $H_R$.

\begin{thm}  \label{Fueter-hyperboloid}
Let $f \in {\cal V}(\BB M)$.
Then, for $\epsilon > 0$ sufficiently small,
$$
f(Y_0) = \frac 1{2\pi^2} \int_{H_R^{\epsilon} - H_{R'}^{-\epsilon}}
\frac {(Z-Y_0)^{-1}} {N(Z-Y_0)} \cdot Dz \cdot f(Z),
\qquad \forall Y_0 \in \BB M,\: N(Y_0)<0,
$$
for any $R,R'>0$. (In particular, the integral converges absolutely.)
\end{thm}

\pf
The map $\pi_l(\gamma)$ sends $H_R^{\epsilon}$ and $H_{R'}^{-\epsilon}$ into
some cycles in $\HC$ which we can call $C_R^{\epsilon}$ and
$C_{R'}^{-\epsilon}$ with compact supports. (The closure of the images
$\pi_l(\gamma) (H_R^{\epsilon})$ and $\pi_l(\gamma) (H_{R'}^{-\epsilon})$
will contain the sphere $\{ X \in SU(2) ; \re X =0 \}$.)
The orientations were chosen so that the chains $C_R^{\epsilon}$ and
$C_{R'}^{-\epsilon}$ have no boundary and are, in fact, cycles.
Changing the variables
$Z = (Z'-i)(Z'+i)^{-1}$, $Y_0 = (Y_0'-i)(Y_0'+i)^{-1}$
and using Lemma \ref{X-Y_difference} with Proposition \ref{Dx-pullback},
we can rewrite
\begin{multline*}
\int_{H_R^{\epsilon} - H_{R'}^{-\epsilon}}
\frac {(Z-Y_0)^{-1}} {N(Z-Y_0)} \cdot Dz \cdot f(Z)  \\
= N(Y_0'+i) \cdot (Y_0'+i) \int_{C_R^{\epsilon} - C_{R'}^{-\epsilon}}
\frac {(Z' - Y_0')^{-1}} {N(Z' - Y_0')} \cdot Dz \cdot
\frac{(Z'+i)^{-1}}{N(Z'+i)} \cdot f \bigl( (Z'-i)(Z'+i)^{-1} \bigr)  \\
= 2\pi^2 \cdot f \bigl( (Y_0'-i)(Y_0'+i)^{-1} \bigr) = 2\pi^2 \cdot f(Y_0)
\end{multline*}
since the Fueter formula applies here.
\qed

This version of the Fueter formula suggests a natural polarization of
${\cal V}(\BB M)$ given by integrals over a single hyperboloid.
This polarization has a natural representation theoretic interpretation
in terms of conformal group $SU(2,2) \subset U(2,2)$
which will be discussed in the next subsection.

\subsection{Cayley Transform and Polarization}

The Lie algebra of $U(2,2)$ has already appeared in
Theorem \ref{unitary-u(2,2)-action} as the real form of
$\mathfrak{gl}(4,\BB C)$ which acts unitarily on ${\cal H}^{\pm}$,
${\cal V}^{\pm}$ and ${\cal V}'^{\pm}$.
The group $U(2,2)$ can be realized as the subgroup of elements of $GL(2,\HC)$
preserving the Hermitian form on $\BB C^4$ given by the $4 \times 4$ matrix
$\begin{pmatrix} 1 & 0 \\ 0 & -1 \end{pmatrix}$. Explicitly,
\begin{align*}
U(2,2) &= \Biggl\{ \begin{pmatrix} a & b \\ c & d \end{pmatrix};\:
a,b,c,d \in \HC,\:
\begin{matrix} a^*a = 1+c^*c \\ d^*d = 1+b^*b \\ a^*b=c^*d \end{matrix}
\Biggr\}  \\
&= \Biggl\{ \begin{pmatrix} a & b \\ c & d \end{pmatrix};\:
a,b,c,d \in \HC,\:
\begin{matrix} a^*a = 1+b^*b \\ d^*d = 1+c^*c \\ ac^*=bd^* \end{matrix}
\Biggr\}.
\end{align*}
The Lie algebra of $U(2,2)$ is
\begin{equation}  \label{u(2,2)-algebra}
\mathfrak{u}(2,2) = \Bigl\{
\begin{pmatrix} A & B \\ B^* & D \end{pmatrix} ;\: A,B,D \in \HC ,\:
A=-A^*, D=-D^* \Bigr\}.
\end{equation}
If $\begin{pmatrix} a & b \\ c & d \end{pmatrix} \in U(2,2)$,
then $\begin{pmatrix} a & b \\ c & d \end{pmatrix}^{-1} =
\begin{pmatrix} a^* & -c^* \\ -b^* & d^* \end{pmatrix}$.

From Lemma \ref{X-Y_difference} we obtain:

\begin{cor}  \label{N(Z-W)}
For an $h = \begin{pmatrix} a' & b' \\ c' & d' \end{pmatrix} \in U(2,2)$
with $h^{-1} = \begin{pmatrix} a & b \\ c & d \end{pmatrix}$,
let $\tilde Z = (aZ+b)(cZ+d)^{-1}$ and $\tilde W = (aW+b)(cW+d)^{-1}$.
Then
$$
N(Z-W)^2 = N(cZ+d) \cdot N(a^*+Zb^*) \cdot N(\tilde Z - \tilde W)^2
\cdot N(cW+d) \cdot N(a^*+Wb^*).
$$
\end{cor}

The fractional linear actions $\pi_l$ and $\pi_r$ of $U(2,2)$
preserve the bounded domain
$$
\BB D^+ = \{ Z \in \HC;\: ZZ^*<1 \},
$$
where the inequality $ZZ^*<1$ means that the matrix $1-ZZ^*$ is positive
definite. The Shilov boundary of $\BB D^+$ is $U(2)$.
Similarly, we can define
$$
\BB D^- = \{ Z \in \HC;\: ZZ^*>1 \}.
$$
The fractional linear actions $\pi_l$ and $\pi_r$ of $U(2,2)$ preserve
$\BB D^-$ and $U(2)$.

A $U(2)$ bi-invariant measure on $U(2)$ is given by the restriction of the
holomorphic 4-form defined on $\HC$
$$
\frac {dZ^4}{N(Z)^2},
\qquad dZ^4 = dz_{11} \wedge dz_{12} \wedge dz_{21} \wedge dz_{22},
$$
to $U(2)$. The tangent space at $1 \in U(2)$ can be identified with $\BB M$.
The space $\BB M$ was oriented by $\{ \tilde e_0,e_1,e_2,e_3 \}$,
which in turn induces an orientation on $U(2)$.

\begin{lem}
With this orientation convention we have:
$$
\int_{U(2)} \frac {dZ^4}{N(Z)^2} = -8\pi^3i.
$$
\end{lem}

We will need the following Jacobian lemma:

\begin{lem}  \label{Jacobian_lemma}
On $\HC$ we have:
$$
dZ^4 = N(cZ+d)^2 \cdot N(a'-Zc')^2 \,d \tilde Z^4,
$$
where $h = \begin{pmatrix} a' & b' \\ c' & d' \end{pmatrix} \in GL(2,\HC)$,
$h^{-1} = \begin{pmatrix} a & b \\ c & d \end{pmatrix}$
and $\tilde Z = (aZ+b)(cZ+d)^{-1}$.
\end{lem}

\pf
The group $GL(2,\HC)$ is connected, so it is sufficient to verify
the identity on the Lie algebra level only. Let
$$
\begin{pmatrix} a & b \\ c & d \end{pmatrix} =
\exp \Bigl( t \begin{pmatrix} A & B \\ C & D \end{pmatrix} \Bigr),
$$
then, modulo terms of order $t^2$ and higher, we have:
$$
\tilde Z = (aZ +b)(cZ +d)^{-1}
= Z + t(AZ+B-ZCZ-ZD).
$$
Hence
$$
d \tilde Z = dZ + t(A \cdot dZ - dZ \cdot CZ - ZC \cdot dZ - dZ \cdot D).
$$
It follows that
$$
d \tilde Z^4 = \bigl( 1 + 2 t \tr (A - CZ - ZC - D) \bigr) \,dZ^4.
$$
On the other hand,
$$
N(cZ+d) = N(1+t(CZ+D)) = 1 +t \tr (CZ+D),
$$
and
$$
N(a'-Zc') = N(1+t(-A+ZC)) = 1 +t \tr (-A+ZC),
$$
and the result follows.
\qed

We also consider a subgroup $U(2,2)'$ of $GL(2,\HC)$ conjugate to $U(2,2)$
and preserving the Hermitian form on $\BB C^4$ given by
$\begin{pmatrix} 0 & 1 \\ 1 & 0 \end{pmatrix}$:
\begin{align*}
U(2,2)' &= 
\biggl\{ \begin{pmatrix} a & b \\ c & d \end{pmatrix} \in GL(2,\HC) ;\:
\begin{pmatrix} a & b \\ c & d \end{pmatrix}
\begin{pmatrix} 0 & 1 \\ 1 & 0 \end{pmatrix}
\begin{pmatrix} a^* & c^* \\ b^* & d^* \end{pmatrix}
= \begin{pmatrix} 0 & 1 \\ 1 & 0 \end{pmatrix} \biggr\}  \\
&= \Biggl\{ \begin{pmatrix} a & b \\ c & d \end{pmatrix};\:
a,b,c,d \in \HC,\:
\begin{matrix} ab^*+ba^*=0 \\ cd^*+dc^*=0 \\ ad^*+bc^*=1 \end{matrix}
\Biggr\}  \\
&= \Biggl\{ \begin{pmatrix} a & b \\ c & d \end{pmatrix};\:
a,b,c,d \in \HC,\:
\begin{matrix} a^*c+c^*a=0 \\ b^*d+d^*b=0 \\ a^*d+c^*b=1 \end{matrix}
\Biggr\}.
\end{align*}
Its Lie algebra is
$$
\mathfrak{u}(2,2)' = 
\biggl\{ \begin{pmatrix} A & B \\ C & -A^* \end{pmatrix} ;\:
A,B,C \in \HC,\: B^*=-B,\: C^*=-C \biggr\}.
$$
If $\begin{pmatrix} a & b \\ c & d \end{pmatrix} \in U(2,2)'$, then
$\begin{pmatrix} a & b \\ c & d \end{pmatrix}^{-1}
= \begin{pmatrix} d^* & b^* \\ c^* & a^* \end{pmatrix}$.

From Lemma \ref{X-Y_difference} we obtain:

\begin{cor} \label{N(Z-W)'}
For an $h = \begin{pmatrix} a' & b' \\ c' & d' \end{pmatrix} \in U(2,2)'$
with $h^{-1} = \begin{pmatrix} a & b \\ c & d \end{pmatrix}$,
let $\tilde Z = (aZ+b)(cZ+d)^{-1}$ and $\tilde W = (aW+b)(cW+d)^{-1}$.
Then
$$
N(Z-W)^2 = N(cZ+d) \cdot N(d^*-Zc^*) \cdot N(\tilde Z - \tilde W)^2
\cdot N(cW+d) \cdot N(d^*-Wc^*).
$$
\end{cor}
Let $\C^+$ denote the open cone
\begin{equation}  \label{C^+}
\C^+ = \{ Y \in \BB M ;\: N(Y)<0, \: i\tr Y <0 \},
\end{equation}
and define the tube domains in $\HC$
$$
\BB T^+ = \BB M + i\C^+,
\qquad \BB T^- = \BB M - i\C^+.
$$
Then $\BB M$ is the Shilov boundary of $\BB T^+$ and $\BB T^-$.
The group $U(2,2)'$ acts meromorphically on $\BB T^+$, $\BB T^-$
and the Minkowski space $\BB M = \{ Z \in \HC ;\: Z^*=-Z \}$.

\begin{lem}
Recall
$\gamma = \begin{pmatrix} i & 1 \\ i & -1 \end{pmatrix} \in GL(2,\HC)$
with $\gamma^{-1} = \begin{pmatrix} -i & -i \\ 1 & -1 \end{pmatrix}$.
The fractional linear map on $\HC$
$$
\pi_l(\gamma): \: Z \mapsto -i(Z+1)(Z-1)^{-1}
$$
maps $\BB T^+ \to \BB D^+$, $\BB T^- \to \BB D^-$, and its inverse 
$$
\pi_l(\gamma^{-1}): \: Z \mapsto (Z-i)(Z+i)^{-1}
$$
maps $\BB D^+ \to \BB T^+$, $\BB D^- \to \BB T^-$.

Also, $\pi_l(\gamma): \BB M \to U(2)$ and
$\pi_l(\gamma^{-1}): \: U(2) \to \BB M$ preserve the orientations.
\end{lem}

We denote by ${\cal H}(\BB M)^+$ the space of holomorphic functions $\phi$
on $\BB T^+$ such that $\square \phi =0$ and $\pi^0_l (\gamma^{-1})\phi$
extends to a holomorphic function defined on a neighborhood of the closure
of $\BB D^+$.
Similarly, we can define ${\cal H}(\BB M)^-$ as the space of
holomorphic functions $\phi$ on $\BB T^-$ such that $\square \phi =0$ and
$\pi^0_l (\gamma^{-1})\phi$ extends to a holomorphic function defined
in a neighborhood of the closure of $\BB D^-$ and regular at infinity.
By the same fashion we can define the spaces of left- and right-regular
functions ${\cal V}(\BB M)^{\pm}$ and ${\cal V}'(\BB M)^{\pm}$.
Thus we obtain polarizations:
$$
{\cal H}(\BB M) = {\cal H}(\BB M)^+ \oplus {\cal H}(\BB M)^-,
\qquad
{\cal V}(\BB M) = {\cal V}(\BB M)^+ \oplus {\cal V}(\BB M)^-,
\qquad
{\cal V}'(\BB M) = {\cal V}'(\BB M)^+ \oplus {\cal V}'(\BB M)^-.
$$
In terms of these spaces, Theorem \ref{Fueter-hyperboloid} can be restated
as follows:

\begin{thm}
For all $f \in {\cal V}(\BB M)^+$, all $Y_0 \in \BB T^+$ and all $R>0$ we have:
$$
f(Y_0) = \frac 1{2\pi^2} \int_{H_R}
\frac {(Z-Y_0)^{-1}} {N(Z-Y_0)} \cdot Dz \cdot f(Z).
$$
Similarly, for all $f \in {\cal V}(\BB M)^-$, all $Y_0 \in \BB T^-$
and all $R>0$ we have:
$$
f(Y_0) = -\frac 1{2\pi^2} \int_{H_R}
\frac {(Z-Y_0)^{-1}} {N(Z-Y_0)} \cdot Dz \cdot f(Z).
$$
\end{thm}

\subsection{Solutions of the Wave Equation}   \label{Poisson-M}

In this subsection we prove an analogue of Theorem \ref{Poisson}
for the solutions of $\square_{3,1} \phi = 0$ on $\BB M$.
Recall that $\widetilde \deg$ is the degree operator plus identity:
$$
\widetilde{\deg} = 1 + y^0 \frac {\partial}{\partial y^0} +
y^1 \frac {\partial}{\partial y^1} + y^2 \frac {\partial}{\partial y^2} +
y^3 \frac {\partial}{\partial y^3}.
$$
For an $R>0$, $H_R$ denotes the two-sheeted hyperboloid
$\{ Y \in \BB M ;\: N(Y)=-R^2 \}$ oriented so that $\{e_1,e_2,e_3\}$
form positively oriented bases of the tangent spaces of $H_R$ at $\pm iR$.
This way $\pi_l(\gamma): H_1 \to SU(2)$ preserves the orientations.


\begin{thm}  \label{Poisson-M-thm}
For all $\phi \in {\cal H} (\BB M)^+$, all $Z \in \BB T^+$ and all $R>0$
we have:
$$
\phi \bigl( Z \bigr) =
- \frac 1{2\pi^2} \int_{Y \in H'_R}
\frac {R^2 + N(Z)}{N(Y-Z)^2} \cdot \phi(Y) \,\frac {Dy}{Y}
= \frac 1{2\pi^2} \int_{Y \in H'_R}
\frac {\bigl( \widetilde{\deg} \phi \bigr)(Y)}{N(Y-Z)} \,\frac {Dy}{Y}.
$$
Similarly, for all $\phi \in {\cal H} (\BB M)^-$, all $Z \in \BB T^-$
and all $R>0$ we have:
$$
\phi \bigl( Z \bigr) =
- \frac 1{2\pi^2} \int_{Y \in H'_R}
\frac {R^2 + N(Z)}{N(Y-Z)^2} \cdot \phi(Y) \,\frac {Dy}{Y}
= \frac 1{2\pi^2} \int_{Y \in H'_R}
\frac {\bigl( \widetilde{\deg} \phi \bigr)(Y)}{N(Y-Z)} \,\frac {Dy}{Y}.
$$
The 3-form $Dy/Y$ equals $-i\frac {dS}{\|Y\|}$ on
$\{ Y \in H_R ;\: i\tr Y>0 \}$ and $i\frac {dS}{\|Y\|}$ on
$\{ Y \in H_R ;\: i\tr Y<0 \}$,
where $dS$ denontes the Euclidean volume form on $H_R$.
\end{thm}

\pf
First we check the formula for $R=1$:
The map $\pi_l(\gamma)$ sends the cycle $H_1$ into the sphere
$S^3 \subset \BB H$ of radius one.
Changing the variables
$Z = (Z'-i)(Z'+i)^{-1}$, $Y = (Y'-i)(Y'+i)^{-1}$,
and using Lemma \ref{X-Y_difference} with Proposition \ref{Dx-pullback},
we can rewrite
\begin{multline*}
\int_{Y \in H_1}
\frac {1 + N(Z)}{N(Y-Z)^2} \cdot \frac {Dy}Y \cdot \phi(Y) =  \\
\int_{Y' \in S^3}
\frac {1 + N \bigl( \frac{Z'-i}{Z'+i} \bigr)} {N(Y'-Z')^2} \cdot
N(Z'+i)^2 \cdot \frac{N(1-iY')^2}2 \cdot \frac{Y'+i}{Y'-i} \cdot
\frac {(1-iY')^{-1}}{N(1-iY')} \cdot Dy \cdot \frac {(Y'+i)^{-1}}{N(Y'+i)}
\cdot \phi(Y)  \\
=
\frac {N(Z'+i)}2 \int_{Y' \in S^3}
\frac {N(Z'+i) + N(Z'-i)} {N(Y'-Z')^2} \cdot
\frac{N(1-iY')}{(Y'-i) \cdot (1-iY')} \cdot \frac {Dy}{N(Y'+i)}
\cdot \phi(Y)  \\
=
N(Z'+i) \int_{Y' \in S^3}
\frac {N(Z')-1} {N(Y'-Z')^2} \cdot
\frac{N(1-iY')}{(Y'-i) \cdot (1-iY')} \cdot Dy \cdot \frac {1}{N(Y'+i)}
\cdot \phi(Y)  \\
=
N(Z'+i) \int_{Y' \in S^3}
\frac {N(Z')-1} {N(Y'-Z')^2} \cdot \frac {Dy}{Y'} \cdot \frac {1}{N(Y'+i)}
\cdot \phi \biggl( \frac{Y'-i}{Y'+i} \biggr)  \\
= -2\pi^2 \phi \bigl( (Z'-i)(Z'+i)^{-1} \bigr) = -2\pi^2 \phi (Z).
\end{multline*}
This proves the first formula for $R=1$.
Then the proof proceeds in exactly the same way as that of
Theorem \ref{Poisson}.
\qed

\begin{rem}
A direct proof of Theorem \ref{Poisson-M-thm} can also be given by adapting
the proof of its Euclidean analogue -- Theorem \ref{Poisson}.
\end{rem}

\subsection{Hydrogen Atom: the Continuous Spectrum}

In this subsection we revisit the spectral decomposition \cite{BI}(II)
of the three-dimensional Laplacian with the Coulomb potential and identify
the positive spectrum of (\ref{Schrodinger}) and the eigenfunctions.
As in the Euclidean case, we start with a function $\phi$ satisfying
the wave equation $\square_{3,1} \phi=0$ and
the Minkowski analogue of the Poisson formula from Theorem \ref{Poisson-M-thm}
with $R=1$.
We need a class of functions of fixed homogeneity degree.
This class can be constructed as follows: pick
$\lambda \in \BB R$, $P \in \BB M$ with $N(P)=0$ and define
$$
\phi_{\lambda}(Y) = \langle Y, -P \rangle^{-1+i\lambda}
= (y^0p^0- y^1p^1 - y^2p^2 - y^3p^3)^{-1+i\lambda}.
$$
Then clearly
$$
(\widetilde{\deg} \phi_{\lambda})(Y) = i\lambda \phi_{\lambda}(Y).
$$
Now, for each $\lambda \in \BB R$, we construct a family of functions
of homogeneity degree $i\lambda$.
Consider the hyperplane section of the light cone
$$
p^0 = 1, \qquad (p^1)^2 + (p^2)^2 + (p^3)^2 = 1
$$
which is a two-dimensional sphere $S^2$.
Choose an orthogonal or orthonormal basis on $S^2$
$$
t_n(P), \qquad n=1,2,3,\dots.
$$
For example, it can consist of spherical harmonics $Y_l^m$'s related to the
matrix coefficients $t^l_{m \, \underline{0}}$'s, where $l=0,1,2,\dots$ and
$-l \le m \le l$ (cf. \cite{V}),
but the specific choice of basis is not essential.
Now we define our family of functions of homogeneity degree $i\lambda$ as
$$
t_n^{\lambda}(Y) = \int_{S^2} t_n(P) \cdot \langle Y, -P \rangle^{-1+i\lambda}
\,d^2P,
\qquad Y \in \BB M.
$$
These functions will play the role of harmonic polynomials
$t^l_{n \, \underline{m}}(X)$ in the Euclidean case.

\begin{rem}
Note that, for a fixed $l$, the span of $t^l_{n \, \underline{m}}(X)$'s
was an irreducible representation of $SO(4) = SU(2) \times SU(2) / \BB Z_2$.
In the same way, the span of $t_n^{\lambda}(Y)$'s is an irreducible
representation of $SO(3,1)$.
In fact these representations have unitary structures and can be
characterized among all unitary representations of $SO(3,1)$ as those having
a fixed vector for $SO(3)$ embedded diagonally in $SO(3,1)$.
These are usually called the type I representations.
\end{rem}

Once we have the space of homogeneous functions spanned by $t_n^{\lambda}(Y)$,
the construction of the eigenfunctions of the hydrogen atom can be done
the same way as in the Euclidean case. Fix a $\lambda \in \BB R$ and let
$\phi$ be a homogeneous solution of $\square_{3,1} \phi=0$ of degree
$i\lambda$. Then we have:
\begin{equation}  \label{Poisson-hydro-M}
\phi \bigl( Z \bigr) = \frac {i\lambda}{2\pi^2} \int_{Y \in H_1}
\frac {\phi(Y)}{N(Y-Z)} \,\frac {Dy}{Y}.
\end{equation}
Fix a $\rho \in \BB R \setminus \{0\}$ and apply the Cayley transform
$$
\pi_l \begin{pmatrix} i\rho & -\rho \\ -1 & i \end{pmatrix}: \HC \to \HC,
\qquad X \mapsto Y= i\frac {X-i\rho}{X+i\rho}.
$$
This transform preserves $\BB M$, and $Y \in H_1$ if and only if $x^0=0$.
Thus we can allow in (\ref{Poisson-hydro-M})
$$
Y = i\frac {X-i\rho}{X+i\rho}, \quad x^0 = 0
\qquad \text{and} \qquad 
Z = i\frac {W-i\rho}{W+i\rho}, \quad w^0 \ne 0.
$$
As in the Euclidean case, we introduce a function $\psi$:
$$
\psi(Z) = \frac 1{(N(Z)-\rho^2)^2} \cdot
\phi \biggl( i\frac {Z-i\rho}{Z+i\rho} \biggr),
\qquad
\psi(W) = \frac 1{(N(W)-\rho^2)^2} \cdot
\phi \biggl( i\frac {W-i\rho}{W+i\rho} \biggr).
$$
Then we apply the Fourier transform and eventually obtain
$$
\biggl( -\frac 12 \Delta + \frac {(i\rho)^2}2 \biggr) \hat\psi_0 (\xi)
= \frac {(i\lambda)(i\rho)}{|\xi|} \cdot \hat\psi_0(\xi).
$$
which is the same as (\ref{Schrod}) in the Euclidean case, except $\rho$ and
$2l+1$ are replaced by $i\rho$ and $i\lambda$ respectively.
We let $\kappa = (i\lambda)(i\rho) = -\lambda\rho$, and we get
\begin{equation}  \label{Schrodinger-M}
- \biggl( \frac 12 \Delta - \frac {\kappa}{|\xi|} \biggr) \hat\psi_0 (\xi)
= E \cdot \hat\psi_0(\xi),
\end{equation}
where
\begin{equation*}
E = \frac {\rho^2}2 = \frac {\kappa^2}{2\lambda^2}.
\end{equation*}
Thus we obtain the eigenfunctions and eigenvalues of (\ref{Schrodinger-M}),
just as in the Euclidean case.
Note that in the Minkowski space the eigenvalues are positive and we
have a continuous spectrum.
Also note that the sign of $\kappa$ depends on the signs of $\lambda$ and
$\rho$. The case $\kappa>0$ is called attractive, and the case $\kappa<0$
is called repulsive.

\section{Middle Series and Quaternionic Analogues of the Second Order Pole}
\label{middle}

\subsection{Reproducing Formula for Functions on $\HC$}

The group $GL(2,\BB H)$ can act on harmonic functions with isolated
singularities by $\pi_l^0$ and $\pi_r^0$ (equations (\ref{left_action})
and (\ref{right_action})).
Differentiating these two actions we obtain two representations of
$\mathfrak{gl}(2, \BB H)$ which preserve ${\cal H}^+$
and agree on $\mathfrak{sl}(2, \BB H)$.
We denote by $(\pi_{lr}^0, {\cal H}^+ \otimes {\cal H}^+)$
the tensor product representation $\pi_l^0 \otimes \pi_r^0$ of
$\mathfrak{gl}(2, \BB H)$. Complexifying, we get a representation of
$\mathfrak{gl}(2, \HC) = \mathfrak{gl}(2, \BB H) \otimes \BB C$,
and we can restrict it to a real form
$\mathfrak{u}(2,2) \subset \mathfrak{gl}(2, \HC)$.
From Theorem \ref{unitary-u(2,2)-action} we know that the representations
$(\pi_{l}^0, {\cal H}^+)$ and $(\pi_{r}^0, {\cal H}^+)$ of
$\mathfrak{u}(2,2)$ are unitary.
It follows that the representation
$(\pi_{lr}^0, {\cal H}^+ \otimes {\cal H}^+)$ of $\mathfrak{u}(2,2)$
is unitary as well and decomposes into a direct sum of irreducible
subrepresentations.

We have a map $M_1$ on ${\cal H}^+ \otimes {\cal H}^+$
which is determined on pure tensors by multiplication
\begin{equation}   \label{M_1}
M_1 \bigl( \phi_1(Z_1) \otimes \phi_2(Z_2) \bigr) = \phi_1(Z) \cdot \phi_2(Z),
\qquad \phi_1,\phi_2 \in {\cal H}^+.
\end{equation}
To make $M_1$ a $\mathfrak{gl}(2, \BB H)$-equivariant map, we define an action
of $GL(2,\BB H)$ on functions on $\HC$ with isolated singularities by
\begin{multline}  \label{1-action}
\rho_1(h): \: F(Z) \mapsto \bigl( \rho_1(h)F \bigr)(Z) =
\frac {F \bigl( (aZ+b)(cZ+d)^{-1} \bigr)}{N(cZ+d) \cdot N(a'-Zc')},  \\
h = \begin{pmatrix} a' & b' \\ c' & d' \end{pmatrix},\:
h^{-1} = \begin{pmatrix} a & b \\ c & d \end{pmatrix} \in GL(2, \BB H);
\end{multline}
since
$$
(aZ+b)(cZ+d)^{-1} = (a'-Zc')^{-1}(-b'+Zd'),
$$
$\rho_1$ is a well-defined action.
Differentiating this action and complexifying,
we obtain an action of $\mathfrak{gl}(4, \BB C)$.
Recall that $\partial = \begin{pmatrix} \partial_{11} & \partial_{21} \\
\partial_{12} & \partial_{22} \end{pmatrix} = \frac 12 \nabla$.

\begin{lem}  \label{rho_1-algebra-action}
The Lie algebra action $\rho_1$ of $\mathfrak{gl}(4,\BB C)$ on
functions on $\HC$ is given by
\begin{align*}
\rho_1 \begin{pmatrix} A & 0 \\ 0 & 0 \end{pmatrix} &:
F \mapsto \tr \bigl( A \cdot (-Z \cdot \partial F - F) \bigr)  \\
\rho_1 \begin{pmatrix} 0 & B \\ 0 & 0 \end{pmatrix} &:
F \mapsto \tr \bigl( B \cdot (-\partial F ) \bigr)  \\
\rho_1 \begin{pmatrix} 0 & 0 \\ C & 0 \end{pmatrix} &:
F \mapsto \tr \Bigl( C \cdot \bigl(
Z \cdot (\partial F) \cdot Z + 2 ZF \bigr) \Bigr)  \\
&: F \mapsto \tr \Bigl( C \cdot \bigl(
Z \cdot \partial (ZF) \bigr) \Bigr)  \\
\rho_1 \begin{pmatrix} 0 & 0 \\ 0 & D \end{pmatrix} &:
F \mapsto \tr \Bigl( D \cdot \bigl( (\partial F) \cdot Z + F \bigr) \Bigr)  \\
&: F \mapsto \tr \Bigl( D \cdot \bigl( \partial (ZF) - F \bigr) \Bigr).
\end{align*}
\end{lem}

This lemma implies that $\mathfrak{gl}(4, \BB C)$ preserves the space
$$
\Zh^+ = \{\text{polynomial functions on $\HC$}\} = \BB C[z_{ij}].
$$
Theorem \ref{JV2-thm} implies that $(\rho_1, \Zh^+)$ is irreducible.
It is the first example of a middle series representation.
In general, the middle series is formed from the lowest component in the
tensor product of dual representations from the most degenerate series.
The middle series is another degenerate series of irreducible unitary
representations of $SU(2,2)$ (see \cite{KS} for a complete classification
of the irreducible unitary representations).

Define
$$
\Zh =
\{\text{polynomial functions on $\HC^{\times}$}\}  \\
= \BB C[z_{ij},N(Z)^{-1}]
$$
and
\begin{multline*}
\Zh^- =
\{\text{polynomials on $\HC^{\times}$ regular at infinity}\}  \\
= \{ F \in \BB C[z_{ij}, N(Z)^{-1}] ;\: \text{$F$ is regular at infinity} \}.
\end{multline*}
(Functions $F \in \Zh$ are defined to be regular at infinity as in
Definition \ref{reg-infinity-def}.)
Note that $(\rho_1, \Zh^+) \simeq (\rho_1,\Zh^-)$ and
$\Zh^+ \oplus \Zh^-$ is a proper subspace of $\Zh$.

\begin{prop}  \label{pairing_n=1}
The representation $(\rho_1, \Zh)$ of $\mathfrak{gl}(4, \BB C)$
has a non-degenerate symmetric bilinear pairing
\begin{equation}  \label{U2-pairing}
\langle F_1,F_2 \rangle_1 =
\frac i{8\pi^3} \int_{U(2)_R} F_1(Z) \cdot F_2(Z) \,dZ^4,
\qquad F_1, F_2 \in \Zh,
\end{equation}
where $R>0$ and
$$
U(2)_R = \{ RZ ;\: Z \in U(2) \}.
$$
This bilinear pairing is $\mathfrak{u}(2,2)$-invariant and
independent of the choice of $R>0$.
\end{prop}

\pf
Since the integrand is a closed form, the pairing is independent of the
choice of $R>0$.
We prove the invariance of the bilinear pairing on $(\rho_1, \Zh)$
by showing that, for all $h \in U(2,2)$ sufficiently
close to the identity element, we have
$$
\langle F_1, F_2 \rangle_1 =
\langle \rho_1(h) F_1,\rho_1(h) F_2 \rangle_1,
\qquad F_1,F_2 \in \Zh.
$$
If $h^{-1}= \begin{pmatrix} a & b \\ c & d \end{pmatrix} \in U(2,2)$,
then $h = \begin{pmatrix} a^* & -c^* \\ -b^* & d^* \end{pmatrix}$.
Writing $\tilde Z = (aZ+b)(cZ+d)^{-1}$ and using Lemma \ref{Jacobian_lemma},
we obtain:
\begin{multline*}
-8\pi^3i \cdot \langle \rho_1(h) F_1,\rho_1(h) F_2 \rangle_1  \\
= \int_{U(2)} \frac {F_1(\tilde Z)}{N(cZ+d) \cdot N(a^*+Zb^*)} \cdot
\frac {F_2(\tilde Z)}{N(cZ+d) \cdot N(a^*+Zb^*)} \,dZ^4  \\
= \int_{U(2)} F_1 (\tilde Z) \cdot F_2(\tilde Z) \,d \tilde Z^4
= -8\pi^3i \cdot \langle F_1, F_2 \rangle_1.
\end{multline*}

It remains to prove that the pairing is non-degenerate.
For $F_1 \in \Zh$, define
$$
F_2(Z) =
\biggl( \B{\rho_1 \begin{pmatrix} 0 & 1 \\ 1 & 0 \end{pmatrix} F_1} \biggr)
(Z^*) = \frac 1{N(Z)^2} \cdot \B{F_1} \bigl( (Z^*)^{-1} \bigr).
$$
Then
$$
\langle F_1, F_2 \rangle_1 = 
\frac i{8\pi^3} \int_{U(2)} |F_1 (Z)|^2  \,\frac{dZ^4}{N(Z)^2} \quad >0,
$$
unless $F_1 =0$.
\qed

Next we prove a reproducing formula for functions on $\HC$.

\begin{thm}  \label{2pole1}
For all $F(Z) \in \BB C[z_{ij}]$ with no differential conditions
imposed whatsoever, all $R>0$ and all $W \in \BB H$ with $N(W)<R^2$, we have
$$
F(W) = \Bigl\langle F(Z), \frac 1{N(Z-W)^2} \Bigr\rangle_1 =
\frac i{8\pi^3} \int_{Z \in U(2)_R} \frac {F(Z)}{N(Z-W)^2} \,dZ^4.
$$
\end{thm}

\begin{lem}
The map on $\Zh^+$
\begin{equation}  \label{2poles-lemma}
F(Z) \mapsto \int_{Z \in U(2)} \frac {F(Z)}{N(Z-W)^2} \,dZ^4
\end{equation}
is $U(2,2)$-equivariant.
\end{lem}

\pf
If $h^{-1}= \begin{pmatrix} a & b \\ c & d \end{pmatrix} \in U(2,2)$,
then $h = \begin{pmatrix} a^* & -c^* \\ -b^* & d^* \end{pmatrix}$.
Writing
$\tilde Z = (aZ+b)(cZ+d)^{-1}$, $\tilde W = (aW+b)(cW+d)^{-1}$
and using Corollary \ref{N(Z-W)} and Lemma \ref{Jacobian_lemma}, we obtain:
\begin{multline*}
\frac {\rho_1(h)F(Z)}{N(Z-W)^2} \,dZ^4
= \frac {F (\tilde Z)} {N(cZ+d) \cdot N(a^*+Zb^*) \cdot N(Z-W)^2} \,dZ^4  \\
= \frac {F (\tilde Z)}
{N(cZ+d)^2 \cdot N(a^*+Zb^*)^2 \cdot N( \tilde Z- \tilde W)^2
\cdot N(cW+d) \cdot N(a^*+Wb^*)} \,dZ^4  \\
= \frac {F (\tilde Z)}
{N( \tilde Z- \tilde W)^2 \cdot N(cW+d) \cdot N(a^*+Wb^*)} \,d \tilde Z^4.
\end{multline*}
The $U(2,2)$-equivariance then follows.
\qed

\noindent {\it Proof of Theorem \ref{2pole1}.}
The group $U(2)\times U(2)$ is a maximal compact subgroup of $U(2,2)$,
and the map (\ref{2poles-lemma}) must preserve the space of
$U(2)\times U(2)$-finite vectors.
Hence by Schur's Lemma there exists a $\lambda \in \BB C$ such that the map 
(\ref{2poles-lemma}) is given by multiplication by $\lambda$.
To pin down the value of $\lambda$ we substitute $P(Z) \equiv 1$ and $W=0$,
and we immediately see that $\lambda=-8\pi^3i$.
\qed

\begin{rem}
Another proof of Theorem \ref{2pole1} can be given using the matrix
coefficient expansions from Proposition \ref{N-square-expansion}.
\end{rem}

We have the following Minkowski counterparts of Proposition \ref{pairing_n=1}
and Theorem \ref{2pole1}. They are proved the same way
Theorems \ref{Fueter-hyperboloid} and \ref{Poisson-M-thm} were.

We denote by $\Zh (\BB M)$ the space of functions $F$ on $\BB M$
such that $\rho_1(\gamma^{-1})F$ extends holomorphically to
an open neighborhood of $U(2)$.
We also define the spaces $\Zh (\BB M)^{\pm}$ of holomorphic functions $F$
on $\BB T^{\pm}$ such that $\rho_1 (\gamma^{-1})F$
extends to a holomorphic function defined on a neighborhood of the closure
of $\BB D^{\pm}$ and, in the case of $\Zh (\BB M)^-$, regular at infinity.

\begin{prop}
The representation $(\rho_1, \Zh(\BB M))$ of $\mathfrak{gl}(4, \BB C)$
has a non-degenerate symmetric bilinear pairing
\begin{equation}  \label{M-pairing}
\langle F_1,F_2 \rangle_1 =
\frac i{8\pi^3} \int_{\BB M} F_1(Z) \cdot F_2(Z) \,dZ^4,
\qquad F_1, F_2 \in \Zh(\BB M),
\end{equation}
which is $\mathfrak{u}(2,2)'$-invariant.
\end{prop}

\begin{thm}  \label{2pole1-M}
For all $F \in \Zh (\BB M)^+$ and all $W \in \BB T^+$, we have:
$$
F(W) = \Bigl\langle F(Z), \frac 1{N(Z-W)^2} \Bigr\rangle_1 =
\frac i{8\pi^3} \int_{Z \in \BB M} \frac {F(Z)}{N(Z-W)^2} \,dZ^4.
$$

Similarly, for all $F \in \Zh (\BB M)^-$ and all $W \in \BB T^-$, we have:
$$
F(W) = \Bigl\langle F(Z), \frac 1{N(Z-W)^2} \Bigr\rangle_1 =
\frac i{8\pi^3} \int_{Z \in \BB M} \frac {F(Z)}{N(Z-W)^2} \,dZ^4.
$$
\end{thm}

We conclude this subsection with a discussion of Hardy spaces.
The main references are \cite{VR}, \cite{FK} and \cite{AU}.
First we define the Hardy spaces of holomorphic functions on $\BB D^{\pm}$:
$$
\ha (\BB D^+) = \biggl\{ f \in {\cal O}(\BB D^+) ;\:
\sup_{0<r<1} \int_{U(2)} |f(rZ)|^2 \,\frac{dZ^4}{N(Z)^2} < \infty \biggr\},
$$
and $\ha (\BB D^-)$ is defined similarly; these spaces are completions of
$\Zh^+$ and $\Zh^-$.
Then we define
$$
\ha_{harm} (\BB D^+) = \biggl\{ f \in {\cal O}(\BB D^+) ;\:
\square f =0,\: \sup_{0<r<1} \int_{SU(2)} 
(\widetilde{\deg}f)(rZ) \cdot \B{f(rZ)} \,dS < \infty \biggr\},
$$
and similarly for $\ha_{harm} (\BB D^-)$; these spaces are completions of
${\cal H}^+$ and ${\cal H}^-$.
Next we define the Hardy spaces on $\BB T^{\pm}$:
$$
\ha (\BB T^+) = \biggl\{ f \in {\cal O}(\BB T^+) ;\:
\sup_{W \in \C^+}
\int_{Y \in \BB M} |f(Y+iW)|^2 \,dY^4 < \infty \biggr\},
$$
and similarly for $\ha (\BB T^-)$; these spaces are completions of
$\Zh(\BB M)^+$ and $\Zh(\BB M)^-$.
Then we define
$$
\ha_{harm} (\BB T^+) = \biggl\{ f \in {\cal O}(\BB T^+) ;\:
\square f =0,\: \sup_{W \in \C^+} \int_{H'_1}
(\widetilde{\deg}f)(Y+iW) \cdot \B{f(Y+iW)} \,\frac{dS}{\|Y\|}
< \infty \biggr\},
$$
and similarly for $\ha_{harm} (\BB T^-)$; these spaces are completions of
${\cal H}(\BB M)^+$ and ${\cal H}(\BB M)^-$.

\begin{thm} [\cite{FK}, \cite{AU}]  \label{Hardy-iso}
The Cayley transform provides isomorphisms of Hardy spaces
$$
\rho_1(\gamma): \ha (\BB D^{\pm}) \simeq \ha (\BB T^{\pm}),  \qquad
\pi^0_l(\gamma): \ha_{harm} (\BB D^{\pm}) \simeq \ha_{harm} (\BB T^{\pm}).
$$
And the Fourier transform provides isomorphisms
$$
\ha (\BB T^{\pm}) \simeq L^2( \C^+, dZ^4),
\qquad
\ha_{harm} (\BB T^{\pm}) \simeq L^2( \partial\C^+, dS/\|Y\|),
$$
where $\C^+$ is the open cone in $\BB M$ defined by (\ref{C^+}).
\end{thm}

The Poisson and reproducing formulas on $\BB H$ and $\BB M$ naturally extend
to the Hardy spaces on $\BB D^{\pm}$ and $\BB T^{\pm}$ respectively.

\subsection{Quaternionic Analogue of the Cauchy Formula
for the Second Order Pole}  \label{repro2-section}

In this subsection we prove a quaternionic analogue of the Cauchy
formula for the second order pole (\ref{2pole-intro})
for quaternionic-valued functions.

We introduce the following notations:
$$
{\cal V}^+ =
\{ \text{polynomial left-regular spinor-valued functions on $\HC$} \},
$$
$$
{\cal V}'^+ =
\{ \text{polynomial right-regular spinor-valued functions on $\HC$} \},
$$
$$
{\cal W}^+ = \{\text{$\HC$-valued polynomial functions on $\HC$}\}
= \HC [z_{ij}].
$$
The group $GL(2,\BB H)$ acts on the spaces of left- and right-regular
spinor-valued functions via $\pi_l$ and $\pi_r$ given by
(\ref{spinor_left_action}) and (\ref{spinor_right_action}).
Differentiating the actions $\pi_l$ and $\pi_r$ of $GL(2,\BB H)$ on
spinor-valued regular functions we get representations of the Lie algebra
$\mathfrak{gl}(2,\BB H)$ which we still denote by $\pi_l$ and $\pi_r$.
We denote by $(\pi_{lr}, {\cal V}^+ \otimes {\cal V}'^+)$
the tensor product representation $\pi_l \otimes \pi_r$ of
$\mathfrak{gl}(2, \BB H)$.
Complexifying, we get a representation of
$\mathfrak{gl}(2, \HC) = \mathfrak{gl}(2, \BB H) \otimes \BB C$,
and we can restrict it to a real form
$\mathfrak{u}(2,2) \subset \mathfrak{gl}(2, \HC)$.
From Theorem \ref{unitary-u(2,2)-action} we know that the representations
$(\pi_{l}, {\cal V}^+)$ and $(\pi_{r}, {\cal V}'^+)$ of $\mathfrak{u}(2,2)$
are unitary.
It follows that the representation
$(\pi_{lr}, {\cal V}^+ \otimes {\cal V}'^+)$ of $\mathfrak{u}(2,2)$ is also
unitary and decomposes into a direct sum of irreducible subrepresentations.
The lowest component ${\cal W}^+$ is the second example of a middle series
representation.

We have a map $M: {\cal V}^+ \otimes {\cal V}'^+ \to {\cal W}^+$
which is determined on pure tensors by
$$
M \bigl( f(Z_1) \otimes g(Z_2) \bigr) = f(Z) \cdot g(Z),
\qquad f \in {\cal V}^+,\: g \in {\cal V}'^+.
$$
This map $M$ becomes $\mathfrak{gl}(2, \BB H)$-equivariant if we define
an action of $GL(2,\BB H)$ on $\HC$-valued functions on $\BB H$ with
isolated singularities by
\begin{multline}  \label{action}
\rho_2(h): \: F(Z) \mapsto \bigl( \rho_2(h)F \bigr)(Z) =
\frac {(cZ+d)^{-1}}{N(cZ+d)} \cdot F \bigl( (aZ+b)(cZ+d)^{-1} \bigr) \cdot
\frac {(a'-Zc')^{-1}}{N(a'-Zc')},  \\
h = \begin{pmatrix} a' & b' \\ c' & d' \end{pmatrix},\:
h^{-1} = \begin{pmatrix} a & b \\ c & d \end{pmatrix} \in GL(2, \BB H);
\end{multline}
since $(aZ+b)(cZ+d)^{-1} = (a'-Zc')^{-1}(-b'+Zd')$,
$\rho_2$ is a well-defined action.
Differentiating this action and complexifying, we obtain an action of
$\mathfrak{gl}(4, \BB C)$ which preserves the space ${\cal W}^+$.
Then $M$ becomes an intertwining operator between
$\mathfrak{gl}(4, \BB C)$-representations
$(\pi_{lr}, {\cal V}^+ \otimes {\cal V}'^+)$ and $(\rho_2, {\cal W}^+)$.
By Theorem \ref{JV2-thm}, $(\rho_2, {\cal W}^+)$ is irreducible unitary
on $\mathfrak{u}(2,2)$.
We will describe the unitary structure at the end of his subsection.

\begin{lem}
For all differentiable $\HC$-valued functions $F$,
$$
\nabla_Z \frac {F(W)}{N(Z-W)^2} \nabla_Z - \square_Z \frac {F(W)^+}{N(Z-W)^2}
= 24 \cdot
\frac {(Z-W)^{-1}}{N(Z-W)} \cdot F(W) \cdot \frac {(Z-W)^{-1}}{N(Z-W)}.
$$
\end{lem}

\pf
We prove the lemma by direct computation using
$\nabla = 2\begin{pmatrix} \partial_{11} & \partial_{21} \\
\partial_{12} & \partial_{22} \end{pmatrix}$.
First we verify
$$
\nabla_Z \frac 1{N(Z-W)^2} = \frac 1{N(Z-W)^2} \nabla_Z
= - 4 \frac {(Z-W)^+}{N(Z-W)^3}.
$$
Next we compute
\begin{multline*}
- \nabla_Z \frac {F(W) \cdot (Z-W)^+}{N(Z-W)^3}  \\
= 6 \frac {(Z-W)^+ \cdot F(W) \cdot (Z-W)^+}{N(Z-W)^4}
- \frac 1{N(Z-W)^3} \cdot \nabla_Z \bigl( F(W)(Z-W)^+ \bigr),
\end{multline*}
\begin{multline*}
\nabla_Z
\biggr[ \begin{pmatrix} F_{11} & F_{12} \\ F_{21} & F_{22} \end{pmatrix}
\begin{pmatrix} z_{22} & -z_{12} \\ -z_{21} & z_{11} \end{pmatrix} \biggr]  \\
=
2 \begin{pmatrix} \partial_{11} & \partial_{21} \\
\partial_{12} & \partial_{22} \end{pmatrix}
\begin{pmatrix} F_{11}z_{22} - F_{12}z_{21} & -F_{11}z_{12} + F_{12}z_{11} \\
F_{21}z_{22} - F_{22}z_{21} & -F_{21}z_{12} + F_{22}z_{11} \end{pmatrix}  \\
=
2 \begin{pmatrix} -F_{22} & F_{12} \\ F_{21} & -F_{11} \end{pmatrix}
= -2 F(W)^+.
\end{multline*}
On the other hand,
$\square = 4 (\partial_{11}\partial_{22} - \partial_{12}\partial_{21})$, so
\begin{multline*}
\square_Z \frac 1{N(Z-W)^2}
=
-8 \partial_{11} \frac {(z_{11}-w_{11})}{N(Z-W)^3}
-8 \partial_{12} \frac {(z_{12}-w_{12})}{N(Z-W)^3}  \\
=
- \frac {16}{N(Z-W)^3} + 24 \frac
{(z_{11}-w_{11})(z_{22}-w_{22}) - (z_{12}-w_{12})(z_{21}-w_{21})}{N(Z-W)^4}
= \frac {8}{N(Z-W)^3},
\end{multline*}
and lemma follows.
\qed

Let $\M$ denote the differential operator on ${\cal W}^+$ defined by
$$
\M F = \nabla F \nabla - \square F^+.
$$
Theorem \ref{2pole1} holds for $\BB C$-valued functions and can also
be applied to the $\HC$-valued functions in ${\cal W}^+$.
Applying $\M$ to both sides of the equation, we obtain another quaternionic
analogue of the Cauchy formula with the second order pole:

\begin{thm}  \label{2pole2}
For all $F(Z) \in \BB \HC[z_{ij}]$, for all $R>0$ and all $W \in \BB H$
with $N(W)<R^2$, we have
$$
\M F(W) = \frac {3i}{\pi^3} \int_{Z \in U(2)_R}
\frac {(Z-W)^{-1}}{N(Z-W)} \cdot F(Z) \cdot \frac {(Z-W)^{-1}}{N(Z-W)} \,dZ^4.
$$
\end{thm}

Define another meromorphic action of $GL(2,\BB H)$ on
$\HC$-valued functions on $\BB H$:
\begin{multline*}
\rho'_2(h): \: F(Z) \mapsto \bigl( \rho'_2(h)F \bigr)(Z) =
\frac {(a'-Zc')}{N(a'-Zc')} \cdot F \bigl( (aZ+b)(cZ+d)^{-1} \bigr)
\cdot \frac {(cZ+d)}{N(cZ+d)},  \\
h = \begin{pmatrix} a' & b' \\ c' & d' \end{pmatrix},\:
h^{-1} = \begin{pmatrix} a & b \\ c & d \end{pmatrix} \in GL(2, \BB H).
\end{multline*}
Differentiating this action and complexifying,
we obtain a representation $\rho'_2$ of $\mathfrak{gl}(4, \BB C)$.
This action preserves ${\cal W}^+$, and we denote by ${\cal W}'^+$
the space ${\cal W}^+$ with $\rho'_2$ action of $\mathfrak{gl}(4, \BB C)$.

\begin{prop}  \label{Mx-intertwiner}
The differential operator $\M$ is an intertwining operator of
$\mathfrak{gl}(4, \BB C)$-representations
$(\rho'_2, {\cal W}'^+) \to (\rho_2, {\cal W}^+)$.
\end{prop}

\pf
It is sufficient to prove that $\M$ intertwines the actions
of $\mathfrak{u}(2,2)$.
And to do this, it is sufficient to check that $\M$ intertwines the actions
of all $h \in U(2,2)$ sufficiently close to the identity element. We have:
\begin{multline*}
\M \bigl( \rho'_2(h)F \bigr)(W) = \frac {3i}{\pi^3} \int_{U(2)}
\frac {(Z-W)^{-1} \cdot \bigl( \rho'_2(h)F \bigr)(Z) \cdot (Z-W)^{-1}}{N(Z-W)^2} \,dZ^4  \\
= \frac {3i}{\pi^3} \int_{U(2)}
\frac {(Z-W)^{-1} \cdot (a^*+Zb^*) \cdot F(\tilde Z) \cdot (cZ+d)
\cdot (Z-W)^{-1}} {N(a^*+Zb^*) \cdot N(cZ+d) \cdot N(Z-W)^2} \,dZ^4  \\
= \frac {3i}{\pi^3} \int_{U(2)}
\frac {(cW+d)^{-1} \cdot (\tilde Z - \tilde W)^{-1} \cdot F(\tilde Z)
\cdot (\tilde Z - \tilde W)^{-1} \cdot (a^*+Wb^*)^{-1}}
{N(cZ+d)^2 \cdot N(a^*+Zb^*)^2 \cdot N(cW+d) \cdot N(a^*+Wb^*) \cdot
N(\tilde Z - \tilde W)^2} \,dZ^4  \\
= \frac {3i}{\pi^3} \int_{U(2)}
\frac {(cW+d)^{-1} \cdot (\tilde Z - \tilde W)^{-1} \cdot F(\tilde Z)
\cdot (\tilde Z - \tilde W)^{-1} \cdot (a^*+Wb^*)^{-1}}
{N(cW+d) \cdot N(a^*+Wb^*) \cdot N(\tilde Z - \tilde W)^2} \,d\tilde Z^4  \\
= (\rho_2(h) \M F)(W).
\end{multline*}
\qed

\begin{rem}
It is also true that the differential operator $\M$ intertwines the actions
$\rho'_2$ and $\rho_2$ on {\em all} smooth $\HC$-valued functions on $\HC$,
not just on the polynomial ones. This can be verified directly without
the use of Theorem \ref{2pole2} and requires calculations involving
Lie algebra actions.
\end{rem}

Define
$$
{\cal W} =
\{\text{$\HC$-valued polynomials on $\HC^{\times}$}\}  \\
= \{ F \in \BB H \otimes \BB C[z_{ij},N(Z)^{-1}] \}
$$
and
\begin{multline*}
{\cal W}^- =
\{\text{$\HC$-valued polynomials on $\HC^{\times}$ regular at infinity}\}  \\
= \{ F \in \BB H \otimes \BB C[z_{ij},N(Z)^{-1}] ;\:
\text{$F$ is regular at infinity} \}.
\end{multline*}
Note that ${\cal W}^+ \oplus {\cal W}^-$ is a proper subspace of ${\cal W}$.
The Lie algebra $\mathfrak{gl}(4, \BB C)$ acts on ${\cal W}$ and
${\cal W}^-$ by $\rho_2$.
We denote by ${\cal W}'$ and ${\cal W}'^-$ the same spaces with $\rho'_2$
action.

\begin{prop} \label{W-W'_pairing}
There is a $\mathfrak{u}(2,2)$-invariant bilinear pairing between the
representations $(\rho_2, {\cal W})$ and $(\rho'_2, {\cal W}')$
of $\mathfrak{gl}(4, \BB C)$ given by
$$
\langle F_1,F_2 \rangle =
\frac i{8\pi^3} \int_{U(2)_R} \tr \bigl( F_1(Z) \cdot F_2(Z) \bigr) \,dZ^4,
\qquad F_1 \in {\cal W}, \: F_2 \in {\cal W}'.
$$
This pairing is independent of the choice of $R>0$.
Moreover, for each $F_1 \in {\cal W}^+$, $F_1 \ne 0$,
there is an $F_2 \in {\cal W}'^-$ such that $\langle F_1,F_2 \rangle \ne 0$.
\end{prop}

\pf
Since the integrand is a closed form, the pairing is independent of
the choice of $R>0$.
We prove the invariance of the bilinear pairing on ${\cal W} \times {\cal W}'$
by showing that, for all $h \in U(2,2)$ sufficiently
close to the identity element, we have
$$
\langle F_1, F_2 \rangle =
\langle \rho_2(h) F_1, \rho'_2(h) F_2 \rangle,
\qquad F_1 \in {\cal W}, \: F_2 \in {\cal W}'.
$$
If $h^{-1}= \begin{pmatrix} a & b \\ c & d \end{pmatrix} \in U(2,2)$,
then $h = \begin{pmatrix} a^* & -c^* \\ -b^* & d^* \end{pmatrix}$.
Writing $\tilde Z = (aZ+b)(cZ+d)^{-1}$ and using Lemma \ref{Jacobian_lemma},
we obtain:
\begin{multline*}
-8\pi^3i \cdot \langle \rho_2(h) F_1, \rho'_2(h) F_2 \rangle =  \\
\int_{U(2)} \tr \biggl(
\frac {(cZ+d)^{-1} \cdot F_1(\tilde Z) \cdot (a^*+Zb^*)^{-1}}
{N(cZ+d) \cdot N(a^*+Zb^*)} \cdot
\frac {(a^*+Zb^*) \cdot F_2(\tilde Z) \cdot (cZ+d)}{N(cZ+d) \cdot N(a^*+Zb^*)}
\biggr) \,dZ^4  \\
= \int_{U(2)} \tr \bigl( F_1 (\tilde Z) \cdot F_2(\tilde Z) \bigr)
\,d\tilde Z^4
= -8\pi^3i \cdot \langle F_1, F_2 \rangle.
\end{multline*}

To prove the second statement, pick an $F_1 \in {\cal W}^+$, $F_1 \ne 0$,
and define $F'_1(Z) = - Z \cdot F_1(Z) \cdot Z$,
$$
F_2(Z) =
\biggl( \B{\rho'_2 \begin{pmatrix} 0 & 1 \\ 1 & 0 \end{pmatrix} F'_1} \biggr)
(Z^*) = \frac 1{N(Z)^2} \cdot \B{F_1} \bigl( (Z^*)^{-1} \bigr)
\qquad \in {\cal W}^-.
$$
Then
$$
\langle F_1, F_2 \rangle = 
\frac i{8\pi^3} \int_{U(2)} |F_1 (Z)|^2  \,\frac {dZ^4}{N(Z)^2} \quad >0.
$$
\qed

Since $(\rho_2,{\cal W}^+)$ is irreducible, Proposition \ref{Mx-intertwiner}
implies ${\cal W}^+ \simeq {\cal W}'^+/\ker \M$.
Thus we obtain a $\mathfrak{u}(2,2)$-invariant pairing on
${\cal W}'^+ \times {\cal W}'^-$:
\begin{multline} \label{W-pairing}
\langle F_1, F_2 \rangle_{{\cal W}^+} = 
\langle F_1, \M F_2 \rangle = \langle \M F_1, F_2 \rangle  \\
= \frac {-3}{8\pi^6} \int_{W \in U(2)_r} \int_{Z \in U(2)_R}
\tr \biggl( F_1(W) \cdot \frac {(Z-W)^{-1}}{N(Z-W)} \cdot F_2(Z)
\cdot \frac {(Z-W)^{-1}}{N(Z-W)} \biggr) \,dZ^4 \,dW^4,
\end{multline}
where $F_1 \in {\cal W}'^+$, $F_2 \in {\cal W}'^-$and $R>r>0$.
This pairing is zero on the kernels of $\M$ and descends to a non-degenerate
pairing on
$({\cal W}'^+/ \ker \M) \times ({\cal W}'^-/ \ker \M) \simeq
{\cal W}^+ \times {\cal W}^-$.

Next we state the Minkowski counterpart of Theorem \ref{2pole2}.

\begin{thm}
For all $F \in \Zh (\BB M)^+$ and all $W \in \BB T^+$, we have:
$$
\M F(W) = \frac {3i}{\pi^3} \int_{Z \in \BB M}
\frac {(Z-W)^{-1}}{N(Z-W)} \cdot F(Z) \cdot \frac {(Z-W)^{-1}}{N(Z-W)} \,dZ^4.
$$

Similarly, for all $F \in \Zh (\BB M)^-$ and all $W \in \BB T^-$, we have:
$$
\M F(W) = \frac {3i}{\pi^3} \int_{Z \in \BB M}
\frac {(Z-W)^{-1}}{N(Z-W)} \cdot F(Z) \cdot \frac {(Z-W)^{-1}}{N(Z-W)} \,dZ^4.
$$
\end{thm}

We can also consider the Hardy spaces which are completions of the spaces
of left- and right-regular functions on $\BB D^{\pm}$ and $\BB T^{\pm}$:
$$
\ha_{left-reg} (\BB D^+, \BB S) = \biggl\{ f \in {\cal O}(\BB D^+, \BB S) ;\:
\nabla^+ f =0,\: \sup_{0<r<1} \int_{SU(2)} |f(rZ)|^2 \,dS < \infty \biggr\},
$$
and similarly for $\ha_{left-reg} (\BB D^-, \BB S)$ and
$\ha_{right-reg} (\BB D^{\pm}, \BB S')$; then
$$
\ha_{left-reg} (\BB T^+, \BB S) = \biggl\{ f \in {\cal O}(\BB T^+, \BB S) ;\:
\nabla^+ f =0,\: \sup_{W \in \C^+} \int_{H'_1}
|f(Y+iW)|^2 \,dS < \infty \biggr\},
$$
and similarly for $\ha_{left-reg} (\BB T^-, \BB S)$ and
$\ha_{right-reg} (\BB T^{\pm}, \BB S')$.
In light of Theorem \ref{Hardy-iso}, it is natural to expect the following
isomorphisms
$$
\pi_l(\gamma) : \ha_{left-reg} (\BB D^{\pm}, \BB S)
\simeq \ha_{left-reg} (\BB T^{\pm}, \BB S),
$$
$$
\pi_r(\gamma) : \ha_{right-reg} (\BB D^{\pm}, \BB S')
\simeq \ha_{right-reg} (\BB T^{\pm}, \BB S').
$$
The Fueter formulas on $\BB H$ and $\BB M$ naturally extend to these Hardy
spaces on $\BB D^{\pm}$ and $\BB T^{\pm}$ respectively.

\subsection{Maxwell Equations in Vacuum}

In this subsection we identify the operator $\M$ with Maxwell equations
for the gauge potential.
Recall the classical Maxwell equations:
$$
\begin{array} {lcl}
\grad \cdot \overrightarrow{B} = 0 & \qquad &
\grad \cdot \overrightarrow{E} = 0  \\
\grad \times \overrightarrow{B} =
\frac {\partial \overrightarrow{E}}{\partial t} & \qquad &
\grad \times \overrightarrow{E} =
-\frac {\partial \overrightarrow{B}}{\partial t},
\end{array}
$$
where $\overrightarrow{B}$ and $\overrightarrow{E}$ are three-dimensional
vector functions on $\BB R^4$ (called respectively the magnetic and electric
fields) and
$\grad =\bigl( \frac{\partial}{\partial y^1}, \frac{\partial}{\partial y^2},
\frac{\partial}{\partial y^3} \bigr)$, as usual.
One can look for solutions in the form of the gauge field
$A=(A_0,\overrightarrow{A})$, where $\overrightarrow{A}$ is a three-dimensional
vector function on $\BB R^4$, by taking
\begin{equation}  \label{gauge-field}
\overrightarrow{B} = \grad \times \overrightarrow{A},
\qquad
\overrightarrow{E} =
\grad A_0 - \frac {\partial \overrightarrow{A}}{\partial t}.
\end{equation}
This presentation is gauge invariant, i.e. $\overrightarrow{B}$
and $\overrightarrow{E}$ stay unchanged if $A$ is replaced with
$A'=A+ \bigl( \partial \phi /\partial t, \grad \phi \bigr)$,
for any scalar-valued function $\phi$ on $\BB R^4$.
Under the presentation (\ref{gauge-field}), the equations
$\grad \cdot \overrightarrow{B} = 0$ and
$\grad \times \overrightarrow{E} =
-\frac {\partial \overrightarrow{B}}{\partial t}$ are satisfied automatically.
Thus we end up with two equations:
\begin{equation}  \label{Maxwell}
\grad \cdot \biggl(
\grad A_0 - \frac {\partial \overrightarrow{A}}{\partial t} \biggr) =0,
\qquad
\grad \times \bigl( \grad \times \overrightarrow{A} \bigr) =
\frac {\partial}{\partial t}
\biggl( \grad A_0 - \frac {\partial \overrightarrow{A}}{\partial t} \biggr).
\end{equation}

Now we would like to compare these equations with the quaternionic
equation $\M \tilde A =0$, $\tilde A \in {\cal W}^+$.
We identify the $y^0$ coordinate with time $t$, write
$\overrightarrow{A}$ as $A_1e_1 + A_2e_2 + A_3e_3$ and introduce
$$
\tilde A = -A_0 \tilde e_0 + A_1e_1 + A_2e_2 + A_3e_3
= - \tilde e_0 A_0 + \overrightarrow{A}.
$$
We chose this $\tilde A$ to represent $(A_0,\overrightarrow{A})$ so that
replacing $(A_0,\overrightarrow{A})$ with
$(A_0,\overrightarrow{A}) + \bigl(\partial \phi /\partial t, \grad \phi \bigr)$
corresponds to replacing $\tilde A$ with
$$
- \biggl( A_0 + \frac {\partial \phi}{\partial y^0} \biggr) \tilde e_0
+ \biggl( A_1 + \frac {\partial \phi}{\partial y^1} \biggr) e_1
+ \biggl( A_2 + \frac {\partial \phi}{\partial y^2} \biggr) e_2
+ \biggl( A_3 + \frac {\partial \phi}{\partial y^3} \biggr) e_3
= \tilde A + \nabla^+_{\BB M} \phi,
$$
for all $\phi: \BB M \to \BB C$, and
$$
\M \nabla^+_{\BB M} \phi =
\nabla_{\BB M} \nabla^+_{\BB M} \phi \nabla_{\BB M}
- \square_{3,1} (\nabla^+_{\BB M} \phi)^+
= \nabla_{\BB M} (\square_{3,1} \phi)
- \square_{3,1} (\nabla_{\BB M} \phi) =0.
$$
This way the equation $\M \tilde A=0$ becomes gauge invariant.
We need to expand
$$
\M \tilde A =
\nabla_{\BB M} \tilde A \nabla_{\BB M} - \square_{3,1} \tilde A^+ =0.
$$
We have
$\nabla_{\BB M} = i\frac {\partial}{\partial y^0} - \grad$, and so
$$
\tilde A \nabla_{\BB M}
= \bigl( iA_0 + \overrightarrow{A} \bigr)
\Bigl( i\frac {\partial}{\partial y^0} - \grad \Bigr)
=
- \frac {\partial A_0}{\partial y^0}
+ i\frac {\partial \overrightarrow{A}}{\partial y^0}
- i\grad A_0 + \grad \cdot \overrightarrow{A}
+ \grad \times \overrightarrow{A}.
$$
Therefore,
\begin{multline*}
\nabla_{\BB M} \tilde A \nabla_{\BB M} =
\Bigl( i\frac {\partial}{\partial y^0} - \grad \Bigr)
(\tilde A \nabla_{\BB M}) \\
=
- i\frac {\partial^2 A_0}{(\partial y^0)^2}
- \frac {\partial^2 \overrightarrow{A}}{(\partial y^0)^2}
+ \grad \frac {\partial A_0}{\partial y^0}
+ i\grad \cdot \frac {\partial \overrightarrow{A}}{\partial y^0}
+ i\grad \times \frac {\partial \overrightarrow{A}}{\partial y^0}
+ \grad \frac {\partial A_0}{\partial y^0}
\\
+ i\frac {\partial}{\partial y^0} \bigl( \grad \cdot \overrightarrow{A} \bigr)
- i\frac {\partial}{\partial y^0} \bigl( \grad \times \overrightarrow{A} \bigr)
- i \grad \cdot \bigl( \grad A_0 \bigr) + i \grad \times (\grad A_0)
\\
- \grad \bigl( \grad \cdot \overrightarrow{A} \bigr)
+ \grad \cdot \bigl( \grad \times \overrightarrow{A} \bigr)
- \grad \times \bigl( \grad \times \overrightarrow{A} \bigr)  \\
=
- i\frac {\partial^2 A_0}{(\partial y^0)^2}
- \frac {\partial^2 \overrightarrow{A}}{(\partial y^0)^2}
+ 2 \grad \frac {\partial A_0}{\partial y^0}
+ 2i\grad \cdot \frac {\partial \overrightarrow{A}}{\partial y^0}  \\
- i \grad \cdot \bigl( \grad A_0 \bigr)
- \grad \bigl( \grad \cdot \overrightarrow{A} \bigr)
- \grad \times \bigl( \grad \times \overrightarrow{A} \bigr).
\end{multline*}
We also have
$$
- \square_{3,1} \tilde A^+ = \biggl( \frac{\partial^2}{(\partial y^0)^2} -
\frac{\partial^2}{(\partial y^1)^2} - \frac{\partial^2}{(\partial y^2)^2}
- \frac{\partial^2}{(\partial y^3)^2} \biggr)
\bigl( iA_0 - \overrightarrow{A} \bigr).
$$
We consider separately the scalar and vector terms in 
$\nabla_{\BB M} \tilde A \nabla_{\BB M} - \square_{3,1} \tilde A^+$.
The scalar term is:
\begin{multline*}
- i\frac {\partial^2 A_0}{(\partial y^0)^2}
+ 2i\grad \cdot \frac {\partial \overrightarrow{A}}{\partial y^0}
- i \grad \cdot \grad A_0
+ i \biggl( \frac{\partial^2 A_0}{(\partial y^0)^2}
- \frac{\partial^2 A_0}{(\partial y^1)^2}
- \frac{\partial^2 A_0}{(\partial y^2)^2}
- \frac{\partial^2 A_0}{(\partial y^3)^2} \biggr)  \\
= 2i \biggl( \grad \cdot \frac {\partial \overrightarrow{A}}{\partial y^0}
- \grad \cdot \bigl( \grad A_0 \bigr) \biggr)
= -2i \grad \cdot \biggl(\grad A_0 -
\frac {\partial \overrightarrow{A}}{\partial y^0} \biggr),
\end{multline*}
and the vector term is:
\begin{multline*}
- \frac {\partial^2 \overrightarrow{A}}{(\partial y^0)^2}
+ 2 \grad \frac {\partial A_0}{\partial y^0}
- \grad \bigl( \grad \cdot \overrightarrow{A} \bigr)
- \grad \times \bigl( \grad \times \overrightarrow{A} \bigr)  \\
- \biggl( \frac{\partial^2}{(\partial y^0)^2} -
\frac{\partial^2}{(\partial y^1)^2} - \frac{\partial^2}{(\partial y^2)^2}
- \frac{\partial^2}{(\partial y^3)^2} \biggr) \overrightarrow{A}  \\
= 2 \frac {\partial}{\partial y^0} \biggl( \grad A_0 - 
\frac {\partial \overrightarrow{A}}{\partial y^0} \biggr)
- 2 \grad \times \bigl( \grad \times \overrightarrow{A} \bigr),
\end{multline*}
where we are using the vector calculus identity
$\grad \bigl( \grad \cdot \overrightarrow{A} \bigr)
- \bigl( \grad \cdot \grad \bigr) \overrightarrow{A}
= \grad \times \bigl( \grad \times \overrightarrow{A} \bigr)$.
Thus we conclude that our equation $\M \tilde A=0$ is just a concise
quaternionic form of the Maxwell equations for
the gauge potential (\ref{Maxwell})!

\subsection{On Rings and Modules of Quaternionic Functions}

In Subsection \ref{repro2-section} we obtained a quaternionic version of the
Cauchy formula for double pole using the representation theoretic analogy
with the complex case. Now we would like to extend further the analogy
and introduce a suitable candidate for the ring of quaternionic functions.

Let us recall the representations $(\rho_k, V_k)$, $k=0,1,2$,
of $\mathfrak{sl}(2,\BB C)$ defined by (\ref{rho-intro}).
The multiplication defines an intertwining operator $V_1 \otimes V_1 \to V_2$.
Moreover, one has an algebra structure on $V_0$ and module structures
on $V_1$ and $V_2$ compatible with the action of $\mathfrak{sl}(2,\BB C)$
and the intertwining operator (\ref{derivative-intro}) between $V_0$ and $V_2$.
The latter can be expressed using the Cauchy formula (\ref{2pole-intro}).
Also there is an inverse operator $V_2 \to V_0$ given by the residue
$\frac 1{2\pi i} \oint f(z)\,dz$.

We have seen that these representations
and intertwining operators have exact quaternionic counterparts
${\cal W}'$, ${\cal V}$ (or ${\cal V}'$) and ${\cal W}$,
except there is no algebra structure on ${\cal W}'$.
In order to get this structure one should embed ${\cal W}'$ and ${\cal W}$
into a larger complex that has already appeared in the form of the gauge
transformations of the Maxwell equations
\begin{equation}  \label{W-sequence}
\begin{CD}
{\cal W}'_0 @>{\nabla^+}>> {\cal W}' @>{\M}>>
{\cal W} @>{\re \circ \nabla^+}>> {\cal W}_0,
\end{CD}
\end{equation}
where ${\cal W}'_0$ denotes the scalar-valued functions on $\HC^{\times}$ and
${\cal W}_0$ is the dual space with respect to the pairing given by
(\ref{U2-pairing}) (or (\ref{M-pairing}) in Minkowski formulation).
One can check that these maps are equivariant for the following actions
of the conformal group

\begin{align*}
F(Z) &\mapsto F \bigl( (aZ+b)(cZ+d)^{-1} \bigr),
\qquad \text{on ${\cal W}'_0$}, \\
F(Z) &\mapsto \frac {F \bigl( (aZ+b)(cZ+d)^{-1} \bigr)}
{N(a'-Zc')^2 \cdot N(cZ+d)^2},
\qquad \text{on ${\cal W}_0$},
\end{align*}
$h = \begin{pmatrix} a' & b' \\ c' & d' \end{pmatrix}$,
$h^{-1} = \begin{pmatrix} a & b \\ c & d \end{pmatrix} \in GL(2, \HC)$.

Thus the space ${\cal W}'_0$ of quaternionic scalar functions has a natural
algebra structure, and the other spaces ${\cal W}'$, ${\cal W}$, ${\cal W}_0$
have module structures over ${\cal W}'_0$.
Since the kernel of $\nabla^+ \bigr|_{{\cal W}'_0}$ consists of constants only,
we can transfer the algebra structure to the image of ${\cal W}'_0$ in
${\cal W}'$ extended by constants.
This is the most elementary analogue of the ring of holomorphic functions
in the quaternionic case.

It is also interesting to study the other intertwining operators between
the spaces in (\ref{W-sequence}).
We can consider an integral operator $\xm$:
$$
F(Z) \mapsto (\xm F)(W)=
\int_{U(2)} \frac {(Z-W) \cdot F(Z) \cdot (Z-W)}{N(Z-W)^2} \,dZ^4.
$$
This map is an intertwining operator
$(\rho_2,{\cal W}) \to (\rho'_2, {\cal W}')$,
which by Theorem \ref{2pole1} vanishes on ${\cal W}^+$ and ${\cal W}^-$.
A routine calculation shows that
$$
\M_W \biggl( \frac {(Z-W) \cdot F(Z) \cdot (Z-W)}{N(Z-W)^2} \biggr) =0,
$$
which implies $\M \circ \xm =0$.
It is also true that $\xm \circ \M =0$.

The operator $\xm$ can be embedded into a larger complex
$$
\begin{CD}
{\cal W}_0 @>>> {\cal W} @>{\xm}>> {\cal W}' @>>> {\cal W}'_0,
\end{CD}
$$
There is also an intertwining operator
$$
\square \circ \square : {\cal W}'_0 \to {\cal W}_0
$$
(see \cite{JV1}) and its inverse is given by the quaternionic residue:
$$
{\cal W}_0 \to {\cal W}'_0, \qquad F(Z) \mapsto \int_{U(2)} F(Z) \,dZ^4
$$
One can also consider the intertwining operators
$$
\square \circ \nabla^+ : {\cal W}_0 \to {\cal W}',
\qquad
\re \circ \square \circ \nabla^+ : {\cal W} \to {\cal W}'_0
$$
and their inverses.

\subsection{Bilinear Pairing and Polarization of Vacuum}

In Subsection \ref{repro2-section} we have obtained an intertwining
operator $\M: {\cal W}'^+ \to {\cal W}^+$ and the pairing
$\langle \, \cdot \, , \, \cdot \, \rangle : {\cal W}^+ \times {\cal W}'^-
\to \BB C$.
The combination of these two yields a bilinear map
\begin{equation}  \label{pairing-diagram}
\xymatrix{
{\cal W}'^+ \times {\cal W}'^- \ar[dr]_{\M \otimes 1} & \longrightarrow &
\BB C  \\
& {\cal W}^+ \times {\cal W}'^-
\ar[ur]_{\langle \, \cdot \, , \, \cdot \, \rangle} & }
\end{equation}
which is zero on the kernels of $\M$.
This bilinear map can be presented in a more symmetric way
using the square of the Fueter kernel with values in $\HC \otimes \HC$
$$
\Pi(Z-W) =
\frac {(Z-W)^{-1}}{N(Z-W)} \tilde\otimes \frac {(Z-W)^{-1}}{N(Z-W)}
\cdot \sigma
= \sigma \cdot
\frac {(Z-W)^{-1}}{N(Z-W)} \tilde\otimes \frac {(Z-W)^{-1}}{N(Z-W)},
$$
where $\tilde\otimes$ denotes the Kronecker product of matrices,
$\sigma = \text{\tiny $\begin{pmatrix} 1 & 0 & 0 & 0 \\ 0 & 0 & 1 & 0 \\
0 & 1 & 0 & 0 \\ 0 & 0 & 0 & 1 \end{pmatrix}$}$,
and the double integration
\begin{equation}  \label{W-pairing2}
\frac {-3}{8\pi^6} \int_{W \in U(2)_r} \int_{Z \in U(2)_R}
\tr \Bigl( \bigl( F_1(W) \tilde \otimes 1 \bigr)
\cdot \Pi(Z-W) \cdot \bigl(1 \tilde\otimes F_2(Z) \bigr) \Bigr)\,dZ^4dW^4,
\end{equation}
where $F_1 \in {\cal W}'^+$, $F_2 \in {\cal W}'^-$, $R>r>0$
and the trace is taken in $\HC \otimes \HC$.
Theorem \ref{2pole2} implies that the bilinear maps defined by
(\ref{W-pairing}), (\ref{pairing-diagram}) and (\ref{W-pairing2}) coincide.
However, there is a crucial difference between the two expressions
(\ref{pairing-diagram}) and (\ref{W-pairing2}) when we want to extend them
to bilinear pairings on the full space ${\cal W}'$.
Since both the intertwining operator $\M$ and the pairing 
$\langle \, \cdot \, , \, \cdot \, \rangle$ are still valid when
${\cal W}'^+$ and ${\cal W}^+$ are replaced by ${\cal W}'$ and ${\cal W}$
respectively, the pairing defined by extending (\ref{pairing-diagram})
makes perfect sense. On the other hand, the best obvious extension of 
(\ref{W-pairing2}) works only for
${\cal W}'^+ \oplus {\cal W}'^+ \subsetneq {\cal W}'$ and
${\cal W}^+ \oplus {\cal W}^+ \subsetneq {\cal W}$.
Thus we arrive at a very interesting and important problem
of how to generalize the double integral (\ref{W-pairing2}) so that
it expresses the bilinear pairing defined by (\ref{pairing-diagram})
on the full space ${\cal W}'$.
Clearly, this problem is directly related to the extension to ${\cal W}'$
of the quaternionic Cauchy formula for the second order pole.
Recall that the Cauchy-Fueter formula (Theorem \ref{Fueter}) is valid
for the full spaces ${\cal V}$ and ${\cal V}'$, though in this case we have
a perfect polarization ${\cal V} = {\cal V}^+ \oplus {\cal V}'^-$
and ${\cal V} = {\cal V}'^+ \oplus {\cal V}'^-$.

Certainly, the problem of the integral presentation of the bilinear pairing in
${\cal W}'$ persists in the Minkowski picture, where the integration over
$U(2)$ in (\ref{W-pairing2}) is replaced by the integration over $\BB M$.
But this case is deeply related to a fundamental problem in quantum
electrodynamics known as the vacuum polarization. It does admit a
``physical'' solution which requires a better mathematical understanding.
Physicists manage to redefine the kernel $\Pi(Z-W)$
-- usually called the polarization operator -- by considering its
Fourier transform and subtracting unwanted infinite terms
(for example, see \cite{BS}). As a result, they obtain
\begin{equation*}
\operatorname{reg} \Pi(Z-W)=
c \Bigl( \nabla \delta^{(4)}(Z-W)\nabla
- \square \bigl( \delta^{(4)}(Z-W) \bigr)^+ \Bigr),
\end{equation*}
where the constant $c$ depends on the regularization procedure, and
the differential operators are applied to the $\HC \otimes \HC$-valued delta
function in such a way that the integration over $Z \in \BB M$ or $W \in \BB M$
yields the pairing (\ref{pairing-diagram}) defined by the second integration
over $\BB M$ and the Maxwell operator $\M$.
The polarization operator has a standard graphical presentation by means
of the divergent Feynman diagram of the second order

\vskip7pt

\centerline{\includegraphics{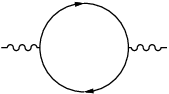}}

\noindent
We will encounter more Feynman diagrams in the next section,
which indicates that their appearance in quaternionic analysis
is not a mere accident but a reflection of a profound relation.

\section{Middle Series and Decomposition of Tensor Product}
\label{decomposition-section}

\subsection{Irreducible Components of ${\cal H}^+ \otimes {\cal H}^+$}

In the previous section we studied the top irreducible component $\Zh^+$
of the tensor product ${\cal H}^+ \otimes {\cal H}^+$.
In this subsection we describe all the irreducible components in the
decomposition of this tensor product, which consist of all middle
series representations.

For $n \in \BB N$, let $\BB C^{n \times n}$ denote the space of complex
$n \times n$ matrices.
The group $GL(2,\BB H)$ acts on $\BB C^{n \times n}$-valued functions on
$\HC$ with isolated singularities via
\begin{multline}  \label{n-action}
\rho_n(h): \: F(Z) \mapsto \bigl( \rho_n(h)F \bigr)(Z) =
\frac {\tau_n(cZ+d)^{-1}}{N(cZ+d)} \cdot
F \bigl( (aZ+b)(cZ+d)^{-1} \bigr) \cdot
\frac {\tau_n(a'-Zc')^{-1}}{N(a'-Zc')},  \\
h = \begin{pmatrix} a' & b' \\ c' & d' \end{pmatrix},\:
h^{-1} = \begin{pmatrix} a & b \\ c & d \end{pmatrix} \in GL(2, \BB H),
\end{multline}
where $cZ+d$ and $a'-Zc'$ are regarded as elements of $GL(2,\BB C)$ and
$\tau_n: GL(2,\BB C) \to\operatorname{Aut}(\BB C^n) \subset \BB C^{n \times n}$
is the irreducible $n$-dimensional representation of $GL(2,\BB C)$.
For $n=1$, $\tau_1 \equiv 1$ and (\ref{n-action}) agrees with (\ref{1-action});
and, for $n=2$, (\ref{n-action}) agrees with (\ref{action}). Let
$$
\Zh_n = \Zh \otimes \BB C^{n \times n}, \qquad
\Zh^+_n = \Zh^+ \otimes \BB C^{n \times n}, \qquad
\Zh^-_n = \Zh^- \otimes \BB C^{n \times n}
$$
be the spaces of polynomial functions with values in $\BB C^{n \times n}$
defined on $\HC^{\times}$, $\HC$ and $\HC^{\times} \cup \{\infty\}$
respectively; $\Zh^+_n \oplus \Zh^-_n$ is a proper subspace of $\Zh_n$.
Differentiating $\rho_n$-action and complexifying,
we obtain a representation of $\mathfrak{gl}(4, \BB C)$
in $\Zh_n$ which preserves $\Zh^+_n$ and $\Zh^-_n$.
Note that $(\rho_n, \Zh_n^+) \simeq (\rho_n, \Zh_n^-)$,
$(\rho_1, \Zh_1) = (\rho_1, \Zh)$ and
$(\rho_2, \Zh_2) = (\rho_2, {\cal W})$.

\begin{thm} [\cite{JV2}]  \label{JV2-thm}
The representations $(\rho_n, \Zh_n^+)$, $n=0,1,2,3,\dots$,
of $\mathfrak{gl}(4, \BB C)$ are irreducible.
They possess inner products which make them unitary
representations of the real form $\mathfrak{u}(2,2)$
of $\mathfrak{gl}(4, \BB C)$.
\end{thm}

As in Subsection \ref{repro2-section}, we can also consider the actions
\begin{multline*}  
\rho'_n(h): \: F(Z) \mapsto \bigl( \rho_n(h)F \bigr)(Z) =
\frac {\tau_n(a'-Zc')}{N(a'-Zc')} \cdot
F \bigl( (aZ+b)(cZ+d)^{-1} \bigr) \cdot \frac {\tau_n(cZ+d)}{N(cZ+d)},  \\
h = \begin{pmatrix} a' & b' \\ c' & d' \end{pmatrix},\:
h^{-1} = \begin{pmatrix} a & b \\ c & d \end{pmatrix} \in GL(2, \BB H).
\end{multline*}
Differentiating these actions and complexifying,
we obtain representations of $\mathfrak{gl}(4, \BB C)$
in $\Zh_n$ which preserves $\Zh^+_n$ and $\Zh^-_n$.
Denote by $\Zh'_n$, $\Zh'^+_n$ and $\Zh'^-_n$ the spaces
$\Zh_n$, $\Zh^+_n$ and $\Zh^-_n$ with $\rho'_n$-action.
Note that $(\Zh'_1, \rho'_1) = (\Zh, \rho_1)$ and
$(\Zh'_2, \rho'_2) = ({\cal W}', \rho'_2)$.
Then $(\Zh^{\pm}_n, \rho_n)$ can be realized as irreducible quotients of
$(\Zh'^{\pm}_n, \rho'_n)$.
The same proof as in Proposition \ref{Mx-intertwiner} shows:
\begin{prop}  \label{Mx-intertwiner-n}
We have an intertwining operator
$$
\M_n:\: (\Zh'^{\pm}_n, \rho'_n) \to (\Zh^{\pm}_n, \rho_n),
$$
where $\M_n$ is defined by
$$
\M_n F(W) = \frac {3i}{\pi^3} \int_{Z \in U(2)_R}
\frac {\tau_n(Z-W)^{-1}}{N(Z-W)} \cdot F(Z) \cdot
\frac {\tau_n(Z-W)^{-1}}{N(Z-W)} \,dZ^4,
\qquad F \in \Zh_n.
$$
(This definition is independent of $R>0$.)
\end{prop}

Similarly, the proof of Proposition \ref{W-W'_pairing} shows:
\begin{prop}
There is a $\mathfrak{u}(2,2)$-invariant bilinear pairing between the
representations $(\rho_n, \Zh_n)$ and $(\rho'_n, \Zh'_n)$ given by
$$
\langle F_1,F_2 \rangle =
\frac i{8\pi^3} \int_{U(2)_R} \tr \bigl( F_1(Z) \cdot F_2(Z) \bigr) \,dZ^4,
\qquad F_1 \in \Zh_n, \: F_2 \in \Zh'_n.
$$
This pairing is independent of the choice of $R>0$.
Moreover, for each $F_1 \in \Zh^+_n$, $F_1 \ne 0$,
there is an $F_2 \in \Zh'^-_n$ such that $\langle F_1,F_2 \rangle \ne 0$.
\end{prop}

Since $(\rho_n,\Zh^{\pm}_n)$ is irreducible, Proposition \ref{Mx-intertwiner-n}
implies that $\Zh^{\pm}_n \simeq \Zh'^{\pm}_n/\ker \M_n$.
Thus we obtain a $\mathfrak{u}(2,2)$-invariant bilinear pairing on
$\Zh'^+_n \times \Zh'^-_n$:
\begin{multline*}
\langle F_1,F_2 \rangle_n =  \\
\frac {-3}{8\pi^6} \int_{W \in U(2)_r} \int_{Z \in U(2)_R}
\tr \biggl( F_1(W) \cdot \frac{\tau_n(Z-W)^{-1}}{N(Z-W)} \cdot F_2(Z) \cdot
\frac{\tau_n(Z-W)^{-1}}{N(Z-W)} \biggr) \,dZ^4 \wedge dW^4.
\end{multline*}
where $F_1 \in \Zh'^+_n$, $F_2 \in \Zh'^-_n$ and $R>r>0$.
This pairing is independent of the choices of $R>r>0$,
vanishes on the kernels of $\M_n$ and descends to a non-degenerate
pairing on
$(\Zh'^+_n / \ker \M_n) \times (\Zh'^-_n / \ker \M_n) \simeq
\Zh^+_n \times \Zh^-_n$.
For $n=1$, $2$ we get (\ref{U2-pairing}), (\ref{W-pairing}) respectively.

According to \cite{JV2}, the representation
$(\pi^0_{lr}, {\cal H}^+ \otimes {\cal H}^+)$ of $\mathfrak{gl}(2,\BB H)$
decomposes into a direct sum of irreducible representations
\begin{equation}  \label{tensor-decomp}
{\cal H}^+ \otimes {\cal H}^+ = \bigoplus_{n=1}^{\infty} \Zh_n^+.
\end{equation}
We describe the intertwining maps
$M_n : {\cal H}^+ \otimes {\cal H}^+ \to \Zh_n^+$.

\begin{thm}
For each $n \ge 1$, the intertwining map
$M_n : {\cal H}^+ \otimes {\cal H}^+ \to \Zh_n^+$ is given by
\begin{equation}  \label{M_n}
\phi_1(Z_1) \otimes \phi_2(Z_2) \mapsto
\Bigl( \tau_n(\partial_{Z_1} - \partial_{Z_2})
\bigl( \phi_1(Z_1) \cdot \phi_2(Z_2) \bigr) \Bigr) \Bigr|_{Z_1=Z_2=Z},
\qquad \phi_1, \phi_2 \in {\cal H}^+.
\end{equation}
\end{thm}

\begin{rem}
Here $\partial = \partial_Z$ is a first order differential operator
$\begin{pmatrix} \partial_{11} & \partial_{21} \\
\partial_{12} & \partial_{22} \end{pmatrix}:
\Zh_1 \to \Zh_2$,
and $\tau_n$ is regarded as a map
$\mathfrak{gl}(2,\BB C) \simeq \BB C^{2\times 2} \to \BB C^{n \times n}$
so that $\tau_n(\partial)$ is an $n \times n$ matrix of differential
operators of order $n-1$.
If $n=1$, $\tau_1$ is defined to be 1, i.e.
$$
M_1: \:
\phi_1(Z_1) \otimes \phi_2(Z_2) \mapsto \phi_1(Z) \cdot \phi_2(Z),
\qquad \phi_1, \phi_2 \in {\cal H}^+,
$$
which is equation (\ref{M_1}).
\end{rem}

\pf
Let $\tilde M_n$ denote the map ${\cal H}^+ \otimes {\cal H}^+ \to \Zh_n^+$
given by (\ref{M_n}). Clearly, $\tilde M_n$ is not trivial.
We need to show that $\tilde M_n$ is $GL(2,\BB H)$-equivariant.
The group $GL(2,\BB H)$ is generated by the translation operators
$\begin{pmatrix} 1 & b \\ 0 & 1 \end{pmatrix}$, diagonal matrices
$\begin{pmatrix} a & 0 \\ 0 & d \end{pmatrix}$, and the inversion
$\begin{pmatrix} 0 & 1 \\ 1 & 0 \end{pmatrix}$, so it is sufficient
to verify that $\tilde M_n$ respects the actions of these elements only.
Clearly, $\tilde M_n$ respects the actions of the translation operators.
A simple computation shows that, for all $F \in \Zh^+$ and
$a,d \in \BB H^{\times}$,
$$
\partial \bigl( F(aZ)) = (\partial F) \bigr|_{aZ} \cdot a, \qquad
\partial \bigl( F(Zd)) = d \cdot (\partial F) \bigr|_{Zd}.
$$
This proves that $\tilde M_2$ respects the actions of the diagonal matrices.
Moreover, these equations immediately imply
$$
\tau_n(\partial) \bigl( F(aZ)) =
\bigl( \tau_n(\partial) F \bigr) \bigr|_{aZ} \cdot \tau_n(a), \qquad
\tau_n(\partial) \bigl( F(Zd)) =
\tau_n(d) \cdot \bigl( \tau_n(\partial) F \bigr) \bigr|_{Zd}.
$$
and hence $\tilde M_n$ respects the actions of the diagonal matrices.

\begin{lem}
For $F \in \Zh^+$, we have:
$$
\partial \bigl( F(Z^{-1}) \bigr) =
- Z^{-1} \cdot (\partial F)(Z^{-1}) \cdot Z^{-1}.
$$
\end{lem}

\pf
We write $Z^{-1} = \frac {Z^+}{N(Z)}$ and also note that
$\partial \bigl( F(Z^+) \bigr) = (\partial^+ F)(Z^+)$.
We have:
\begin{multline*}
\partial \bigl( F(Z^{-1}) \bigr) = \partial \bigl( F(Z^+ / N(Z)) \bigr)  \\
= \frac 1{N(Z)}  \cdot (\partial^+ F)(Z^{-1})
- \frac {Z^+}{N(Z)^2} \cdot \bigl( (z_{22}\partial_{11} - z_{12}\partial_{12}
- z_{21}\partial_{21} + z_{11}\partial_{22}) F \bigr) (Z^{-1})  \\
= \frac {Z^+}{N(Z)^2}  \cdot \biggl(
\begin{pmatrix}
-z_{22} \partial_{11} + z_{21}\partial_{21} &
-z_{11} \partial_{21} + z_{12}\partial_{11}  \\
z_{21} \partial_{22} - z_{22}\partial_{12} &
z_{12} \partial_{12} - z_{11}\partial_{22}
\end{pmatrix}
F \biggr) (Z^{-1}) \\
= - \frac {Z^+}{N(Z)^2}  \cdot (\partial F)(Z^{-1}) \cdot Z^+
= - Z^{-1} \cdot (\partial F)(Z^{-1}) \cdot Z^{-1}.
\end{multline*}
\qed

For $n \ge 2$ we have:
\begin{multline*}
\tau_n(\partial_{Z_1} - \partial_{Z_2}) \biggl(
\frac{F(Z_1^{-1}) \cdot G(Z_2^{-1})}{N(Z_1)\cdot N(Z_2)} \biggr)  \\
=
\frac {\tau_n(Z_1^{-1}) \cdot \bigl(\tau_n(\partial_{Z_1}) F \bigr)(Z_1^{-1})
\cdot G(Z_2^{-1}) \cdot \tau_n(-Z_1^{-1})}{N(Z_1) \cdot N(Z_2)}  \\
-
\frac {\tau_n(Z_2^{-1}) \cdot F (Z_1^{-1}) \cdot
\bigl(\tau_n(\partial_{Z_2}) G \bigr)(Z_2^{-1})
\cdot \tau_n(-Z_2^{-1})}{N(Z_1) \cdot N(Z_2)}
\end{multline*}
plus terms involving
$$
\tau_{n-k}(\partial_{Z_1} - \partial_{Z_2})
\bigl( F(Z_1^{-1}) \cdot G(Z_2^{-1}) \bigr)
\otimes \tau_k(\partial_{Z_1} - \partial_{Z_2})
\Bigl( \frac 1{N(Z_1)N(Z_2)} \Bigr),
\qquad k>0.
$$
But
$$
\tau_k(\partial_{Z_1} - \partial_{Z_2}) \Bigl( \frac 1{N(Z_1)N(Z_2)} \Bigr)
\Bigr|_{Z_1=Z_2=Z} =0,
$$
so the remaining terms cancel out.
This proves that $\tilde M_n$ respects the inversion.
\qed

\subsection{Integral Kernels of Projections onto Irreducible Components}
\label{projectors1}

Now let us consider the projectors
$$
{\cal P}^0_n :\:
{\cal H}^+ \otimes {\cal H}^+ \twoheadrightarrow \Zh^+_n
\hookrightarrow {\cal H}^+ \otimes {\cal H}^+
$$
and
$$
{\cal P} :\:
{\cal V}^+ \otimes {\cal V}'^+ \twoheadrightarrow {\cal W}^+
\hookrightarrow {\cal V}^+ \otimes {\cal V}'^+.
$$

First we observe that, using the Poisson formula (Theorem \ref{Poisson}),
the map $M_1 : {\cal H}^+ \otimes {\cal H}^+ \twoheadrightarrow \Zh^+_1$
can be expressed by the following integral formula:
\begin{multline*}
M_1 \bigl( \phi_1(Z_1) \otimes \phi_2(Z_2) \bigr) (T)  \\
= \frac 1{(2\pi^2)^2} \int_{Z_1 \in S^3_R} \int_{Z_2 \in S^3_R}
\frac 1{N(Z_1-T) \cdot N(Z_2-T)} \cdot
(\widetilde{\deg}_{Z_1} \phi_1)(Z_1) \cdot (\widetilde{\deg}_{Z_2} \phi_2)(Z_2)
\, \frac{dS dS}{R^2}
\end{multline*}
defined for $|T|<R$.
Similarly, the map
$M :\: {\cal V}^+ \otimes {\cal V}'^+ \twoheadrightarrow {\cal W}^+$
can be expressed by the integral formula:
$$
M \bigl( f(Z_1) \otimes g(Z_2) \bigr) (T) =
\frac 1{(2\pi^2)^2} \int_{S^3_R} \int_{S^3_R}
k(Z_1-T) \cdot Dz_1 \cdot f(Z_1) \cdot g(Z_2) \cdot Dz_2 \cdot k(Z_2-T)
$$
also defined for $|T|<R$.
In general, the map
$M_n : {\cal H}^+ \otimes {\cal H}^+ \twoheadrightarrow \Zh^+_n$
can be written as
\begin{multline*}
M_n \bigl( \phi_1(Z_1) \otimes \phi_2(Z_2) \bigr) (T)  \\
= \frac 1{(2\pi^2)^2} \int_{Z_1 \in S^3_R} \int_{Z_2 \in S^3_R}
m_n(Z_1,Z_2,T) \cdot
(\widetilde{\deg}_{Z_1} \phi_1)(Z_1) \cdot (\widetilde{\deg}_{Z_2} \phi_2)(Z_2)
\, \frac{dS dS}{R^2},
\end{multline*}
where $|T|<R$ and
$$
m_n(Z_1,Z_2,T) =
\Bigl( \tau_n(\partial_{T_1} - \partial_{T_2})
\frac 1{N(Z_1-T_1) \cdot N(Z_2-T_2)} \Bigr) \Bigr|_{T_1=T_2=T}.
$$

\begin{lem}  \label{m_n}
For $n \ge 1$ we have:
$$
m_n(Z_1,Z_2,T) = (n-1)! \cdot
\frac {\tau_n \bigl( (T-Z_1)^{-1}(Z_1-Z_2)(Z_2-T)^{-1} \bigr)}
{N(Z_1-T) \cdot N(Z_2-T)}.
$$
\end{lem}

\pf
First we check that, for $n \ge 1$,
$$
\tau_n (\partial_{T_i}) \frac 1{N(Z_i-T_i)} =
\frac{(n-1)!}{N(Z_i-T_i)} \cdot \tau_n (Z_i-T_i)^{-1},
\qquad i=1,2.
$$
Then it follows that
\begin{multline*}
\tau_n(\partial_{T_1} - \partial_{T_2}) \frac 1{N(Z_1-T_1) \cdot N(Z_2-T_2)}
\\
= \frac {(n-1)!}{N(Z_1-T_1) \cdot N(Z_2-T_2)} \cdot \tau_n \bigl(
(Z_1-T_1)^{-1} - (Z_2-T_2)^{-1} \bigr)  \\
= (n-1)! \cdot \frac {\tau_n \bigl(
(Z_1-T_1)^{-1}(Z_2-Z_1-T_2+T_1)(Z_2-T_2)^{-1} \bigr)}
{N(Z_1-T_1) \cdot N(Z_2-T_2)}.
\end{multline*}
\qed

The main result of this section is

\begin{thm}
For each $n \ge 1$, the projector ${\cal P}^0_n$ has integral kernel
\begin{equation}  \label{p_n}
p^0_n(Z_1,Z_2;W_1,W_2) =
\bigl\langle m_n(Z_1,Z_2,T), m_n(W_1,W_2,T) \bigr\rangle_n,
\end{equation}
where the pairing is done with respect to the variable $T$.
Namely, there exists a $\lambda_n \in \BB C$ such that,
for $R>0$ and $|W_1|,|W_2| <R$, we have
\begin{multline} \label{P^0_n-tilde}
({\cal P}^0_n \phi_1 \otimes \phi_2)(W_1,W_2)  \\
= \frac {\lambda_n}{(2\pi^2)^2}
\int_{Z_1 \in S^3_R} \int_{Z_2 \in S^3_R} p^0_n(Z_1,Z_2;W_1,W_2) \cdot
(\widetilde{\deg}_{Z_1} \phi_1)(Z_1) \cdot (\widetilde{\deg}_{Z_2} \phi_2)(Z_2)
\, \frac{dS dS}{R^2}.
\end{multline}
\end{thm}

\pf
Just as we have (\ref{tensor-decomp}),
${\cal H}^- \otimes {\cal H}^-$ decomposes as $\bigoplus_{n=1}^{\infty} \Zh_n^-$.
Let $M_n$ denote the projection
${\cal H}^- \otimes {\cal H}^- \twoheadrightarrow \Zh^-_n$.
Since $\Zh^{\pm}_n$ is irreducible, there exists a $\lambda_n \in \BB C$
such that
\begin{multline*}
\bigl\langle {\cal P}^0_n (\phi_1 \otimes \phi_2),
\phi_1' \otimes \phi_2' \bigr\rangle
= \lambda_n \cdot \bigl\langle M_n (\phi_1 \otimes \phi_2),
M_n (\phi_1' \otimes \phi_2') \bigr\rangle_n,  \\
\forall \phi_1,\phi_2 \in {\cal H}^+, \:
\forall \phi_1',\phi_2' \in {\cal H}^-,
\end{multline*}
where the first pairing is taken inside
$({\cal H}^+ \otimes {\cal H}^+) \times ({\cal H}^- \otimes {\cal H}^-)$
using (\ref{H-pairing}),
and the second -- inside $\Zh_n^+ \times \Zh_n^-$.

On the other hand, let $\widetilde{\cal P}^0_n:
{\cal H}^+ \otimes {\cal H}^+ \to {\cal H}^+ \otimes {\cal H}^+$ be the map
given by the integral operator (\ref{P^0_n-tilde}).
We want to show that $\widetilde{\cal P}^0_n = {\cal P}^0_n$.
It is sufficient to show that
\begin{multline*}
\bigl\langle \widetilde{\cal P}^0_n (\phi_1 \otimes \phi_2),
\phi_1' \otimes \phi_2' \bigr\rangle
= \lambda_n \cdot \bigl\langle M_n (\phi_1 \otimes \phi_2),
M_n (\phi_1' \otimes \phi_2') \bigr\rangle_n,  \\
\forall \phi_1,\phi_2 \in {\cal H}^+, \:
\forall \phi_1',\phi_2' \in {\cal H}^-.
\end{multline*}
Indeed, by Corollary \ref{Poisson 1/X} and the argument preceding
Lemma \ref{m_n}, for $\phi_1',\phi_2' \in {\cal H}^-$
and $0 < r < |W_1|, |W_2|$,
\begin{multline*}
M_n \bigl( \phi_1'(W_1) \otimes \phi_2'(W_2) \bigr) (T)  \\
= \frac 1{(2\pi^2)^2} \int_{W_1 \in S^3_r} \int_{W_2 \in S^3_r}
m_n(W_1,W_2,T) \cdot (\widetilde{\deg}_{W_1} \phi'_1)(W_1)
\cdot (\widetilde{\deg}_{W_2} \phi'_2)(W_2) \, \frac{dS dS}{r^2}.
\end{multline*}
Hence
\begin{multline*}
\lambda_n \cdot \bigl\langle M_n (\phi_1 \otimes \phi_2),
M_n (\phi_1' \otimes \phi_2') \bigr\rangle_n  \\
= \frac {\lambda_n}{(2\pi^2)^4} \cdot \biggl\langle
\int_{Z_1 \in S^3_R} \int_{Z_2 \in S^3_R} m_n(Z_1,Z_2,T) \cdot
(\widetilde{\deg}_{Z_1} \phi_1)(Z_1) \cdot (\widetilde{\deg}_{Z_2} \phi_2)(Z_2)
\,\frac{dS dS}{R^2},  \\
\int_{W_1 \in S^3_r} \int_{W_2 \in S^3_r} m_n(W_1,W_2,T) \cdot
(\widetilde{\deg}_{W_1} \phi'_1)(W_1) \cdot
(\widetilde{\deg}_{W_2} \phi'_2)(W_2) \,\frac{dS dS}{r^2}
\biggr\rangle_n  \\
= \frac {\lambda_n}{(2\pi^2)^4}
\int_{Z_1 \in S^3_R} \int_{Z_2 \in S^3_R} \int_{W_1 \in S^3_r}
\int_{W_2 \in S^3_r}
\bigl\langle m_n(Z_1,Z_2,T), m_n(W_1,W_2,T) \bigr\rangle_n   \\
\cdot (\widetilde{\deg}_{Z_1} \phi_1)(Z_1) \cdot
(\widetilde{\deg}_{Z_2} \phi_2)(Z_2) \cdot 
(\widetilde{\deg}_{W_1} \phi'_1)(W_1) \cdot
(\widetilde{\deg}_{W_2} \phi'_2)(W_2)
\,\frac{dS dS}{r^2}\frac{dS dS}{R^2}  \\
= \frac 1{(2\pi^2)^2} \int_{W_1 \in S^3_r} \int_{W_2 \in S^3_r}
\widetilde{\cal P}^0_n (\phi_1 \otimes \phi_2) (W_1,W_2) \cdot
(\widetilde{\deg}_{W_1} \phi_1')(W_1) \cdot
(\widetilde{\deg}_{W_1} \phi_2')(W_2) \,\frac{dS dS}{r^2}  \\
= \bigl\langle \widetilde{\cal P}^0_n (\phi_1 \otimes \phi_2),
\phi_1' \otimes \phi_2' \bigr\rangle.
\end{multline*}
\qed

\begin{cor}
This theorem combined with Proposition \ref{pairing_n=1}
shows that for $n=1$ the kernel
\begin{multline*}
p^0_1(Z_1,Z_2;W_1,W_2) =
\bigl\langle m_1(Z_1,Z_2,T), m_1(W_1,W_2,T) \bigr\rangle_1  \\
= \frac i{8\pi^3} \int_{U(2)}
\frac {1}{N(Z_1-T) \cdot N(Z_2-T) \cdot N(W_1-T) \cdot N(W_2-T)} \,dT^4.
\end{multline*}
\end{cor}

The last integral can be computed as follows.
For $k=1,2$, write
$$
T'=(T-i)(T+i)^{-1}, \qquad
Z'_k=(Z_k-i)(Z_k+i)^{-1}, \qquad
W'_k=(W_k-i)(W_k+i)^{-1},
$$
$T' \in \BB M$, $Z'_1,Z'_2, W'_1,W'_2 \in \BB T^+$.
Then, using Lemmas \ref{X-Y_difference} and \ref{Jacobian_lemma},
we can rewrite
$$
p^0_1(Z_1,Z_2;W_1,W_2) =
N(Z'_1+i) \cdot N(Z'_2+i) \cdot N(W'_1+i) \cdot N(W'_2+i) \cdot
p'^0_1(Z_1,Z_2;W_1,W_2),
$$
where
$$
p'^0_1(Z_1,Z_2;W_1,W_2) = \frac i{8\pi^3} \int_{\BB M}
\frac {1}{N(Z_1-T) \cdot N(Z_2-T) \cdot N(W_1-T) \cdot N(W_2-T)} \,dT^4.
$$
This is a one-loop Feynman integral, it has been computed in terms of the
dilogarithm function and it plays a significant role in physics.
We will discuss Feynman integrals in Subsection \ref{Feynman}.

\subsection{Action of the Casimir Element}

As will be explained in Subsection \ref{Feynman}, it is very important to
find explicit expressions for the integral kernels $p^0_n(Z_1,Z_2;W_1,W_2)$'s.
For this purpose we propose to study the action of the Casimir element
$\Omega \in {\cal U}(\mathfrak{gl}(4,\BB C))$.
Since the representations $(\rho_n, \Zh^+_n)$ are irreducible,
$\Omega$ acts on them by scalars $\mu_n$.
Moreover, for different $n$ the scalars $\mu_n$ are different.
Thus, in order to prove that a given operator
$\widetilde{\cal P} : {\cal H}^+ \otimes {\cal H}^+ \to
{\cal H}^+ \otimes {\cal H}^+$
has image in $\Zh^+_n$, it is sufficient to show that
$\pi_{lr}^0(\Omega) \circ \widetilde{\cal P} = \mu_n \cdot \widetilde{\cal P}$.

\begin{lem}  \label{dz-relations}
\begin{enumerate}
\item
For any $F(Z) \in \Zh_1$, we have $\partial(ZF) = (\partial F) \cdot Z + 2F$.
\item
For any $G(Z) \in \Zh_2$, we have
$\tr \bigl( \partial(GZ) \bigr) = \tr (Z \cdot \partial G) + 2\tr(G)$.
\end{enumerate}
\end{lem}

We define the Casimir element $\Omega$ of $\mathfrak{gl}(4,\BB C)$ relative to
the invariant symmetric bilinear form $(Z_1,Z_2) = \tr(Z_1Z_2)$.

\begin{prop}
The Casimir element $\Omega$ of $\mathfrak{gl}(4,\BB C)$ acts on
$\Zh_1$ by $F \mapsto -4F$.
\end{prop}

\pf
From Lemma \ref{rho_1-algebra-action}, $\rho_1(\Omega)$ acts on
$F(Z) \in \Zh_1$ by
$$
\rho_1(\Omega)F
= \tr \bigl ((-Z \cdot \partial -1)^2 F + (\partial Z -1)^2 F
- \partial \circ (Z \cdot \partial Z) F
- (Z \cdot \partial Z) \circ \partial F \bigr).
$$
This expression can be easily simplified using Lemma \ref{dz-relations} to
$\rho_1(\Omega)F = -4F$.
\qed

\begin{lem}  \label{pi_lr-algebra-action}
The Lie algebra action of $\mathfrak{gl}(2,\BB H)$ on its representation
$(\pi_{lr}^0, {\cal H}^+ \otimes {\cal H}^+)$ is given by
\begin{align*}
\pi_{lr}^0 \begin{pmatrix} A & 0 \\ 0 & 0 \end{pmatrix} &:
F \mapsto \tr \bigl( A \cdot
(-Z_1 \cdot \partial_{Z_1}F - Z_2 \cdot \partial_{Z_2}F -F) \bigr)
\\
\pi_{lr}^0 \begin{pmatrix} 0 & B \\ 0 & 0 \end{pmatrix} &:
F \mapsto \tr \bigl( B \cdot (-\partial_{Z_1}F - \partial_{Z_2}F) \bigr)
\\
\pi_{lr}^0 \begin{pmatrix} 0 & 0 \\ C & 0 \end{pmatrix} &:
F \mapsto \tr \Bigl( C \cdot \bigl(
Z_1 \cdot (\partial_{Z_1}F) \cdot Z_1 + Z_2 \cdot (\partial_{Z_2}F) \cdot Z_2
+ Z_1 F + Z_2F \bigr) \Bigr)
\\
&: F \mapsto \tr \Bigl( C \cdot \bigl(
Z_1 \cdot \partial_{Z_1} (Z_1F) + Z_2 \cdot \partial_{Z_2} (Z_2F)
- Z_1 F - Z_2F \bigr) \Bigr)
\\
\pi_{lr}^0 \begin{pmatrix} 0 & 0 \\ 0 & D \end{pmatrix} &:
F \mapsto \tr \Bigl( D \cdot \bigl(
(\partial_{Z_1}F) \cdot Z_1 + (\partial_{Z_2}F) \cdot Z_2 +F \bigr) \Bigr)
\\
&: F \mapsto \tr \Bigl( D \cdot \bigl(
\partial_{Z_1} (Z_1F) + \partial_{Z_2} (Z_2F) -3F \bigr) \Bigr)
\end{align*}
\end{lem}

The actions of $\begin{pmatrix} A & 0 \\ 0 & 0 \end{pmatrix}$ and
$\begin{pmatrix} 0 & B \\ 0 & 0 \end{pmatrix}$ are obtained
by direct computation and the actions of the other elements are obtained by
writing
$$
\begin{pmatrix} 0 & 0 \\ C & 0 \end{pmatrix} =
\begin{pmatrix} 0 & 1 \\ 1 & 0 \end{pmatrix} \cdot
\begin{pmatrix} 0 & C \\ 0 & 0 \end{pmatrix} \cdot
\begin{pmatrix} 0 & 1 \\ 1 & 0 \end{pmatrix}
\qquad \text{and} \qquad
\begin{pmatrix} 0 & 0 \\ 0 & D \end{pmatrix} =
\begin{pmatrix} 0 & 1 \\ 1 & 0 \end{pmatrix} \cdot
\begin{pmatrix} D & 0 \\ 0 & 0 \end{pmatrix} \cdot
\begin{pmatrix} 0 & 1 \\ 1 & 0 \end{pmatrix}.
$$

\begin{prop}
The Casimir element $\Omega$ of $\mathfrak{gl}(4,\BB C)$ acts on
${\cal H} \otimes {\cal H}$ by
$$
F \mapsto -4F -
\tr \Bigl( (Z_1-Z_2) \partial_{Z_1} \bigl( (Z_1-Z_2) \partial_{Z_2} F
\bigr) \Bigr)
-
\tr \Bigl( (Z_1-Z_2) \partial_{Z_2} \bigl( (Z_1-Z_2) \partial_{Z_1} F
\bigr) \Bigr).
$$
\end{prop}

\pf
From Lemma \ref{pi_lr-algebra-action}, $\pi_{lr}^0(\Omega)$ acts on
$F(Z_1,Z_2) \in {\cal H} \otimes {\cal H}$
by sending it into the trace of
\begin{multline*}
(-Z_1 \cdot \partial_{Z_1} - Z_2 \cdot \partial_{Z_2} -1)^2 F
+ \bigl( \partial_{Z_1} Z_1 + \partial_{Z_2} Z_2 -3 \bigr)^2 F  \\
+ (-\partial_{Z_1} - \partial_{Z_2}) \circ 
\bigl(
Z_1 \cdot \partial_{Z_1} Z_1 + Z_2 \cdot \partial_{Z_2} Z_2
- Z_1 - Z_2 \bigr) F  \\
+ \bigl( Z_1 \cdot \partial_{Z_1} Z_1 + Z_2 \cdot \partial_{Z_2} Z_2
- Z_1 - Z_2 \bigr)
\circ (-\partial_{Z_1} - \partial_{Z_2}) F.
\end{multline*}
This expression can be simplified using Lemma \ref{dz-relations},
and the result follows.
\qed

\begin{lem}
For $T \in \HC$ and $F(Z_1,Z_2) = \frac 1{N(Z_1-T) \cdot N(Z_2-T)}$,
$$
(Z_1-Z_2) \partial_{Z_1} \bigl( (Z_1-Z_2) \partial_{Z_2} F \bigr)
+
(Z_1-Z_2) \partial_{Z_2} \bigl( (Z_1-Z_2) \partial_{Z_1} F \bigr) =0.
$$
In particular, $\pi_{lr}^0(\Omega) F = -4F$.
\end{lem}

\pf
An easy calculation shows that
\begin{multline*}
(Z_1-Z_2) \partial_{Z_1} \bigl( (Z_1-Z_2) \partial_{Z_2} F \bigr)  \\
=
\frac {Z_1-Z_2}{N(Z_1-T) \cdot N(Z_2-T)} \cdot
\bigl( (Z_1-T)^{-1}(Z_1-Z_2)(Z_2-T)^{-1} - 2(Z_2-T)^{-1} \bigr)
\end{multline*}
and
\begin{multline*}
(Z_1-Z_2) \partial_{Z_2} \bigl( (Z_1-Z_2) \partial_{Z_1} F \bigr)  \\
=
\frac {Z_1-Z_2}{N(Z_1-T) \cdot N(Z_2-T)} \cdot
\bigl( (Z_2-T)^{-1}(Z_1-Z_2)(Z_1-T)^{-1} + 2(Z_1-T)^{-1} \bigr).
\end{multline*}
Then the  result follows from the identity
\begin{multline*}
(Z_1-T)^{-1} - (Z_2-T)^{-1} = - (Z_1-T)^{-1}(Z_1-Z_2)(Z_2-T)^{-1}  \\
= - (Z_2-T)^{-1}(Z_1-Z_2)(Z_1-T)^{-1}.
\end{multline*}
\qed

\subsection{Kernels of Projectors in the Minkowski Case and Feynman Integrals}
\label{Feynman}

In Subsection \ref{projectors1} we have derived an explicit expression for
the kernel of the projector ${\cal P}'^0_1$ as a certain integral over $U(2)$.
The equivalent problem in the Minkowski case is to find the projectors
$$
{\cal P}'^0_n :
{\cal H}(\BB M)^+ \otimes {\cal H}(\BB M)^+ \twoheadrightarrow \Zh_n(\BB M)^+
\hookrightarrow {\cal H}(\BB M)^+ \otimes {\cal H}(\BB M)^+.
$$
By the same reasoning we get

\begin{thm}
The projector ${\cal P}'^0_n$ is given by the kernel
\begin{equation}
p'^0_n(Z_1,Z_2;W_1,W_2) =
\bigl\langle m_n(Z_1,Z_2,T), m_n(W_1,W_2,T) \bigr\rangle_n,
\end{equation}
where the pairing is done with respect to the variable $T$.
Namely, there exists a $\lambda_n \in \BB C$ such that,
for $R>0$ and $W_1,W_2 \in \BB T^+$, we have
\begin{multline}
({\cal P}'^0_n \phi_1 \otimes \phi_2)(W_1,W_2) =  \\
\frac {\lambda_n}{(2\pi^2)^2} \int_{Z_1 \in H'_R} \int_{Z_2 \in H'_R}
p'^0_n(Z_1,Z_2;W_1,W_2) \cdot
(\widetilde{\deg}_{Z_1} \phi_1)(Z_1) \cdot (\widetilde{\deg}_{Z_2} \phi_2)(Z_2)
\, \frac{dS}{\|Z_2\|} \frac{dS}{\|Z_1\|}.
\end{multline}

In particular,
\begin{multline*}
p'^0_1(Z_1,Z_2;W_1,W_2) =
\bigl\langle m_1(Z_1,Z_2,T), m_1(W_1,W_2,T) \bigr\rangle_1  \\
= \frac i{8\pi^3} \int_{\BB M}
\frac {1}{N(Z_1-T) \cdot N(Z_2-T) \cdot N(W_1-T) \cdot N(W_2-T)} \,dT^4.
\end{multline*}
\end{thm}

The integrals of this type are well known in four-dimensional
quantum field theory as Feynman integrals and admit convenient graphical
presentations. The integral that gives the kernel $p'^0_1$ of the projector
${\cal P}'^0_1$ can be viewed as the simplest nontrivial example of
a convergent Feynman integral, namely the one-loop integral for four
scalar massless particles represented by the diagram

\vskip7pt

\centerline{\includegraphics{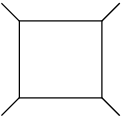}}

\noindent
There are various techniques developed by physicists for computing
Feynman integrals.
In particular, the integral corresponding to the kernel of ${\cal P}'^0_1$
can be expressed in terms of the dilogarithm function and has a remarkable
relation to the volume of the ideal tetrahedra in the three-dimensional
hyperbolic space \cite{DD}, \cite{W}.

One can generalize our integral formula for the kernel of ${\cal P}'^0_{n=1}$
to an arbitrary $n$. In order to do that one needs an explicit integral
expression for the invariant bilinear form on $\Zh^+_n$.
Then the expression for $p'^0_n$ will be given by the corresponding integral
representing this form. One can conjecture that these integrals are
also given by certain Feynman integrals.

The natural candidates for the Feynman diagrams representing $p'^0_n$ are the
$n$-loop integrals for four scalar massless particles.
However, there is an apparent difficulty that the number of such
Feynman diagrams is growing rapidly with $n$ and there are no obvious
criteria for the ``right'' choice. Fortunately, in the recent paper
``Magic identities for conformal four-point integrals''
J.~M.~Drummond, J.~Henn, V.~A.~Smirnov and E.~Sokatchev (\cite{DHSS})
show that all $n$-loop Feynman integrals for four scalar massless particles
are identical!
In particular, one can choose the so-called ladder diagrams for any $n$
consisting of a chain of $n$ boxes

\vskip7pt

\centerline{\includegraphics{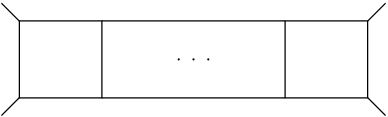}}

\noindent
The explicit expressions for these integrals $h^{(n)}(Z_1,Z_2;W_1,W_2)$
were obtained in terms of polylogarithms \cite{UD}.
Thus we end our paper with the conjecture about a relation between
the kernels of the projectors ${\cal P}'^0_n(Z_1,Z_2;W_1,W_2)$ and the $n$-loop
Feynman integrals $h^{(n)}(Z_1,Z_2;W_1,W_2)$.
We view this conjecture together with other relations that have appeared
in this paper as a beginning of a profound theory unifying
quaternionic analysis, representation theory and four-dimensional physics.

\separate

\end{document}